\documentclass[a4paper,10pt]{amsart}

\usepackage[utf8x]{inputenc}
\usepackage{amssymb}
\usepackage{amsthm}
\usepackage{amsmath}
\usepackage{amsfonts}
\usepackage{bbm}
\usepackage[T1]{fontenc}
\usepackage{graphicx}
\usepackage[dvips]{epsfig}
\usepackage[english]{babel}
\usepackage{array}
\graphicspath{{images/}}
\usepackage[a4paper]{geometry}
\usepackage{color}
\usepackage{enumerate}
\usepackage{mdwlist}
\makecompactlist{enumeratec}{enumerate}

\newcommand{\N}{\mathbb N}
\newcommand{\Z}{\mathbb Z} 

\newcommand{\R}{\mathbb R}
\newcommand{\C}{\mathbb C}

\newcommand{\la}{\lambda}
\newcommand{\Ac}{\mathcal{A}}

\newcommand{\Fc}{\mathcal{F}}
\newcommand{\cM}{\mathcal{M}}
\newcommand{\Sc}{\mathcal{S}}

\newcommand{\Rc}{\mathcal{R}}
\newcommand{\Lc}{\mathcal{L}}
\newcommand{\Jc}{\mathcal{J}}

\newcommand{\bs}{\boldsymbol}

\newcommand{\ep}{\varepsilon}
\newcommand{\La}{\Lambda}
\newcommand{\Aug}{\mathcal{A}ug}

\DeclareMathOperator{\ind}{ind}
\DeclareMathOperator{\id}{id}

\DeclareMathOperator{\mfm}{\mathfrak{m}}

\DeclareMathOperator{\Ob}{Ob}
\DeclareMathOperator{\Cth}{Cth}
\DeclareMathOperator{\Fuk}{\mathcal{F}uk}

\newtheorem{lem}{Lemma}
\newtheorem{teo}{Theorem}
\newtheorem{coro}{Corollary}
\newtheorem{prop}{Proposition}

\theoremstyle{definition}
\newtheorem{defi}{Definition}

\newtheorem{ex}{Example}
\newtheorem{rem}{Remark}

\begin{document}
\title{Product structures in Floer theory for Lagrangian cobordisms}
\author{No\'emie Legout}

\begin{abstract}
   We construct a product on the Floer complex associated to a pair of Lagrangian cobordisms. More precisely, given three exact pairwise transverse Lagrangian cobordisms in the symplectization of a contact manifold, we define a map $\mfm_2$ by a count of rigid pseudo-holomorphic disks with boundary on the cobordisms and having punctures asymptotic to intersection points and Reeb chords of the negative Legendrian ends of the cobordisms. More generally, to a $(d+1)$-tuple of exact transverse Lagrangian cobordisms we associate a map $\mfm_d$ such that the family $(\mfm_d)_{d\geq1}$ are maps satisfying the $A_\infty$ equations. Finally, we extend the Ekholm-Seidel isomorphism to an $A_\infty$-morphism, giving in particular that it is a ring isomorphism.
\end{abstract}
\maketitle

\section{Introduction}
\subsection{Background}
A \textit{contact manifold} $(Y,\xi)$ is a smooth manifold $Y$ equipped with a completely non-integrable plane field $\xi$ called a \textit{contact structure}. We consider $\xi$ cooriented, which means that there is a $1$-form $\alpha$ such that $\xi=\ker(\alpha)$ and $\alpha\wedge d\alpha\neq0$. The form $\alpha$ is called a \textit{contact form} for $(Y,\xi)$. In particular, $Y$ is an odd dimensional manifold. The \textit{Reeb vector field} $R_\alpha$ associated to $(Y,\alpha)$ is the unique vector field on $Y$ satisfying $d\alpha(R_\alpha,\cdot)=0$ and $\alpha(R_\alpha)=1$. In this article, we consider a particular type of contact manifold which is the contactization of a Liouville manifold.

A \textit{Liouville domain} $(\widehat{P},\widehat{\omega},X)$ is a compact symplectic manifold with boundary equipped with a vector field $X$ satisfying: 
 \begin{enumerate}
  \item $\Lc_X\widehat{\omega}=\widehat{\omega}$
  \item $X$ is outwards pointing along $\partial\widehat{P}$.
 \end{enumerate}
where $\Lc_X$ is the Lie derivative. Condition (1) can be rewritten $d\iota_X\widehat{\omega}=\widehat{\omega}$ because $\widehat{\omega}$ is closed, and thus this implies that it is an exact form $\widehat{\omega}=d\beta$, with $\beta=\iota_X\widehat{\omega}$. The $1$-form $\beta$ restricted to $\partial\widehat{P}$ is a contact form, and the completion of $(\widehat{P},\widehat{\omega})$ is the non compact exact symplectic manifold $(P,\omega=d\theta)$ defined by $$P=\widehat{P}\cup_{\partial\widehat{P}}\big([0,\infty)\times\partial\widehat{P}\big)$$
and $\theta$ is equal to $\beta$ on $\widehat{P}$, and to $e^\tau\beta_{|\partial\widehat{P}}$ on $[0,\infty)\times\partial\widehat{P}$, with $\tau$ the coordinate on $[0,\infty)$. The Liouville vector field $X$ on $\widehat{P}$ can be extended to the whole $(P,d\theta)$. The manifold $(P,\theta)$ is called a \textit{Liouville manifold}. The well-known symplectic manifolds $(\R^{2n},\sum_idx_i\wedge dy_i)$ and $(T^*M,-d\lambda)$, the cotangent fiber bundle of a smooth manifold $M$ equipped with the standard Liouville form, are examples of Liouville manifolds.

The \textit{contactization} of a Liouville manifold $(P,d\theta)$ is the contact manifold $(P\times\R,dz+\theta)$ where $z$ is the coordinate on the $\R$-factor. For example, the contactization of $(T^*M,-d\la)$ is the $1$-jet space $J^1(M)$.
From now, we fix a $(2n+1)$-dimensional contact manifold $(Y,\alpha)$ which is the contactization of a $2n$-dimensional Liouville manifold $(P,d\theta)$. Remark that for this special type of contact manifold, the Reeb vector field is $\partial_z$, in particular there are no closed Reeb orbits.

A \textit{Legendrian submanifold} $\La\subset Y$ is a submanifold of dimension $n$ such that $\alpha_{|T\La}=0$, which means that for all $x\in\La$, $T_x\La\subset\xi_x$. The \textit{Reeb chords} of a Legendrian submanifold $\La$ are Reeb-flow trajectories that start and end on $\La$. Compact Legendrian submanifolds in the contactization of a Liouville manifold generically have a finite number of Reeb chords. These chords correspond to vertical lines which start and end on $\La$. Let $\gamma$ be a Reeb chord of length $\ell$ which starts at a point $x^-\in\La$ and ends at $x^+\in\La$, and let us denote $\varphi^R_t$ the Reeb flow. If $d_{x^-}\varphi^R_\ell(T_{x^-}\La)$ and $T_{x^+}\La$ intersect transversely, we say that the Reeb chord $\gamma$ is \textit{non-degenerate}, and then $\La$ is called \textit{chord generic} if all its Reeb chords are non-degenerate. From now, we will only consider compact chord generic Legendrian submanifolds, and we denote by $\Rc(\La)$ the set of Reeb chords of $\La$. If $\La_1,\La_2,\dots,\La_d$ are $d$ Legendrian submanifolds of $Y$, we can consider the union $\La=\La_1\cup\dots\cup\La_d$. Reeb chords of $\La$ from $\La_i$ to itself are called \textit{pure} Reeb chords while those from $\La_i$ to $\La_j$ with $i\neq j$ are called \textit{mixed} Reeb chords. We denote by $\Rc(\La_i,\La_j)$ the set of Reeb chords from $\La_i$ to $\La_j$.
 
The \textit{Lagrangian projection} of a Legendrian submanifold $\La\subset Y=P\times\R$ is the image of $\La$ under the projection $\Pi_P\colon P\times\R\to P$. Reeb chords of $\La$ are then in bijection with intersection points of $\Pi_P(\La)$. In the particular case where the contact manifold is the $1$-jet space of a manifold $M$ (i.e. $Y=J^1(M)=T^*M\times\R$), the \textit{front projection} of $\La$ is the image of $\La$ under $\Pi_F\colon J^1(M)\to M\times\R$. In this case, Reeb chords are in bijection with vertical segments in $M\times\R$ beginning and ending respectively on points $c^-,c^+\in\Pi_F(\La)$, and such that the tangent spaces $T_{c^-}\Pi_F(\La)$ and $T_{c^+}\Pi_F(\La)$ are equal.

One natural question when studying Legendrian submanifolds is to understand whether two Legendrian submanifolds $\La_0, \La_1\subset Y$ are Legendrian isotopic or not (i.e. is there a smooth function $F:[0,1]\times\La\to Y$ such that $\La$ is a $n$-dimensional manifold and $F(t,\La)$ is a Legendrian submanifold of $Y$ for all $t\in[0,1]$, with $F(0,\La)=\La_0$ and $F(1,\La)=\La_1$?). A lot of work has been achieved in order to answer this question of classification under Legendrian isotopy of Legendrian submanifolds. There exists a lot of Legendrian isotopy invariants, among which the first were the classical ones, namely the smooth isotopy type, the Thurston-Bennequin invariant and the rotation class (see for example \cite{ej,EES3}). The development then of non-classical invariants gave new directions in order to better understand the Legendrians. One of the first non-classical invariants is a relative version of contact homology \cite{EGH} called the Legendrian contact homology. It was defined by Eliashberg in \cite{E} using pseudo-holomorphic curves techniques. Independently, it was defined combinatorially by Chekanov for Legendrian links in $(\R^3,dz-ydx)$ in \cite{Che}, and this combinatorial description was generalized in higher dimension by Ekholm, Etnyre and Sullivan \cite{EES1,EES2}. These two definitions were then shown to compute the same invariant, by Etnyre, Ng and Sabloff \cite{ENS} in dimension $3$, and by Dimitroglou-Rizell \cite{DR} in all dimension. This is a very powerful Legendrian isotopy invariant which gave rise to numerous other invariants, as for example the linearized and ``multi-linearized'' versions of Legendrian contact homology, using augmentations of the differential graded algebra introduced by Chekanov. Then there are higher algebraic structures on linearized Legendrian contact cohomology that are Legendrian isotopy invariants, as a product structure and an $A_\infty$-algebra structure (see \cite{CKESW}), and more generally, there are $A_\infty$-categories $\Aug_-(\La)$ and $\Aug_+(\La)$, called the augmentation categories of a Legendrian submanifold (see \cite{BCh}, \cite{NRSSZ}, and Subsection \ref{linearisation} for $\Aug_-(\La)$). In parallel to these invariants defined by pseudo-holomorphic curves counts, other types of Legendrian isotopy invariants have been defined, by generating functions techniques. We will not go through these invariants in this article, nevertheless, even if the definition of this two types (pseudo-holomorphic curves vs generating functions) of invariants are constructed using completely different techniques, they are closely related. Indeed, the existence of a (linear at infinity) generating family for a Legendrian knot $\La$ in $\R^3$ implies the existence of an augmentation such that the linearized contact homology of $\La$ is isomorphic to the generating family homology of $\La$ (see \cite{FR}). In higher dimension, the relation is not so clear. However, there are parallel results in the Legendrian contact homology side and the generating family homology side, as for example a duality exact sequence (\cite{EESa} \cite{ST}), and also results relying the Legendrian invariants and the topology of a Lagrangian filling of the Legendrian (\cite{Ch1} \cite{Gol} \cite{E1} \cite{DR} \cite{ST}). To continue along this path, we could imagine to define the algebraic structures appearing in this paper in the generating family setting, that is to say a generating family Floer complex, and a generating family product on Floer complexes, that could be related through the cobordism maps to the product structure on generating family homology defined by Myer (\cite{ZM}).
\vspace{2mm}

In this article, we will be interested in the relation of exact Lagrangian cobordisms between Legendrian submanifolds, introduced by Chantraine in \cite{Ch1}. These objects live in the symplectization of the contact manifold $(Y,\alpha)$, which is the symplectic manifold $(\R\times Y,d(e^t\alpha))$ where $t$ is the coordinate on $\R$.
\begin{defi}
 An \textit{exact Lagrangian cobordism} from $\La^-$ to $\La^+$, denoted $\La^-\prec_\Sigma\La^+$, is a properly embedded submanifold $\Sigma\subset\R\times Y$ satisfying the following:
 \begin{enumerate}
  \item there exists a constant $T>0$ such that:
  \begin{itemize}
    \item[$-$] $\Sigma\cap(-\infty,-T)\times Y=(-\infty,-T)\times\La^-$,
    \item[$-$] $\Sigma\cap(T,\infty)\times Y=(T,\infty)\times\La^+$,
    \item[$-$] $\Sigma\cap[-T,T]\times Y$ is compact.
    \end{itemize}
 \item there exists a smooth function $f\colon\Sigma\to\R$ such that:
    \begin{enumerate}
    \item $e^t\alpha_{|T\Sigma}=df$,
    \item $f_{|(-\infty,-T)\times\La^-}$ is constant,
    \item $f_{|(T,\infty)\times\La^+}$ is constant.
    \end{enumerate}
 \end{enumerate}
\end{defi}
\begin{rem} 
 Condition (a) above says by definition that $\Sigma$ is an exact Lagrangian submanifold of $(\R\times Y,d(e^t\alpha))$. Moreover, using the fact that $\Sigma$ is a cylinder over $\La^-$ in the negative end and a cylinder over $\La^+$ in the positive end, condition (a) implies that $f$ is constant on each connected component of the negative and the positive ends of $\Sigma$. Conditions (b) and (c) imply that $f$ is in fact globally constant on each end (the constant on the positive end is not necessarily the same as the constant on the negative end). Thus, if $\La^\pm$ are connected, conditions (b) and (c) are automatically satisfied.
\end{rem}
We denote by $\overline{\Sigma}:=[-T,T]\times\Sigma$ the compact part of the cobordism and the boundary components $\partial_-\overline{\Sigma}=\{-T\}\times\Lambda^-$ and $\partial_+\overline{\Sigma}:=\{T\}\times\Lambda^+$. In the case where $\Sigma$ is diffeomorphic to a cylinder, we call it a \textit{Lagrangian concordance} from $\La^-$ to $\La^+$ and  denote it simply $\La^-\prec\La^+$, and when $\La^-=\emptyset$, $\Sigma$ is called an \textit{exact Lagrangian filling} of $\La^+$. A Legendrian isotopy between two Legendrian submanifolds $\La_1$ and $\La_2$ induces Lagrangian concordances $\La_1\prec\La_2$ and $\La_2\prec\La_1$ \cite{EG,Ch1}, however, if such concordances exist it is not known if it implies that $\La_1$ and $\La_2$ are Legendrian isotopic. In general, as evoked above, some Legendrian isotopy invariants give obstructions to the existence of Lagrangian cobordisms (see for example \cite{Ch1,E1,ST,CNS,Yu}, which is absolutely not an exhaustive list). In the same vein, Chantraine, Dimitroglou-Rizell, Ghiggini and Golovko (\cite{CDGG2}) have defined a Floer-type complex associated to a pair of Lagrangian cobordisms, the Cthulhu complex, in order to understand better the topology of a Lagrangian cobordism between two given Legendrians. The goal of this article is to provide a richer algebraic structure associated to Lagrangian cobordisms.

\subsection{Results}
Let $(Y,\alpha)$ be the contactization of a Liouville manifold $(P^{2n},d\theta)$ and consider four Legendrian submanifolds $\La_1^-,\La_1^+,\La_2^-,\La_2^+\subset Y$ such that the Chekanov-Eliashberg algebras (Legendrian contact homology algebras) $\Ac(\La_1^-)$ and $\Ac(\La_2^-)$ admit augmentations $\ep_1^-$ and $\ep_2^-$ respectively. Throughout the paper, the coefficient field is $\Z_2$, i.e. the algebras and vector spaces we consider are all over $\Z_2$. Assume there exist two exact transverse Lagrangian cobordisms $\La_1^-\prec_{\Sigma_1}\La_1^+$ and $\La_2^-\prec_{\Sigma_2}\La_2^+$. The Cthulhu complex $(\Cth(\Sigma_1,\Sigma_2),\mathfrak{d}_{\ep_1^-,\ep_2^-})$ associated to the pair $(\Sigma_1,\Sigma_2)$ is generated by Reeb chords from $\La_2^+$ to $\La_1^+$, intersection points in $\Sigma_1\cap\Sigma_2$, and by Reeb chords from $\La_2^-$ to $\La_1^-$. Given $\ep_1^+$ and $\ep_2^+$ augmentations of $\Ac(\La_1^+)$ and $\Ac(\La_2^+)$ respectively induced by $\ep_1^-$ and $\ep_2^-$ (see Section \ref{morphisme_induit}), the differential of the Cthulhu complex is a linear map defined by a count of rigid pseudo-holomorphic curves with boundary on $\Sigma_1$ and $\Sigma_2$ (see Section \ref{cthulhu}). This complex admits a quotient complex $CF_{-\infty}(\Sigma_1,\Sigma_2)$, generated only by intersection points and Reeb chords from $\La_2^-$ to $\La_1^-$, called the \textit{Floer complex} of the pair $(\Sigma_1,\Sigma_2)$. The main result of this article is the following: 

\begin{teo}\label{teo_prod}
Let $\Sigma_1, \Sigma_2$ and $\Sigma_3$ be three pairwise transverse exact Lagrangian cobordisms from $\La_i^-$ to $\La_i^+$, for $i=1,2,3$, where $\La_i^\pm$ are Legendrian submanifolds of $P\times\R$ such that the Chekanov-Eliashberg algebras $\Ac(\La_i^-)$ admit augmentations. Then,
\begin{enumerate}
	\item for any choice of augmentation $\ep_i^-$ of $\Ac(\La_i^-)$, for $i=1,2,3$, there exists a map:
 \begin{alignat*}{1}
  \mfm_2\colon CF_{-\infty}(\Sigma_2,\Sigma_3)\otimes CF_{-\infty}(\Sigma_1,\Sigma_2)\to CF_{-\infty}(\Sigma_1,\Sigma_3)  
 \end{alignat*}
which satisfies the Leibniz rule $\partial_{-\infty}\circ \mfm_2(-,-)+\mfm_2(\partial_{-\infty},-)+\mfm_2(-,\partial_{-\infty})=0$.
	\item in the case where $\La_1^-=\emptyset$, and $\Sigma_2$ and $\Sigma_3$ are small Hamiltonian perturbations of $\Sigma_1$ such that the pairs $(\Sigma_1,\Sigma_2)$, $(\Sigma_2,\Sigma_3)$ and $(\Sigma_1,\Sigma_3)$ are directed (see Section \ref{perturbations}), the product $\mfm_2$ is equal to the cup product on $H^*(\Sigma_1,\La_1^+)$ after appropriate identifications.
\end{enumerate}
\end{teo}

Now, the Cthulhu homology is invariant by a certain type of Hamiltonian isotopy which permits to displace the Lagrangian cobordisms. This implies the acyclicity of the complex. But the Cthulhu complex is in fact the cone of a map
\begin{alignat*}{1}
	\Fc_{21}^1\colon CF_{-\infty}(\Sigma_1,\Sigma_2)\to C(\La_1^+,\La_2^+)
\end{alignat*}
from the Floer complex to a complex generated by Reeb chords from $\La_2^+$ to $\La_1^+$ (the bilinearized Legendrian contact cohomology complex of $\La_1^+\cup\La_2^+$ restricted to chords from $\La_2^+$ to $\La_1^+$). The acyclicity implies that this map is a quasi-isomorphism. When $\Sigma_1$ is a Lagrangian filling of $\La_1^+$ and $\ep_1^+$ is the augmentation induced by this filling, take $\Sigma_2$ a small perturbation of $\Sigma_1$ such that the pair $(\Sigma_1,\Sigma_2)$ is directed, then the quasi-isomorphism $\Fc^1$ recovers Ekholm-Seidel isomorphism (\cite{E1,DR}). We will show that the map induced by $\Fc^1$ in homology preserves the product structures, that is to say, the product $\mfm_2$ on Floer complexes is mapped to the product $\mu^2_{\ep_{3,2,1}^+}$ of the augmentation category $\Aug_-(\La_1^+\cup\La_2^+\cup\La_3^+)$, where $\ep^+_{3,2,1}$ is the diagonal augmentation on the algebra $\Ac(\La_1^+\cup\La_2^+\cup\La_3^+)$ induced by $\ep_1^+,\ep_2^+$ and $\ep_3^+$ (see Section \ref{linearisation}). More precisely, we have:
\begin{teo}\label{teofoncteur}
 Let $\Sigma_1, \Sigma_2$ and $\Sigma_3$ be three transverse exact Lagrangian cobordisms from $\La_i^-$ to $\La_i^+$, such that the Chekanov-Eliashberg algebras $\Ac(\La_i^-)$ admit augmentations. Then, for any choice of augmentation $\ep_i^-$ of $\Ac(\La_i^-)$ we have:
 \begin{alignat*}{1}
  \big[\mu^2_{\ep_{3,2,1}^+}(\Fc_{32}^1,\Fc_{21}^1)]=[\Fc_{31}^1\circ\mfm_2\big]
 \end{alignat*}
\end{teo}
In the same setting as part \textit{(2)} of Theorem \ref{teo_prod}, Theorem \ref{teofoncteur} implies that the Ekholm-Seidel isomorphism is a ring morphism:
\begin{coro}\label{isoES}
	Let $\Sigma$ be an exact Lagrangian filling of a Legendrian submanifold $\La\subset Y$, and denote $\ep_{\Sigma}$ the augmentation of $\Ac(\La)$ induced by $\Sigma$. Then the Ekholm-Seidel isomorphism
	\begin{alignat*}{1}
	 LCH_{\ep_{\Sigma}}^*(\La)\simeq H^{*+1}(\Sigma,\Lambda)
	 \end{alignat*}
	 is an isomorphism of non-unital rings.
\end{coro}
A related result appears in a paper of Ekholm and Lekili (see \cite[Theorem 53]{EL} and Remark \ref{remEL}).
The product $\mfm_2$ is in fact part of an $A_\infty$-structure defined in Section \ref{infini}. There, we define maps $\{\mfm_k\}_{1\leq k\leq d}$ for any $(d+1)$-tuple of pairwise transverse Lagrangian cobordisms such that the $\mfm_1$ map is the differential $\partial_{-\infty}$. We refer the reader to Section \ref{infini} for the precise domains and codomains of the maps $\mfm_i$. These maps are defined on tensor products of Floer complexes between \textit{transverse} Lagrangians, and so in particular we do not have an $A_\infty$-algebra structure. Then we prove the following:
\begin{teo}\label{teo3}
 The maps $\{\mfm_k\}_{1\leq k\leq d}$ satisfy the $A_\infty$ equations, i.e. for all $1\leq k\leq d$:
 \begin{alignat*}{1}
   \sum\limits_{\substack{1\leq j\leq k\\ 0\leq n\leq k-j}}\mfm_{k-j+1}(\id^{\otimes k-j-n}\otimes \mfm_j\otimes\id^{\otimes n})=0
  \end{alignat*}
\end{teo}

A direct corollary of Theorem \ref{teo3} is that the product $\mfm_2$ is associative in homology. Similarly, the map $\Fc^1$ extends in a family of maps $\Fc=\{\Fc^k\}_{0\leq k\leq d}$ defined for any $(d+1)$-tuple of pairwise transverse Lagrangian cobordisms and we have:
\begin{teo}\label {ffon}
	The maps $\Fc=\{\Fc^k\}_{0\leq k\leq d}$ satisfy the $A_\infty$-functor equations. 
\end{teo}

Of course we would like to be able to define the maps $\mfm_d$ and $\Fc_d$ for families of not necessarily transverse cobordisms. We conjecture that given a pair of Legendrian submanifolds $\La_-,\La_+\subset Y$, there are unital $A_\infty$-categories of cobordisms from $\La^-$ to $\La^+$ denoted $\Fuk_-(\La^+,\La^-)$ and $\Fuk_+(\La^+,\La^-)$ which can be defined by localisation (as used in \cite{GPS} to define the wrapped Fukaya category of a Liouville sector). Moreover, we would obtain cohomologically full and faithful unital $A_\infty$-functors
	\begin{alignat*}{1}
	\Fc_\pm:\Fuk_\pm(\La^+,\La^-)\to \Aug_\pm(\La^+)
	\end{alignat*}
adding a unit to $\Aug_-(\La^+)$ to make it unital, as explained in \cite[Remark 26]{EL}. The definition of these categories and functors will be done in a forthcoming paper. Let us remark that the categories of fillings $\Fuk_\pm(\La^+,\emptyset)$ could probably also be defined with another algebraic approach in the same spirit as Ekholm-Lekili \cite[Section 3]{EL}, using coefficients in chains in the based loop space of $\La^+$, but we will not develop this in this article. Also, in \cite[Theorem 1.6]{Yu} it is proven that the functor $\Aug_+(\La^-)\to\Aug_+(\La^+)$ induced by a cobordism on augmentation categories is cohomologically faithful. The faithfulness of $\Fc_+$ would in particular recover this, but then to get fullness it is needed to take into account intersection point generators.
\vspace{3mm}

The paper is organized as follows. In Section \ref{moduli} we set up the definition and notations of all types of moduli spaces that are involved in the definition of all maps in the rest of the paper. In Sections \ref{LCH} and \ref{CTH}, we review the definitions of Legendrian contact homology and Cthulhu homology. In Section \ref{PROD}, we construct the product structure on the Floer complexes and prove Theorem \ref{teo_prod} and Theorem \ref{teofoncteur}. We give a very basic example of computation of the product in Section \ref{example}, and finally in Section \ref{infini} we define the $A_\infty$-structure on Floer complexes and prove Theorem \ref{teo3} and Theorem \ref{ffon}.

\subsection*{Acknowledgements}
This work is part of my PhD thesis which I did at the Université de Nantes under the supervision of Frédéric Bourgeois and Baptiste Chantraine, who I warmly thank for their guidance, help and support. I would also like to thank the referees of my PhD thesis Jean-François Barraud and Lenny Ng, for many helpful comments and remarks, and Paolo Ghiggini and François Laudenbach for helpful discussions. Finally, I thank the anonymous referee for the corrections and comments. Part of this paper has been written at the CIRGET which I thank for its hospitality.

\section{Moduli spaces}\label{moduli}
In this section we describe the different types of moduli spaces of pseudo-holomorphic curves which will be necessary to define the Legendrian contact homology complex, the Cthulhu complex and the product structure. The first three subsections contain useful material from \cite{S} and \cite{CDGG2}, in order to define the moduli spaces.

\subsection{Deligne-Mumford space}\label{deligne}

Let us denote $$\Rc^{d+1}=\{(y_0,\dots,y_d)\,|\,y_i\in S^1,\,y_i\in(y_{i-1},y_{i+1})\}/Aut(D^2)$$ the space of $(d+1)$-tuples of points cyclically ordered on the boundary of the disk $D^2$, where $y_{-1}:=y_d$ and $y_{d+1}:=y_0$, and denote by $\Sc^{d+1}$ the universal curve: 
\begin{alignat*}{1}
 \Sc^{d+1}=\{(z,y_0,\dots,y_d)\,/\,z\in D^2,\,y_i\in S^1 \mbox{ and } y_i\in(y_{i-1},y_{i+1})\}/Aut(D^2)
\end{alignat*}
The projection $\pi\colon\Sc^{d+1}\to\Rc^{d+1}$ given by $\pi(z,y_0,\dots,y_d)=(y_0,\dots,y_d)$ is a fibration with fiber a disk. For all $r\in\Rc^{d+1}$ we denote $\widehat{S}_r=\pi^{-1}(r)$ and $S_r=\widehat{S}_r\backslash\{y_0,\dots,y_d\}$.

Given $r\in\Rc^{d+1}$, to each marked point $y_i$ of $\widehat{S}_r$, $i\geq1$, one can associate a neighborhood $V_i\subset\widehat{S}_r$ and a biholomorphism $\ep_i\colon(-\infty,0)\times[0,1]\to V_i\backslash\{y_i\}$. For the puncture $y_0$, we choose a neighborhood $V_0$ and a biholomorphism $\ep_0\colon(0,+\infty)\times[0,1]\to V_0\backslash\{y_0\}$. These biholomorphisms are called \textit{strip-like ends}. A \textit{universal choice of strip-like ends} for $\Rc^{d+1}$ corresponds to maps
\begin{alignat*}{1}
 \ep_0^{d+1}\colon\Rc^{d+1}\times(0,+\infty)\times[0,1]\to\Sc^{d+1}
\end{alignat*}
and
\begin{alignat*}{1}
 \ep_i^{d+1}\colon\Rc^{d+1}\times(-\infty,0)\times[0,1]\to\Sc^{d+1}
\end{alignat*}
for $1\leq i\leq d$, such that for all $r\in\Rc^{d+1}$, $\ep_i^{d+1}(r,\cdot,\cdot)$ is a choice of strip-like ends for $y_i\in\widehat{S}_r$.

The space $\Rc^{d+1}$, for $d\geq2$, admits a compactification which can be described in terms of trees. In fact, we have $\overline{\Rc}^{d+1}=\sqcup_T\Rc^T$ which is a disjoint union over all stable planar rooted trees $T$ with $d$ leaves, and with $\Rc^T=\sqcup\Rc^{|v_i|}$ where the union is over all interior vertices (vertices which are neither leaves nor the root) $v_i$ of $T$. Here, $|v_i|$ denotes the degree of the vertex $v_i$, and recall that a tree is called \textit{stable} if each interior vertex has degree at least $3$. The space $\Rc^{d+1}$ corresponds to $\Rc^{T_{d+1}}$ where $T_{d+1}$ is the planar rooted tree with $d$ leaves and one interior vertex. An \textit{interior edge} of a planar rooted tree $T$ is an edge between two interior vertices. We denote by $Ed^{int}(T)$ the set of interior edges of $T$. Given $T$ and $T'$ two stable planar rooted trees with $d$ leaves, if $T'$ can be obtained from $T$ by removing one or several interior edges (i.e. contracting an edge until the two corresponding vertices are identified), it gives rise to a gluing map:
\begin{alignat*}{1}
 \gamma^{T,T'}\colon\Rc^{T}\times(-1,0]^{Ed^{int}(T)}\to\Rc^{T'}
\end{alignat*}
If $e$ is an interior edge from the vertices $v^-$ to $v^+$ to remove of $T$ to obtain $T'$, this gluing map consists in gluing the two disks $S_{r_{v^-}}$ and $S_{r_{v^+}}$ along $e$ with a certain gluing parameter. Let us denote $\ep_-$ and $\ep_+$ the strip-like ends of $r_{v^-}$ and $r_{v^+}$ for the marked points connected by $e$. Given a real $l_e\in(0,\infty)$, the gluing operation is given by the connected sum
\begin{alignat*}{1}
 S_{r_{v^-}}\backslash\ep_-((-\infty,l_e)\times[0,1])\bigcup S_{r_{v^+}}\backslash\ep_+((l_e,\infty)\times[0,1])/\sim
\end{alignat*}
where we identify $\ep_-(l_e-s,t)\sim\ep_+(s,t)$. The map $\gamma^{T,T'}$ glues each interior edge of $T$ using the parameter $\rho_e=-e^{-\pi l_e}\in(-1,0]$ instead of $l_e$. If $\rho=0$, the edge is not modified (see \cite{S}).

Now suppose that $S\in\Sc^{d+1}$ is obtained from $S_{r_1},S_{r_2},\dots S_{r_k}$ by gluing, then $S$ admits a \textit{thin-thick decomposition}. The thin part $S^{thin}$ corresponds to strip-like ends of $S$ and to strips of length $l_e$ coming from the identification of strip-like ends in the gluing of two disks $S_{r_{v_i}}$ and $S_{r_{v_j}}$ along an edge $e$. The thick part is then $S\backslash S^{thin}$. If $r\in\Rc^{d+1}$ is in the image of $\gamma^{T,T_{d+1}}$, then it admits two sets of strip-like ends: one coming from the universal choice on $\Rc^{d+1}$ and the other one coming from the universal choice on $\Rc^{|v_i|}$ for all vertices $v_i$ of $T$ and the gluing operation. A universal choice of strip-like ends on $\Rc^{d+1}$ is said \textit{consistent} if there exists a neighborhood $U\subset\overline{\Rc}^{d+1}$ of $\partial\overline{\Rc}^{d+1}$ such that the two choices of strip-like ends coincide on $U\cap\Rc^{d+1}$.
\begin{teo}\cite[Lemma 9.3]{S}
 Consistent universal choices of strip-like ends exist.
\end{teo}

\begin{rem}
 In the cases $d=0,1$, a punctured disk in $\Sc^1$ is biholomorphic to a half-plane and a punctured disk in $\Sc^{2}$ is biholomorphic to the strip $Z=\R\times[0,1]$ with standard coordinates $(s,t)$.
\end{rem}

\subsection{Lagrangian labels}\label{etiquette}
The holomorphic disks we will consider are holomorphic maps from a disk with some marked points removed to the manifold $\R\times Y$, with boundary on Lagrangian submanifolds of $\R\times Y$. The corresponding Lagrangian submanifolds are called a \textit{Lagrangian label} for the disk and is defined as follows. 

The boundary of $S_r$, for $r\in\Rc^{d+1}$ is subdivided into $d+1$ components. We denote by $\partial_iS_r$ for $1\leq i\leq d+1$ the part of the boundary between the marked points $y_{i-1}$ and $y_i$. A Lagrangian label for $S_r$ is a choice of Lagrangian submanifolds $L_i\subset\R\times Y$ for each component $\partial_iS_r$ of the boundary of $S_r$.
If $L_i$ is the Lagrangian submanifold associated to $\partial_iS_r$, we will denote by $\underline{L}=(L_1,\dots,L_{d+1})$ the Lagrangian label for $S_r$.
A natural compatibility condition for Lagrangian labels is clearly necessary in order to apply the gluing maps $\gamma^{T,T'}$.


\subsection{Almost complex structure}
In this subsection we recall the different types of almost complex structures that will be useful in order to achieve transversality for moduli spaces. Recall that on a symplectic manifold $(X,\omega)$, an almost complex structure is a map $J\colon TX\to TX$ such that $J^2=-\id$. We say that $J$ is \textit{compatible} with $\omega$ (or \textit{$\omega$-compatible}) if:
 \begin{enumerate}
  \item $\omega(v,Jv)>0$ for all $v\in TX$ such that $v\neq0$,
  \item $\omega(Ju,Jv)=\omega(u,v)$ for all $x\in X$ and $u,v\in T_xX$.
 \end{enumerate}
 
\subsubsection{Cylindrical almost complex structure}
Let us go back to the case where the symplectic manifold is the symplectization of a contact manifold $(Y,\alpha)$. An almost complex structure $J$ on $(\R\times Y,d(e^t\alpha))$ is \textit{cylindrical} if:
 \begin{itemize}
  \item[$\bullet$] $J$ is $d(e^t\alpha)$-compatible,
  \item[$\bullet$] $J$ is invariant under $\R$-action by translation on $\R\times Y$,
  \item[$\bullet$] $J(\partial_t)=R_\alpha$,
  \item[$\bullet$] $J$ preserves the contact structure, i.e. $J(\xi)=\xi$.  
 \end{itemize}
Following notations of \cite{CDGG2}, we denote by $\Jc^{cyl}(\R\times Y)$ the set of cylindrical almost complex structures on $\R\times Y$.

In our setting, the contact manifold is the contactization of a Liouville manifold, $Y=P\times\R$, and recall that a Liouville manifold $P$ can be viewed as the completion of a Liouville domain $(\widehat{P},d\beta)$. An almost complex structure $J_P$ on $P$ is \textit{admissible} if it is cylindrical on $P\backslash\widehat{P}$ outside of a compact subset $K\subset P\backslash\widehat{P}$. We denote by $\Jc^{adm}(P)$ the set of admissible almost complex structures on $P$. Now, if $J_P\in\Jc^{adm}(P)$ and $\pi_P\colon\R\times(P\times\R)\to P$ is the projection on $P$, then there exists a unique cylindrical almost complex structure $\widetilde{J}_P$ on $\R\times(P\times\R)$ such that $\pi_P$ is holomorphic, that is to say $d\pi_P\circ\widetilde{J}_P=J_P\circ d\pi_P$. Such an almost complex structure is called the \textit{cylindrical lift} of $J_P$ and we denote by $\Jc^{cyl}_\pi(\R\times Y)$ the set of cylindrical almost complex structures on $\R\times Y$ which are cylindrical lifts of admissible almost complex structures on $P$.

Let $J^-,J^+\in\Jc^{cyl}(\R\times Y)$ such that $J^-$ and $J^+$ coincide outside of a cylinder $\R\times K$ where $K\subset Y$ is compact. For all $T>0$ we consider an almost complex structure $J$ on $\R\times Y$ equals to $J^-$ on $(-\infty,-T)\times Y$, $J^+$ on $(T,\infty)\times Y$ and equals to the cylindrical lift of an admissible complex structure on $P$ in $[-T,T]\times(Y\backslash K)$. The reason for considering such almost complex structures is that transversality holds generically for moduli spaces of Legendrian contact homology with a cylindrical almost complex structure (see Section \ref{LCH}), and that cylindrical lifts of admissible almost complex structures on $P$ are useful to prevent pseudo-holomorphic curves to escape at infinity (the projection on $P\times\R$ must be compact).

We denote by $\Jc^{adm}_{J^-,J⁺,T}(\R\times Y)$ the set of almost complex structures on $\R\times Y$ described above, and $\Jc^{adm}(\R\times Y)=\bigcup\limits_{J^-,J^+,T}\Jc^{adm}_{J^-,J^+,T}(\R\times Y)$. 

\subsubsection{Domain dependent almost complex structure}\label{str_dep}
Considering domain dependent almost complex structures is a way to achieve transversality for moduli spaces of pseudo-holomorphic curves. A domain dependent almost complex structure on $\R\times Y$ is the data, for each $r\in\Rc^{d+1}$, of an almost complex structure parametrized by $S_r$, that is to say a map in $C^{\infty}(S_r,\Jc^{adm}(\R\times Y))$. Then, we need some special behavior of the almost complex structure in strip-like ends in order to get some compatibility with the gluing map.

Fix a $r\in\Rc^{d+1}$, and let $L_1,\dots,L_{d+1}$ be transverse exact Lagrangian cobordisms in $\R\times Y$ such that $\underline{L}=(L_1,\dots,L_{d+1})$ is a choice of Lagrangian label for $S_r$. Let $T>0$ such that all the $L_i$'s are cylindrical out of $L_i\cap[-T,T]\times Y$, and take $J^\pm\in\Jc^{cyl}(\R\times Y)$.

For each pair $(L_i,L_{i+1})$, we consider a path $J_t^{L_i,L_{i+1}}$ for $t\in[0,1]$ of almost complex structures in $\Jc_{J^-,J^+,T}^{adm}(\R\times Y)$, such that it is constant near $t=0$ and $t=1$. The type of domain dependent almost complex structures we consider are maps $$J_{r,\underline{L}}:S_r\to\Jc_{J^-,J^+,T}^{adm}(\R\times Y)$$ such that $J_{r,\underline{L}}(\ep_i(s,t))=J_t^{L_i,L_{i+1}}$, where $\ep_i$ is a choice of strip-like ends for $S_r$.

Now, consider a universal choice of strip-like ends. A \textit{universal choice of domain dependent almost complex structures} is the data, for all $r\in\Rc^{d+1}$ and Lagrangian label $\underline{L}=(L_1,\dots,L_{d+1})$, of maps $J_{r,\underline{L}}$ as above that fit into a smooth map
$$\Jc_{d,\underline{L}}:\Sc^{d+1}\to\Jc_{J^-,J^+,T}^{adm}(\R\times Y)$$
defined by $\Jc_{d,\underline{L}}(z)=J_{r,\underline{L}}(z)$ if $z\in S_r$. Moreover, $\Jc_{d,\underline{L}}$ must satisfy $\Jc_{d,\underline{L}}(\ep_i^{d+1}(r,s,t))=J_t^{L_i,L_{i+1}}$ where $\ep_i^{d+1}$ is part of the universal choice of strip-like ends.

Again, if $S\in\Sc^{d+1}$ is obtained from $S_{r_1},S_{r_2},\dots S_{r_k}$ by gluing, we need compatibility conditions between the almost complex structure induced by the universal choice and the one induced by the gluing map. The two choices of almost complex structures are said \textit{consistent} if there exists a neighborhood $U\subset\Rc^{d+1}$ of $\partial\overline{\Rc}^{d+1}$ such that the choice of strip-like ends is consistent, the choices of almost complex structures coincide on the thin parts for each $r\in U$, and for every sequence $\{r^n\}_{n\in\N}$ in $\Rc^{d+1}$ converging to a point $r\in\partial\overline{\Rc}^{d+1}$, the almost complex structures on the thick parts must converge to the almost complex structure on the thick part of $S_r$.

The latter condition on thick part is analogous to the condition on thin parts, the difference is that we ask for convergence of almost complex structures instead of equality because the almost complex structure on thick parts is not fixed, whereas it is on thin parts. Indeed, a universal choice of almost complex structures depends on fixed paths $J_t^{L_i,L_{i+1}}$  for each pair of Lagrangian submanifolds.
\begin{teo}\cite[Lemma 9.5]{S}
 Consistent choices of almost complex structures exist.
\end{teo}

\subsection{Moduli spaces of holomorphic curves}\label{esp_de_mod}
We are now ready to define the moduli spaces we will use in the next sections.

\subsubsection{General definition}\label{gen-def}
Let $\underline{\Sigma}=(\Sigma_1,\dots,\Sigma_{d+1})$ be a choice of Lagrangian label such that for all $1\leq i\leq d+1$, $\Sigma_i$ is an exact Lagrangian cobordism from $\La_i^-$ to $\La_i^+$. We assume that the cobordisms are pairwise transverse. We consider then a set $A(\underline{\Sigma})$ of \textit{asymptotics} consisting of intersection points in $\Sigma_i\cap\Sigma_j$ for all $1\leq i\neq j\leq d+1$, Reeb chords from $\La_i^+$ to $\La_j^+$, and Reeb chords from $\La_i^-$ to $\La_j^-$ for all $1\leq i, j\leq d+1$. Let $J$ be an almost complex structure on $\R\times Y$ (we will explain later the properties needed to achieve transversality in each case), and $j$ the standard almost complex structure on the disk $D^2\subset \C$, which induces an almost complex structure on each $S_r$, $r\in\Rc^{d+1}$. For $r\in\Rc^{d+1}$ and $x_0,\dots,x_d$ in $A(\underline{\Sigma})$, we define the moduli space $\cM^r_{\underline{\Sigma},J}(x_0;x_1,\dots,x_d)$ as the set of smooth maps:
\begin{alignat*}{1}
 u\colon(S_r,j)\to(\R\times Y,J) 
\end{alignat*}
satisfying:
\begin{enumerate}
 \item $du(z)\circ j=J(z)\circ du(z)$, for all $z\in S_r\backslash\partial S_r$,
 \item $u(\partial_iS_r)\subset\Sigma_i$,
 \item if $x_0$ is an intersection point then $\lim\limits_{z\to y_0}u(z)=x_0$ and $x_0$ is required to be a jump from $\Sigma_{d+1}$ to $\Sigma_1$ when traversing the boundary counterclockwise,
 \item if $x_i$, $1\leq i\leq d$, is an intersection point then $\lim\limits_{z\to y_i}u(z)=x_i$,
 \item if $x_0$ is a Reeb chord with a parametrization $\gamma_0\colon[0,1]\to x_0$, then every $z\in S_r$ sufficiently close to $y_0$ is in $\ep_0((0,+\infty)\times[0,1])$ and we have either
 \begin{itemize}
 	\item[$\bullet$] $\lim\limits_{s\to+\infty}u(\ep_0(s,t))=(+\infty,\gamma_0(t))$, and in this case we say that $u$ has a \textit{positive asymptotic} to $x_0$ at $y_0$, or
 	\item[$\bullet$] $\lim\limits_{s\to+\infty}u(\ep_0(s,t))=(-\infty,\gamma_0(1-t))$ and we say that $u$ has a \textit{negative asymptotic} to $x_0$ at $y_0$.
 \end{itemize}
 \item if $x_i$ for $i>0$ is a Reeb chord with parametrization $\gamma_i\colon[0,1]\to x_i$, then either
 \begin{itemize}
  \item[$\bullet$] $\lim\limits_{s\to-\infty}u(\ep_i(s,t))=(-\infty,\gamma_i(t))$ and $u$ has a negative asymptotic to $x_i$ at $y_i$, or
  \item[$\bullet$] $\lim\limits_{s\to-\infty}u(\ep_i(s,t))=(+\infty,\gamma_i(1-t))$ and $u$ has a positive asymptotic to $x_i$ at $y_i$.
 \end{itemize}
\end{enumerate}
Then we denote 
\begin{alignat*}{1}
 \cM_{\underline{\Sigma},J}(x_0;x_1,\dots,x_d)=\bigsqcup_r\big(\cM^r_{\underline{\Sigma},J}(x_0;x_1,\dots,x_d)/Aut(S_r)\big)
\end{alignat*}
The moduli space $\cM_{\underline{\Sigma},J}(x_0;x_1,\dots,x_d)$ can be viewed as the kernel of a section of a Banach bundle. The linearization of this section at a point $u\in\cM_{\underline{\Sigma},J}(x_0;x_1,\dots,x_d)$ is a Fredholm operator. Then, the almost complex structure $J$ is called \textit{regular} if this operator is surjective. In this case, $\cM_{\underline{\Sigma},J}(x_0;x_1,\dots,x_d)$  is a smooth manifold whose dimension is the Fredholm index of the linearized operator. We will denote by $\cM^i_{\underline{\Sigma},J}(x_0;x_1,\dots,x_d)$ the moduli space of pseudo-holomorphic curves of index $i$ satisfying the conditions (1)-(6) above. Moreover, for $u\in\cM^i_{\underline{\Sigma},J}(x_0;x_1,\dots,x_d)$, we denote $\ind(u):=i$. 

In the following subsections, in order to simplify notations we will not indicate the almost complex structure we use to define the moduli spaces.

\subsubsection{Pseudo-holomorphic curves with boundary on cylindrical cobordisms}
Let us consider $d+1$ Legendrian submanifolds $\La_1,\dots,\La_{d+1}$. The choice of Lagrangian label for disks takes values in the set of Lagrangian cylinders $\{\R\times\La_1,\R\times\La_2,\dots,\R\times\La_{d+1}\}$ and the set of asymptotics consists of Reeb chords from $\La_i$ to $\La_j$ for $1\leq i,j\leq d+1$. To simplify notations for Lagrangian labels we will denote $\R\times\La_{1,\dots,d+1}=(\R\times\La_{1},\dots,\R\times\La_{d+1})$. Moreover this label will indicate only the Lagrangians associated to mixed Reeb chords.
Let $\gamma_{d+1,1}\in\Rc(\La_{d+1},\La_1)$, $\gamma_i\in\Rc(\La_i,\La_{i+1})\cup\Rc(\La_{i+1},\La_i)$ and $\bs{\delta}_i$ be words of pure Reeb chords of $\La_i$, for $1\leq i\leq d+1$. Considering a cylindrical almost complex structure on $\R\times Y$, we define
\begin{alignat*}{1}
\cM_{\R\times\La_{1,\dots,d+1}}(\gamma_{d+1,1};\bs{\delta}_1,\gamma_1,\bs{\delta}_2,\gamma_2,\dots,\bs{\delta}_d,\gamma_d,\bs{\delta}_{d+1})
\end{alignat*}
to be the moduli space of pseudo-holomorphic disks with boundary on $\R\times\La_{1,\dots,d+1}$ that have a positive asymptotic to $\gamma_{d+1,1}$, positive asymptotic to the chord $\gamma_i$ if $\gamma_i$ respects the ordering of Legendrians (i.e. if $\gamma_i\in\Rc(\La_i,\La_{i+1})$), negative asymptotic to the chord $\gamma_i$ if $\gamma_i$ does not respect the ordering of Legendrians, (i.e. $\gamma_i\in\Rc(\La_{i+1},\La_i)$), and negative asymptotics to the pure chords forming the words $\bs{\delta}_i$. There is an action of $\R$ by translation on this moduli space and we denote the quotient by:
\begin{alignat*}{1}
\widetilde{\cM}_{\R\times\La_{1,\dots,d+1}}(\gamma_{d+1,1};\bs{\delta}_1,\gamma_1,\dots,\bs{\delta}_d,\gamma_d,\bs{\delta}_{d+1}):=\cM_{\R\times\La_{1,\dots,d+1}}(\gamma_{d+1,1};\bs{\delta}_1,\gamma_1,\dots,\bs{\delta}_d,\gamma_d,\bs{\delta}_{d+1})/\R
\end{alignat*}
\begin{rem} We have the following particular cases:
	\begin{itemize}\item when $d=0$: the asymptotics are all pure chords of $\La_1$. In the case where only the first asymptotic is positive, such moduli spaces are used to define the differential of the Legendrian contact homology of $\La_1$ (Section \ref{LCHdef}).
		\item when $d=1$: by denoting $\xi_{i,j}$ a chord from $\La_i$ to $\La_j$, the disks in the moduli spaces $\widetilde{\cM}_{\R\times\La_{1,2}}(\gamma_{2,1};\bs{\delta}_1,\xi_{2,1},\bs{\delta}_2)$ (i.e. when only the first asymptotic is positive) are involved in the definition of the Legendrian contact homology of $\La_1\cup\La_2$, and the disks in the moduli spaces $\widetilde{\cM}_{\R\times\La_{1,2}}(\gamma_{2,1};\bs{\delta}_1,\xi_{1,2},\bs{\delta}_2)$ (i.e. when the two mixed chords are positive asymptotics) are called \textit{bananas} and are involved in the definition of the Cthulhu complex (Section \ref{cthulhu}).
	\end{itemize} 
\end{rem}
We will also consider the same type of moduli spaces but with the condition that the first asymptotic is this time a negative Reeb chord asymptotic. Namely, the moduli spaces
\begin{alignat*}{1}
\cM_{\R\times\La_{1,\dots,d+1}}(\gamma_{1,d+1};\bs{\delta}_1,\gamma_1,\bs{\delta}_2,\gamma_2,\dots,\bs{\delta}_d,\gamma_d,\bs{\delta}_{d+1})
\end{alignat*}
where $\gamma_{1,d+1}\in\Rc(\La_{1},\La_{d+1})$ is a negative asymptotic, $\gamma_i\in\Rc(\La_i,\La_{i+1})\cup\Rc(\La_{i+1},\La_i)$ are positive or negative asymptotics depending if the chord respects or not the ordering of Legendrians, and the pure Reeb chords forming the words $\bs{\delta}_i$ are negative asymptotics. As we will see in Section \ref{action_energie}, for energy reasons such moduli spaces are empty if all the $\gamma_i$ are negative asymptotics.

\subsubsection{Pseudo-holomorphic curves with boundary on non cylindrical cobordisms}
\label{esp_de_mod_prod}\label{moduli_spaces}
We consider now moduli spaces of pseudo-holomorphic curves with boundary on the Lagrangians $\Sigma_1,\dots,\Sigma_{d+1}$, where $\La_i^-\prec_{\Sigma_i}\La_i^+$. The choice of Lagrangian label $\underline{\Sigma}$ takes values in $\{\Sigma_1,\dots,\Sigma_{d+1}\}$ and the set of asymptotics consists of Reeb chords from $\La_i^\pm$ to $\La_j^\pm$ for $1\leq i,j\leq d+1$, and intersection points in $\Sigma_i\cap\Sigma_j$ for $1\leq i\neq j\leq d+1$. Again to simplify notations, for Lagrangian labels we will now denote $\Sigma_{1,\dots,d+1}=(\Sigma_{1},\dots,\Sigma_{d+1})$ and this indicates only the Lagrangians associated to mixed asymptotics, i.e. intersection points and chords from a Legendrian to another one.

Given $\gamma_{d+1,1}\in\Rc(\La_{d+1}^+,\La_1^+)$, $a_i\in\{\Sigma_i\cap\Sigma_{i+1}\}\cup\Rc(\La_{i+1}^-,\La_i^-)$ and $\bs{\delta}_i$ words of pure Reeb chords of $\La_i^-$, we consider the moduli spaces
\begin{alignat}{1}
\cM_{\Sigma_{1,\dots,d+1}}(\gamma_{d+1,1};\bs{\delta}_1,a_1,\bs{\delta}_2,a_2,\dots,\bs{\delta}_d,a_d,\bs{\delta}_{d+1})\label{t1}
\end{alignat}
of pseudo-holomorphic disks with boundary on $\Sigma_{1,\dots,d+1}$ that have a positive asymptotic to $\gamma_{d+1,1}$, and asymptotic to intersection points or Reeb chords $a_i$ (negative asymptotics), and negative asymptotics to chords forming the words $\bs{\delta}_i$.
\begin{rem}
	When $d=0$, the Lagrangian label consists of one cobordism $\Sigma_1$ and so the set of asymptotics consists only on Reeb chords of $\La_1^\pm$. In this case, the moduli spaces above are involved in the definition of the differential graded algebra map induced by an exact Lagrangian cobordism from the Chekanov-Eliashberg algebra of $\La_1^+$ to the Chekanov-Eliashberg algebra of $\La_1^-$ (see Section \ref{morphisme_induit}).
\end{rem}	
We consider also for $x\in\Sigma_1\cap\Sigma_{d+1}$ the moduli spaces of pseudo-holomorphic disks
\begin{alignat}{1}
\cM_{\Sigma_{1,\dots,d+1}}(x;\bs{\delta}_1,a_1,\bs{\delta}_2,a_2,\dots,\bs{\delta}_d,a_d,\bs{\delta}_{d+1})\label{t2}
\end{alignat}
with the same asymptotic conditions as before for the $a_i$'s and $\bs{\delta}_i$'s (intersection point or negative Reeb chord asymptotic at $a_i$ and negative Reeb chords asymptotics at chords forming the words $\bs{\delta}_i$).

Finally, given $\gamma_{1,d+1}\in\Rc(\La_1^-,\La_{d+1}^-)$, we consider the moduli spaces
\begin{alignat}{1}
\cM_{\Sigma_{1,\dots,d+1}}(\gamma_{1,d+1};\bs{\delta}_1,a_1,\bs{\delta}_2,a_2,\dots,\bs{\delta}_d,a_d,\bs{\delta}_{d+1})\label{t3}
\end{alignat}
of pseudo-holomorphic disks with a negative asymptotic to $\gamma_{1,d+1}$, and again the same asymptotic conditions as before for the $a_i$'s and $\bs{\delta}_i$'s. Again in this last case, if the $a_i$'s are all negative Reeb chords asymptotics, this moduli space will be empty for energy reasons.
\begin{rem}
	When $d=1$, these three types of moduli spaces \eqref{t1}, \eqref{t2}, \eqref{t3} are involved in the definition of the Cthulhu complex. In particular, the map $\Fc^1$ of Theorem \ref{teofoncteur} is defined by a mod-$2$ count of curves in moduli spaces of type \eqref{t1} (See Section \ref{cthulhu}). For $d=2$, the curves in moduli spaces of type \eqref{t2} and \eqref{t3} are useful to define the product structure of Theorem \ref{teo_prod} (Section \ref{defprod}) and the curves in moduli spaces of type \eqref{t1} appear in the definition of an order-2 map $\Fc^2$ in the proof of Theorem \ref{teofoncteur} (Section \ref{teo2}).
\end{rem}

\subsection{Action and energy}\label{action_energie}
Consider $d+1$ transverse exact Lagrangian cobordisms $(\Sigma_1,\dots,\Sigma_{d+1})$. Recall that by definition, associated to each cobordism there is a function $f_i\colon\Sigma_i\to\R$, primitive of the form $e^t\alpha_{|\Sigma_i}$, and this function is constant on the cylindrical ends of $\Sigma_i$. Without loss of generality, we can consider that the constants in the negative ends of the cobordisms are zero, and we denote $\mathfrak{c}_i$ the constant for the positive end of $\Sigma_i$. We also denote $T>0$ and $\epsilon>0$ such that the cobordisms $\Sigma_i$ are all cylindrical out of $\Sigma_i\cap([-T+\epsilon,T-\epsilon]\times Y)$. To each asymptotic, we can associate a quantity called \textit{action} as follows.
For an intersection point $x\in\Sigma_i\cap\Sigma_j$ with $i>j$, the action of $x$ is given by:
\begin{alignat*}{1}
 \mathfrak{a}(x)=f_i(x)-f_j(x)
\end{alignat*}
For a Reeb chord $\gamma$, the \textit{length} of $\gamma$ is given by $\ell(\gamma):=\int_\gamma\alpha$ and then the action of $\gamma_{i,j}^+\in\Rc(\La_i^+,\La_j^+)$ is defined by:
\begin{alignat*}{1}
 \mathfrak{a}(\gamma_{i,j}^+)=e^T\ell(\gamma_{i,j}^+)+\mathfrak{c}_i-\mathfrak{c}_j
\end{alignat*}
and for a Reeb chord $\gamma_{i,j}^-\in\Rc(\La_i^-,\La_j^-)$ we set:
\begin{alignat*}{1}
 \mathfrak{a}(\gamma_{i,j}^-)=e^{-T}\ell(\gamma_{i,j}^-)\\
\end{alignat*}
Remark that Reeb chords have always a positive action whereas intersection points can be of negative action.
Then, to a pseudo-holomorphic curve $u$ in $\cM_\Sigma(x_0;x_1,\dots,x_d)$ is associated an \textit{energy}, which is the analogue of the area for the case of pseudo-holomorphic curves in compact symplectic manifolds. To define it, let $\chi\colon\R\to\R$ be a function such that:

\begin{alignat*}{1}
 \left\{\begin{array}{ll} 
		 \chi(t)=e^{t} & \mbox{ if }t\in[-T+\epsilon,T-\epsilon]\\
		 \lim\limits_{t\to+\infty}\chi(t)=e^{T} &\\
		 \lim\limits_{t\to-\infty}\chi(t)=e^{-T} &\\
                 \chi'(t)>0 & 
                \end{array}\right.
\end{alignat*}
We define then the $d(\chi\alpha)$-energy of a pseudo-holomorphic curve $u\colon S_r\to\R\times Y$ by:
\begin{alignat*}{1}
 E_{d(\chi\alpha)}(u)=\int_{S_r}u^*d(\chi\alpha)
\end{alignat*}
We have the following very standard result:
\begin{lem}\label{energie_pos}
 $E_{d(\chi\alpha)}(u)\geq0$
\end{lem}
\proof
The $d(e^t\alpha)$-compatibility of the almost complex structure $J$ implies the $d\alpha_{|\xi}$-compatibility of the restriction of $J$ to the contact structure $(\xi,(d\alpha)_{|\xi})$. This permits to show that $E_{d(\chi\alpha)}(u)=\frac{1}{2}\int_{S_r}|du|^2$, where $|v|^2=d(\chi\alpha)(v,Jv)$ is strictly positive if $v\neq0$.
\endproof

Now, the energy of a pseudo-holomorphic curve can be expressed in terms of the actions of its asymptotics.
                                                                                                               
\begin{prop}\label{energie}
	We have the following:
	\begin{enumerate}
		\item if $u\in\cM_{\R\times\La_{1,\dots,d+1}}(\gamma_{d+1,1};\bs{\delta}_1,\gamma_1,\bs{\delta}_2,\gamma_2,\dots,\bs{\delta}_d,\gamma_d,\bs{\delta}_{d+1})$, let us consider the following partition of $\{1,\dots,d\}$ into two subsets:
		\begin{alignat*}{1}
		&I^+=\{i\,|\,\gamma_i\mbox{ positive Reeb chord asymptotic of }u\}\\
		&I^-=\{i\,|\,\gamma_i\mbox{ negative Reeb chord asymptotic of }u\}
		\end{alignat*}
		then we have
		\begin{alignat*}{1}
		E_{d(\chi\alpha)}(u)=\mathfrak{a}(\gamma_{d+1,1})+\sum\limits_{i\in I^+}\mathfrak{a}(\gamma_i)-\sum\limits_{i\in I^-}\mathfrak{a}(\gamma_i)-\sum\limits_{i=1}^{d+1}\mathfrak{a}(\bs{\delta}_i),                  \end{alignat*}
	\item if $u\in\cM_{\R\times\La_{1,\dots,d+1}}(\gamma_{1,d+1};\bs{\delta}_1,\gamma_1,\bs{\delta}_2,\gamma_2,\dots,\bs{\delta}_d,\gamma_d,\bs{\delta}_{d+1})$,
		\begin{alignat*}{1}
		E_{d(\chi\alpha)}(u)=-\mathfrak{a}(\gamma_{1,d+1})+\sum\limits_{i\in I^+}\mathfrak{a}(\gamma_i)-\sum\limits_{i\in I^-}\mathfrak{a}(\gamma_i)-\sum\limits_{i=1}^{d+1}\mathfrak{a}(\bs{\delta}_i),                  \end{alignat*}
	\item if $u\in\cM_{\Sigma_{1,\dots,d+1}}(\gamma_{d+1,1};\bs{\delta}_1,a_1,\bs{\delta}_2,a_2,\dots,\bs{\delta}_d,a_d,\bs{\delta}_{d+1})$,
		\begin{alignat*}{1}
		E_{d(\chi\alpha)}(u)=\mathfrak{a}(\gamma_{d+1,1})-\sum\limits_{i=1}^d\mathfrak{a}(a_i)-\sum\limits_{i=1}^{d+1}\mathfrak{a}(\bs{\delta}_i),
		\end{alignat*}
	\item if $u\in\cM_{\Sigma_{1,\dots,d+1}}(x;\bs{\delta}_1,a_1,\bs{\delta}_2,a_2,\dots,\bs{\delta}_d,a_d,\bs{\delta}_{d+1})$,
		\begin{alignat*}{1}
		E_{d(\chi\alpha)}(u)=\mathfrak{a}(x)-\sum\limits_{i=1}^d\mathfrak{a}(a_i)-\sum_{i=1}^{d+1}\mathfrak{a}(\bs{\delta}_i),
		\end{alignat*}
	\item if $u\in	\cM_{\Sigma_{1,\dots,d+1}}(\gamma_{1,d+1};\bs{\delta}_1,a_1,\bs{\delta}_2,a_2,\dots,\bs{\delta}_d,a_d,\bs{\delta}_{d+1})$,
		\begin{alignat*}{1}
		E_{d(\chi\alpha)}(u)=-\mathfrak{a}(\gamma_{1,d+1})-\sum\limits_{i=1}^{d}\mathfrak{a}(a_i)-\sum\limits_{i=1}^{d+1}\mathfrak{a}(\bs{\delta}_i)
		\end{alignat*}
	\end{enumerate}
\end{prop}
Lemma \ref{energie_pos} and Proposition \ref{energie} give thus some constraints on the action of asymptotics of pseudo-holomorphic curves. These will be useful in order to cancel some pseudo-holomorphic configurations in Section \ref{PROD}.

\subsection{Compactness}\label{comp}
When transversality holds, i.e. when the almost complex structure is regular for moduli spaces, these are smooth manifolds which are not necessarily compact. However, they admit a compactification in the sense of Gromov (\cite{Gr}), by adding broken curves called \textit{pseudo-holomorphic buildings}. Compactness results together with transversality results imply that the compactification of a moduli space is a compact manifold whose boundary components are in bijection with pseudo-holomorphic buildings arising as degeneration of pseudo-holomorphic curves in the moduli space. We recall below the definition of pseudo-holomorphic buildings whose components are pseudo-holomorphic disks with boundary on Lagrangian cobordisms with cylindrical ends (see \cite{BEHWZ} and \cite{Abbas}).

Given again $d+1$ transverse exact Lagrangian cobordisms $\La_i^-\prec_{\Sigma_i}\La_i^+$, we consider the following Lagrangian labels $\underline{\Sigma}=(\Sigma_1,\dots,\Sigma_{d+1})$ and $\underline{\R\times\Lambda}^\pm=(\R\times\La_1^\pm,\dots,\R\times\La_{d+1}^\pm)$.
Given a planar rooted tree $T$ with $d-1$ leaves, $d\geq2$, we can associate to each interior vertex $v$ a triple $(S_{r_v},I_v,L_v)$, where:
\begin{itemize}
	\item $S_{r_v}$ is the Riemann disk associated to the vertex $v$ (Section \ref{deligne}),
	\item $I_v$ is the set of boundary marked points of $S_{r_v}$ called \textit{nodes} corresponding to interior edges of $T$,
	\item $L_v$ is the set of boundary marked points of $S_{r_v}$ corresponding to edges connecting $v$ to the leaves and the root of $T$.
\end{itemize}
\begin{ex}
	 For $T=T_{d+1}$, $T$ is a rooted tree with one interior vertex $v$ and $d$ leaves. In this case $I_v=\emptyset$ and $L_v$ consists of $d+1$ elements, one corresponding to the root, the other to the leaves. For a tree having an interior vertex $v$ adjacent only to interior vertices, we have $L_v=\emptyset$ and $I_v$ contains $|v|$ elements.
\end{ex}

Moreover, the set of all nodes $\bigcup\limits_{v} I_v$ contains an even number of elements organized in pair. Indeed, each interior edge $e_i$ of $T$ is by definition connecting two interior vertices $v$ and $\bar{v}$ of $T$. Denote by $p_i$ the node on the boundary of $S_{r_v}$ corresponding to $e_i$, and by $\bar{p}_i$ the node on the boundary of $S_{r_{\bar{v}}}$ corresponding to $e_i$. So we get that $\bigcup\limits_{v} I_v$ admits a partition 
\begin{alignat*}{1}
	\bigcup\limits_{v} I_v=\{p_1,\bar{p}_1\}\cup\{p_2,\bar{p}_2\}\cup\dots\cup\{p_k,\bar{p}_k\}
\end{alignat*}
where $k$ is the number of interior edges of $T$.
\begin{defi}\label{phb}
	Given integers $k^-,k^+\geq0$, a \textit{pseudo-holomorphic building of height $k^-|1|k^+$} in $\R\times Y$ with boundary on $\underline{\Sigma}$ is the data of:
	\begin{enumeratec}
		\item a planar rooted tree $T$ and the corresponding union of triples $\bigcup\limits_{v}(S_{r_v},I_v,L_v)$,
		\item a choice of asymptotic in $A(\underline{\Sigma})$ for each node in $\bigcup I_v$. We require that for each pair $\{p_j,\bar{p}_j\}$, the same asymptotic is assigned to $p_j$ and $\bar{p}_j$.
		\item a choice of asymptotic for each marked point in $\bigcup L_v$,
		\item a pair $(u_v,\rho_v)$ for each interior vertex $v$ of $T$, where $u_v:S_{r_v}\to\R\times Y$ is a pseudo-holomorphic disk asymptotic to the given asymptotics assigned to elements in $I_v\cup L_v$, and $\rho_v$ is an integer called the \textit{floor} of $v$, satisfying $-k^-\leq\rho_v\leq k^+$, such that
		\begin{enumeratec}
			\item each floor $-k^-\leq\rho\leq k^+$ except $\rho=0$ admits at least one non trivial disk,
			\item if $\rho_v>0$: $u_v$ has boundary on $\underline{\R\times\La}^+$, and thus the asymptotics assigned to elements of $I_v\cup L_v$ are only Reeb chords. The disk $u_v$ is said to live in the \textit{top level}.
			\item if $\rho_v=0$: $u_v$ has boundary on the compact parts of $\underline{\Sigma}$ (i.e. on $\overline{\Sigma}_1\cup\overline{\Sigma}_2\cup\dots\cup\overline{\Sigma}_{d+1}$) and can have asymptotics to intersection points and Reeb chords in the negative and positive ends. The disk $u_v$ is said to live in the \textit{middle level}.
			\item if $\rho_v<0$: $u_v$ has boundary on $\underline{\R\times\La}^-$, and has asymptotics only to Reeb chords. The disk $u_v$ is said to live in the \textit{bottom level}.
	\suspend{enumeratec}
	Moreover, we require the following conditions on nodes:
	\resume{enumeratec}
			\item if $S_{r_v}$ has a boundary puncture at a node $p_i$ such that $u_v$ has a positive (resp negative) asymptotic to a Reeb chord $\gamma\in A(\underline{\Sigma})$ at $p_i$, then the corresponding pair $(u_{\bar{v}},\rho_{\bar{v}})$ in the building (where $\bar{v}$ is the interior vertex of $T$ such that the interior edge $e_i$ connects $v$ and $\bar{v}$) satisfies $\rho_{\bar{v}}=\rho_v+1$ (resp $\rho_{\bar{v}}=\rho_v-1$) and $u_{\bar{v}}$ has a negative (resp positive) asymptotic to $\gamma$ at $\bar{p}_i$,
			\item if $S_{r_v}$ has a boundary puncture at $p_i$ such that $u_v$ is asymptotic at $p_i$ to an intersection point $x\in A(\underline{\Sigma})$, which is a jump from $\Sigma_l$ to $\Sigma_m$ when traversing the boundary of $u_v$ counterclockwise, then the corresponding pair $(u_{\bar{v}},\rho_{\bar{v}})$ in the building satisfies $\rho_{\bar{v}}=0$ and $u_{\bar{v}}$ has an asymptotic to $x$ at $\bar{p}_i$, where $x$ is a jump from $\Sigma_m$ to $\Sigma_l$ when traversing the boundary of $u_{\bar{v}}$ counterclockwise,
	\suspend{enumeratec}
	and the following conditions on marked points which are not nodes:
	\resume{enumeratec}
			\item if the middle level ($\rho=0$) of the building is not empty, then for each element in $L_v$ asymptotic to $y_v$,
			\begin{enumeratec}
				\item if $y_v$ is a Reeb chord in $A(\underline{\R\times\La^+})$, then $y_v$ is a positive Reeb chord asymptotic of $u_v$ and $\rho_v\geq0$,
				\item if $y_v$ is an intersection point, then $\rho_v=0$,
				\item if $y_v$ is a Reeb chord in $A(\underline{\R\times\La^-})$, then $y_v$ is a negative Reeb chord asymptotic of $u_v$ and $\rho_v\leq0$.
			\end{enumeratec}
	\end{enumeratec}
\end{enumeratec}
\end{defi}
\begin{rem}
	If the middle level is empty, then we have a building with boundary on the cylindrical ends of the Lagrangians. It is a building of height $0|0|k^+$ if all components have boundary on $\underline{\R\times \La}^+$, and a building of height $k^-|0|0$ if all components have boundary on $\underline{\R\times \La}^-$.
\end{rem}

\begin{defi}[Equivalence of pseudo-holomorphic buildings]\label{def:eqphb}
	Two pseudo-holomorphic buildings are \textit{equivalent} if they become the same after the removal of an appropriate number of trivial cylinders together with the obvious deformation of the underlying planar rooted tree, to each of them in the bottom and top levels.
	In other words, two pseudo-holomorphic buildings are equivalent if they become the same after:
	\begin{itemize}
		\item removing all possible trivial cylinders in the bottom and top levels attached to asymptotics assigned to $\cup L_v$,
		\item removing simultaneously trivial cylinders attached to all the positive ends and/or all the negative ends corresponding to nodes in $\cup I_v$ of a component in the bottom or top level.
	\end{itemize}
\end{defi}
	In particular, two equivalent pseudo-holomorphic buildings have the same components in the middle level, and the same non trivial components in the bottom and top levels but connected to each other by a certain number of trivial cylinders which can vary from one building to the other.

\begin{figure}[ht]     
      \begin{center}\includegraphics[width=7cm]{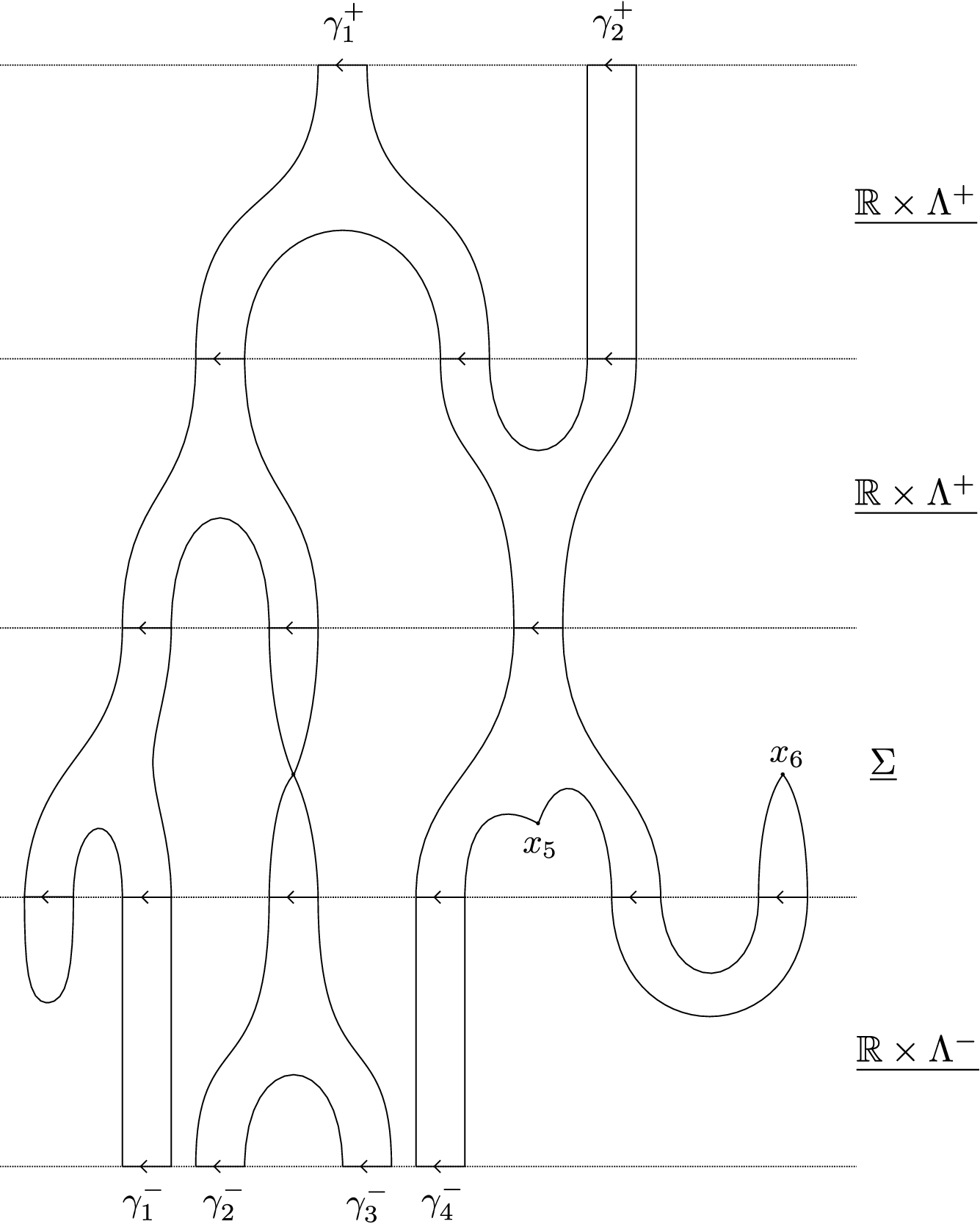}\end{center}
      \caption{Example of a pseudo-holomorphic building of height $1|1|2$.}
      \label{bat_holo}
\end{figure}
Performing the gluing operation on each interior edge of $T$ corresponds to do the connected sums of the disks $S_{r_{v}}$ at each node, as described in Section \ref{deligne} (i.e. identifying standard neighbourhoods of $p_i$ and $\bar{p}_i$, for $1\leq i\leq k$). The boundary marked points of the disk resulting from the connected sum are asymptotic to the asymptotics assigned to the marked points in $\bigcup L_v$. Given a set of asymptotics $(x_0,x_1,\dots,x_d)$, we denote by $\cM_{\underline{\Sigma}}^{k^-|1|k^+}(x_0;x_1,\dots,x_d)$ the set of pseudo-holomorphic buildings of height $k^-|1|k^+$ with boundary on $\underline{\Sigma}$ such that the disk obtained after boundary connected sum of the domains at nodes is asymptotic to $(x_0,x_1,\dots,x_d)$.
Moduli spaces of pseudo-holomorphic disks with boundary on non-cylindrical Lagrangian cobordisms as described in Section \ref{esp_de_mod} can be viewed as pseudo-holomorphic buildings of height $0|1|0$ with boundary on $\underline{\Sigma}$, in other words we have
\begin{alignat*}{1}
 \cM_{\underline{\Sigma}}(x_0;x_1,\dots,x_d)\subset\cM_{\underline{\Sigma}}^{0|1|0}(x_0;x_1,\dots,x_d)
\end{alignat*}

By Gromov's compactness, a sequence of pseudo-holomorphic disks $u_s$ in $\cM_{\underline{\Sigma}}(x_0;x_1,\dots,x_d)$ admits a subsequence which converges to a pseudo-holomorphic building with boundary on $\underline{\Sigma}$.
The pseudo-holomorphic buildings obtained this way are of two types:
\begin{enumerate}
 \item \textit{Stable breaking}: pseudo-holomorphic building such that each component is a curve having at least three mixed asymptotics. For example, a pseudo-holomorphic building in a product $$\cM(x_0;x_1,...,x_{i-1},x_i',x_{i+j},...,x_d)\times\cM(x_i';x_i,...,x_{i+j-1})$$
 with $1\leq i\leq d-1$ and $j\geq2$.
 \item \textit{Unstable breaking}: pseudo-holomorphic building having at least a curve with at most two mixed asymptotics. Such a curve is either a pseudo-holomorphic half-plane (so without mixed asymptotic), or a pseudo-holomorphic strip.
\end{enumerate}
The important result is that the set of buildings asymptotic to $x_0,x_1,\dots,x_d$ gives a compactification of the moduli space $\cM_{\underline{\Sigma}}(x_0;x_1,\dots,x_d)$, i.e. the disjoint union 
\begin{alignat*}{1}
 \bigsqcup_{k^-,k^+\geq0}\cM_{\underline{\Sigma}}^{k^-|1|k^+}(x_0;x_1,\dots,x_d)
\end{alignat*}
is compact.

Assume that we have an admissible regular almost complex structure. Consider a pseudo-holomorphic building in $\cM_{\underline{\Sigma}}^{k^-|1|k^+}(x_0;x_1,\dots,x_d)$ given by a tree $T$ and pseudo-holomorphic disks $\{u_i\}$ with choices of asymptotics for nodes. Gluing results (\cite{EES1}) imply that there exists sufficiently small gluing parameters associated to each interior edge of $T$ (see Section \ref{deligne}) and a unique pseudo-holomorphic curve $u$ in $\cM_{\underline{\Sigma}}(x_0;x_1,\dots,x_d)$ depending on these parameters such that $u$ converges to the pseudo-holomorphic building $\{u_i\}$ when the gluing parameters converge to $0$ (see \cite{Abbas} and \cite{BEHWZ} for the topology on the space of pseudo-holomorphic buildings). In the rest of the paper, when we are given a pseudo-holomorphic building consisting of components $\{u_i\}$, by abuse of language we will say that we glue its components meaning that we consider a corresponding pseudo-holomorphic disk $u$ as above.

\section{Legendrian contact homology}\label{LCH}
Legendrian contact homology is a Legendrian isotopy invariant which has been defined by Chekanov \cite{Che} and Eliashberg \cite{E} independently. Eliashberg gave a definition of Legendrian contact homology in the setting of Symplectic Field Theory (SFT, see \cite{EGH}). Chekanov defined it combinatorially for Legendrian links in $\R^3$, by a count of certain types of convex polygons with boundary on the Lagrangian projection of the Legendrian. Then, Ekholm, Etnyre and Sullivan in \cite{EES1,EES2} generalized the definition of Chekanov for Legendrian submanifolds in $\R^{2n+1}$ and $P\times\R$, by counting pseudo-holomorphic disks with boundary on the Lagrangian projection. In fact, it has been proven by Etnyre, Ng and Sabloff \cite{ENS} in dimension $3$ and then by Dimitrolgou-Rizell \cite{DR} in every dimension that the SFT-version of Legendrian contact homology computes the same invariant as the combinatorial version of Chekanov and its generalization in higher dimension. With this in mind, we recall below the SFT-definition of Legendrian contact homology, which is more in the spirit of this article.

\subsection{The differential graded algebra}\label{LCHdef}
Given a Legendrian submanifold $\La\subset P\times\R$, we denote by $C(\La)$ the $\Z_2$-vector space generated by Reeb chords of $\La$, and $\Ac(\La)=\bigoplus_iC(\La)^{\otimes i}$ the tensor algebra of $C(\La)$, called the \textit{Chekanov-Eliashberg algebra} of $\La$. There is a grading associated to Reeb chords and defined from the Conley-Zehnder index by the following: if $\La$ is connected and $c\in\Rc(\La)$ then we set $|c|:=\nu_{\gamma_c}(c)-1$, where $\gamma_c$ is a capping path for $c$. This is a well-defined grading in $\Z$ modulo the Maslov number of $\La$ (because of the choice of the capping path) and twice the first Chern class of $TP$ (because of the choice of a symplectic trivialization of $TP$ along $\Pi_P(\gamma_c)$ to compute $\nu_{\gamma_c}(c)$), see \cite{EES2} for more details. This induces a grading for each word of Reeb chords in $\Ac(\La)$ by $|b_1b_2\dots b_m|=\sum_i|b_i|$ for Reeb chords $b_i$.
If $\La$ is not connected and $c$ is a mixed chord with ends $c^+\in\La^+$ and $c^-\in\La^-$, where $\La^\pm$ are connected components of $\La$, in order to define the grading we choose some points $p^\pm\in\La^\pm$ and some paths $\gamma_c^+\subset\La^+$ and $\gamma_c^-\subset\La^-$ from $c^+$ to $p^+$ from $p^-$ to $c^-$ respectively. Then we choose a path $\gamma_{+-}$ from $p^+$ to $p^-$ and so if we denote $\Gamma_c=\gamma_c^+\cup\gamma_{+-}\cup\gamma_c^-$ the concatenation of the three paths, the degree of $c$ is defined to be $|c|=\nu_{\Gamma_c}(c)-1$. The grading of mixed chords depend on the paths $\gamma_{+-}$ but for two mixed chords $c_1,c_2$ from $\La^-$ to $\La^+$, the difference $|c_1|-|c_2|$ does not depend on $\gamma_{+-}$.

Let $J$ be a cylindrical almost complex structure on $\R\times Y$. The differential on $\Ac(\La)$ is a map $\partial\colon\Ac(\La)\to\Ac(\La)$ which is defined by a count of pseudo-holomorphic disks in $\R\times Y$ with boundary on $\R\times\La$ and asymptotic to Reeb chords. More precisely, if $a\in\Rc(\La)$:
\begin{alignat}{1}
 \partial(a)=\sum\limits_{m\geq0}\sum\limits_{\substack{\bs{b}=b_1\cdots b_m\\ |\bs{b}|=|a|-1}}\#\widetilde{\cM}_{\R\times\La}(a;\bs{b})\cdot \bs{b} \label{differential}
\end{alignat}
where $a$ is a positive asymptotic and the chords $b_i$ are negative Reeb chord asymptotics. When $m=0$ we set $\bs{b}=1$. We then extend this to the whole algebra by the Leibniz rule.

About transversality results, Dimitroglou-Rizell proved in \cite{DR2} that generically, a cylindrical almost complex structure on $\R\times Y$ is regular for the moduli spaces $\cM_{\R\times\La}(a;b_1,\dots,b_m)$ which are thus manifolds of dimension
\begin{alignat*}{1}
 \dim\cM_{\R\times\La}(a;b_1,\dots,b_m)=|a|-\sum|b_i|
\end{alignat*}
and so
\begin{alignat*}{1}
 \dim\widetilde{\cM}_{\R\times\La}(a;b_1,\dots,b_m)=|a|-\sum|b_i|-1
\end{alignat*}
This is done by generalizing a result of Dragnev (\cite{Dra}) to the case of pseudo-holomorphic disks, using the fact that as pseudo-holomorphic curves in the moduli spaces above have only one positive Reeb chord asymptotic, it is always possible to find an injective point.
These dimension formula imply that in the definition of the differential \eqref{differential}, this is a mod-$2$ count of pseudo-holomorphic disks in $0$-dimensional moduli spaces. Then, by Gromov's compactness these $0$-dimensional moduli spaces are compact and thus the differential $\partial$ is a well-defined map of degree $-1$.
Transversality also holds for almost complex structures that are cylindrical lifts of regular compatible almost complex structures on $P$ (satisfying a technical condition near the intersection points of $\Pi_P(\La)$, see \cite{DR2,EES2}).
\begin{teo}\cite{Che,EES1,EES2,DR}
 \begin{itemize}
  \item[$\bullet$] $\partial^2=0$,
  \item[$\bullet$] The Legendrian contact homology $LCH_*(\La,J)$ does not depend on a generic choice of cylindrical almost complex structure $J$ and is a Legendrian isotopy invariant.
 \end{itemize}
\end{teo}

Legendrian contact homology being generally of infinite dimension, we recall in the next section the linearization process introduced by Chekanov, in order to extract finite dimensional (and so more computable) invariants from Legendrian contact homology.

\subsection{Linearization and the augmentation category}\label{linearisation}
We begin this section by recalling the fundamental tool for the linearization: augmentations.
\begin{defi}
 An \textit{augmentation} for $(\Ac(\La),\partial)$ is a unital DGA-map $\ep\colon\Ac(\La)\to\Z_2$ where $\Z_2$ is viewed as a DGA with vanishing differential. In other words, $\ep$ is a map satisfying:
 \begin{itemize}
 \item $\ep(1)=1$,
  \item $\ep(a)=0$ if $|a|\neq0$,
  \item $\ep(ab)=\ep(a)\ep(b)$,
  \item $\ep\circ\partial=0$.
 \end{itemize}
\end{defi}
A Legendrian submanifold does not necessarily admit an augmentation. Typically, once there is an element of the algebra $a\in\Ac(\La)$ such that $\partial a=1$, the third condition in the definition above cannot be satisfied and hence there is no augmentation. For example, loose Legendrians (see \cite{M}) do not admit augmentation. In this paper, we will only focus on Legendrians whose Chekanov-Eliashberg algebra can be augmented. So let us consider a Legendrian submanifold $\La\subset\R\times Y$ such that $\Ac(\La)$ admits an augmentation, then it is possible to associate to $\La$ a new complex $(C(\La),\partial_1^\ep)$, with $\partial_1^\ep$ defined on Reeb chords by:
\begin{alignat*}{1}
 \partial_1^\ep(a)=\sum\limits_{m\geq0}\sum\limits_{\substack{\bs{b}=b_1\cdots b_m\\ |\bs{b}|=|a|-1}}\sum_{i=1}^m\#\widetilde{\cM}_{\R\times\La}(a;\bs{b})\cdot \ep(b_1)\dots\ep(b_{i-1})\ep(b_{i+1})\dots\ep(b_m)\cdot b_i
\end{alignat*}
In fact, conjugating the differential $\partial$ by the DGA-morphism $g_\ep$ defined on chords by $g_\ep(c)=c+\ep(c)$ gives a new differential $\partial^\ep$ on $\Ac(\La)$, the differential $\partial$ \textit{twisted by $\ep$}, such that the restriction on $C(\La)$ can be decomposed as $\partial^\ep_{|C(\La)}=\sum_{i\geq0}\partial^\ep_i$, with $\partial^\ep_i\colon C(\La)\to C(\La)^{\otimes i}$. But the differential $\partial^\ep_{|C(\La)}$ does not admit any constant term (i.e. $\partial^\ep_0=0$) due to the properties of $\ep$, and so $(\partial^\ep_{|C(\La)})^2=0$ implies that $(\partial_1^\ep)^2=0$. The homology of the complex $(C(\La),\partial_1^\ep)$ is by definition the Legendrian contact homology of $\La$ linearized by $\ep$.
\begin{teo}\cite{Che}
 The set $\{LCH_*^\ep(\La),\,\ep\}$ of linearized Legendrian contact homologies is a Legendrian isotopy invariant.
\end{teo}
This linearization process can be done using two augmentations instead of one (see \cite{BCh}), leading to the bilinearized Legendrian contact homology $LCH^{\ep_1,\ep_2}_*(\La)$, which is the homology of the complex $(C(\La),\partial_1^{\ep_1,\ep_2})$ with
\begin{alignat*}{1}
 \partial_1^{\ep_1,\ep_2}(a)=\sum\limits_{m\geq0}\sum\limits_{\substack{\bs{b}=b_1\cdots b_m\\ |\bs{b}|=|a|-1}}\sum_{i=1}^m\#\widetilde{\cM}_{\R\times\La}(a;\bs{b})\cdot \ep_1(b_1)\dots\ep_1(b_{i-1})\ep_2(b_{i+1})\dots\ep_2(b_m)\cdot b_i
\end{alignat*}
The advantage of the bilinearized version in comparison to the linearized one is that it retains some information about the non-commutativity of the Chekanov-Eliashberg DGA.
More generally, given $d+1$ augmentations $\ep_1,\dots,\ep_{d+1}$ of $\Ac(\La)$, there is a map
\begin{alignat*}{1}
 \partial_d^{\ep_1,\dots,\ep_{d+1}}\colon C(\La)\to C(\La)^{\otimes d}
\end{alignat*}
such that $\partial_d^{\ep_1,\dots,\ep_{d+1}}(a)$ is a sum of words of length $d$ coming from words in $\partial a$ to which we keep $d$ letters and augment the others by $\ep_1,\dots,\ep_{d+1}$ in this order (changing the augmentation each time we jump a chord we keep). In all the rest of the article, we will adopt a cohomology point of view, so let us describe the dual maps of the maps $\partial^{\ep_1,\dots,\ep_{d+1}}_d$. As the vector space $C(\La)$ and its dual are canonically isomorphic, by an abuse of notation we will still denote $C(\La)$ the dual vector space. So the dual of $\partial^{\ep_1,\dots,\ep_{d+1}}_d$, denoted $\mu_{\ep_{d+1},\dots,\ep_1}^d$, is defined by:
\begin{alignat*}{1}
 \mu_{\ep_{d+1},\dots,\ep_1}^d(b_d,\dots,b_1)=\sum_{a\in\Rc(\La)}\sum_{\substack{\bs{\delta}_1,...,\bs{\delta}_{d+1}\\ \sum|b_i|+\sum|\bs{\delta}_i|=|a|-1}}\#\widetilde{\cM}_{\R\times\La}(a;\bs{\delta}_1,b_1,\bs{\delta}_2,\dots,\bs{\delta}_d,b_d,\bs{\delta}_{d+1})\ep_1(\bs{\delta}_1)\dots\\
 \dots\ep_{d+1}(\bs{\delta}_{d+1})\cdot a
\end{alignat*}
where $\bs{\delta}_i$ are words of Reeb chords of $\La$. In fact, as already explained above for the dual map, the coefficient $\langle\mu_{\ep_{d+1},\dots,\ep_1}^d(b_d,\dots,b_1),a\rangle$ is computed by considering all words of length at least $d$ in $\partial(a)$ containing the letters $b_1,\dots,b_d$ in this order, and augmenting the (possibly) remaining chords between $b_i$ and $b_{i+1}$ by $\ep_{i+1}$, for all $1\leq i\leq d$. These maps $\{\mu^d_{\ep_{d+1},\dots,\ep_1}\}_{d\geq1}$ satisfy the $A_\infty$-relations, i.e. for all $d\geq1$ and Reeb chords $b_d,\dots,b_1$ we have
\begin{alignat}{1}
 \sum_{\substack{1\leq j\leq d\\ 0\leq n\leq d-j}}\mu^{d-j+1}_{\ep_{d+1},\dots,\ep_{n+j+1},\ep_{n+1},\dots,\ep_1}(b_d,\dots,\mu^j_{\ep_{n+j+1},\dots,\ep_{n+1}}(b_{n+j},\dots,b_{n+1}),b_n,\dots,b_1)=0\label{relAug}
\end{alignat}
and thus the maps $\{\mu^d_{\ep_{d+1},\dots,\ep_1}\}_{d\geq1}$ are $A_\infty$-composition maps of an $A_\infty$-category called the \textit{augmentation category} of $\La$, denoted $\Aug_-(\La)$. This category has been defined by Bourgeois and Chantraine in \cite{BCh} as follows:
\begin{itemize}
 \item $\Ob(\Aug_-(\La)):$ $\ep$ augmentation of $\Ac(\La)$,
 \item $\hom(\ep_1,\ep_2)=(C(\La),\mu^1_{\ep_2,\ep_1})$ the bilinearized Legendrian contact cohomology complex,
 \item the $A_\infty$-composition maps are the maps $\mu^d_{\ep_{d+1},\dots,\ep_1}$ defined above.
\end{itemize}
If we look at the full subcategory generated by one object $\ep$, then we get the $A_\infty$-algebra $(C(\La),\{\mu^d_\ep\}_{d\geq1})$ that appeared first in a work of Civan, Etnyre, Koprowski, Sivek and Walker \cite{CKESW}.

In fact, the $A_\infty$-maps of the augmentation category can be viewed as dual maps of components of the differential of the $(d+1)$-copy of $\La$ twisted by a particular augmentation. This is a way to show that the $A_\infty$-relations are satisfied, using a bijection between moduli spaces with boundary on $\La$ and moduli spaces with boundary on the $k$-copy of $\La$ (see \cite[Theorem 3.6]{EESa}). For $k\geq 1$, the $k$-copy of $\La$ denoted $\La_{(k)}$ is defined as follows. Set $\La_1:=\La$, and for a small $\epsilon>0$  we define $\widetilde{\La}_j:=\varphi^R_{(j-1)\epsilon}(\La)$ for $2\leq j\leq k$, where $\varphi_t^R$ is the Reeb flow (recall $R_\alpha=\partial_z$ here). The Legendrian submanifold $\La_1\cup\widetilde{\La}_2\cup\dots\cup\widetilde{\La}_k$ has an infinite number of Reeb chords, so we have to perturb it to turn it into a chord generic Legendrian. Take Morse functions $f_j\colon\La\to\R$, for $2\leq j\leq k$, such that the functions $f_j-f_i$ are Morse. Then, identify a small tubular neighborhood of $\widetilde{\La_j}$ to a neighborhood of the $0$-section in $J^1(\La)$, and replace $\widetilde{\La_j}$ by the $1$-jet of $f_j$ which is by definition the submanifold $j^1(f_j)=\{(q,d_qf_j,f_j(q))\,|\,q\in\La\}\subset J^1(\La)$. We denote this new Legendrian $\La_j\subset P\times\R$. The $k$-copy of $\La$ is defined to be the union $\La_1\cup\La_2\cup\dots\cup\La_k$. It is a chord generic Legendrian which has four different types of Reeb chords:
\begin{enumerate}
 \item \textit{pure chords}: chords of $\La_j$, for all $1\leq j\leq k$, and there is a bijection between $\Rc(\La)$ and $\Rc(\La_j)$,
 \item \textit{Morse chords}: mixed chords corresponding to critical points of the functions $f_j$ for $j\geq2$ (these are chords from $\La_1$ to $\La_j$) and $f_j-f_i$ (these are chords from $\La_i$ to $\La_j$ for $i<j$),
 \item \textit{small chords}: mixed chords from $\La_j$ to $\La_i$ for $i<j$, in bijection with chords of $\La$,
 \item \textit{long chords}: mixed chords (which are not Morse) from $\La_i$ to $\La_j$ for $i<j$ also in bijection with chords of $\La$.
\end{enumerate}
For a chord $a\in\Rc(\La)$, denote by $a_{i,j}$ the corresponding Reeb chord in $\Rc(\La_i,\La_j)$. Let $(\ep_1,\dots,\ep_k)$ be augmentations of $\Ac(\La)$ and consider the DGA-morphism $\ep_{(k)}\colon\Ac(\La_{(k)})\to\Z_2$ defined on Reeb chords by:
\begin{alignat*}{1}
 &\ep_{(k)}(a_{i,i})=\ep_i(a)\\
 &\ep_{(k)}(a_{i,j})=0\mbox{ for }i\neq j\\
 &\ep_{(k)}(c_M)=0\mbox{ for }c_M\mbox{ Morse chord}
\end{alignat*}
It is shown in \cite{BCh} that $\ep_{(k)}$ is an augmentation of $(\Ac(\La_{(k)}),\partial_{(k)})$, that we call \textit{diagonal augmentation induced by $\ep_1,\dots,\ep_k$}. Also, by denoting $I_M$ the two-sided ideal of $\Ac(\La_{(k)})$ generated by Morse chords, we have $\partial_{(k)}(I_M)\subset I_M$ \cite[Proposition 3.1]{BCh}. This implies that $\partial_{(k)}$ descends to a differential on the quotient $\Ac(\La_{(k)})/I_M$ which by abuse of notation we still denote $\partial_{(k)}$, and $\ep_{(k)}$ descends to an augmentation of $\Ac(\La_{(k)})/I_M$, which we still denote $\ep_{(k)}$. Now, denoting $C(\La_i,\La_j)$ the $\Z_2$-vector space generated by Reeb chords from $\La_j$ to $\La_i$ in $\Ac(\La_{(k)})/I_M$, given a diagonal augmentation as above we have that the twisted differential $\partial^{\ep_{(k)}}$ restricted to $C(\La_1,\La_k)$ is a map:
\begin{alignat*}{1}
 \partial^{\ep_{(k)}}_{|C(\La_1,\La_k)}\colon C(\La_1,\La_k)\to\bigoplus_{\substack{d\geq1\\ 1\leq i_2,\dots,i_d\leq k}} C(\La_{i_d},\La_k)\otimes C(\La_{i_{d-1}},\La_{i_d})\otimes\dots\otimes C(\La_1,\La_{i_2})
\end{alignat*}
The dual of the appropriate component of this map is then $\mu^{d}_{\ep_k,\ep_{i_d},\dots,\ep_{i_2},\ep_1}$ (see \cite[Theorem 3.2]{BCh} and \cite[Theorem 3.6]{EESa} for the correspondence of moduli spaces in the case of the $2$-copy). Then, dualizing the relation $$\big(\partial^{\ep_{(k)}}_{|C(\La_1,\La_k)}\big)^2=0$$ gives all the $A_\infty$-relations for $d\leq k-1$, i.e. the $A_\infty$-relations for each sequence of objects $(\ep_1,\ep_{i_2},\dots,\ep_{i_{d}},\ep_k)$. For example, the two first are:
\begin{alignat*}{1}
 &\big(\mu^1_{\ep_k,\ep_1}\big)^2=0\\
 &\mu^1_{\ep_k,\ep_1}\circ\mu^2_{\ep_k,\ep_{i},\ep_1}+\mu^2_{\ep_k,\ep_{i},\ep_1}\big(\mu^1_{\ep_k,\ep_i}\otimes\id\big)+\mu^2_{\ep_k,\ep_{i},\ep_1}\big(\id\otimes\mu^1_{\ep_i,\ep_1}\big)=0, \mbox{ for all }1\leq i\leq k.
\end{alignat*}

\subsection{Morphism induced by a cobordism}\label{morphisme_induit}
Given an exact Lagrangian cobordism $\La^-\prec_\Sigma\La^+$, there exists a DGA-map $\phi_\Sigma\colon\Ac(\La^+)\to\Ac(\La^-)$ defined on Reeb chords by
\begin{alignat*}{1}
 \phi_\Sigma(\gamma^+)=\sum_{\gamma_1^-,\dots,\gamma_m^-}\#\cM^0_\Sigma(\gamma^+;\gamma_1^-,\dots,\gamma_m^-)\cdot\gamma_1^-\cdots\gamma_m^-
\end{alignat*}
where $\gamma^+\in\Rc(\La^+)$ and $\gamma_1^-,\dots,\gamma_m^-\in\Rc(\La^-)$ (see \cite{EHK}). When $\Sigma$ is an exact Lagrangian filling of $\La^+$ ($\La^-=\emptyset$), then the DGA-map $\phi_\Sigma\colon\Ac(\La^+)\to\Z_2$ is an augmentation of $\Ac(\La^+)$. In this case, the corresponding linearized Legendrian contact cohomology of $\La$ is determined by the topology of the filling by the Ekholm-Seidel isomorphism:
\begin{teo}[\cite{E1,DR}]\label{EkhS}
 If $\La\subset Y$ is a n-dimensional closed Legendrian submanifold which admits a Lagrangian filling $\Sigma$, then $H_*(\Sigma)\simeq LCH_\ep^{n-*}(\La)$, where $\ep$ is the augmentation of $\Ac(\La)$ induced by $\Sigma$.
\end{teo}
This theorem gives a very powerful obstruction to the existence of Maslov $0$ Lagrangian fillings. For example, once the linearized Legendrian contact cohomology of a Legendrian $\La$ has a generator of degree strictly less than $0$ or strictly higher than $n$ for all possible augmentations, it means that $\La$ is not fillable by an exact Lagrangian.
More generally, given an augmentation $\ep^-$ of $\Ac(\La^-)$, its pullback by $\phi_\Sigma$ is an augmentation of $\Ac(\La^+)$ that we denote $\ep^+:=\ep^-\circ\phi_\Sigma$. This is the order-$0$ map of a family of maps defining an $A_\infty$-functor $\Phi_\Sigma\colon\Aug(\La^-)\to\Aug(\La^+)$ as follows (see \cite{BCh}):
\begin{itemize}
 \item on the objects of the category, $\Phi_\Sigma(\ep^-)=\ep^-\circ\phi_\Sigma$,
 \item for each $(d+1)$-tuple $(\ep_1^-,\dots,\ep_{d+1}^-)$ of augmentations of $\Ac(\La^-)$, there is a map
 \begin{alignat*}{1}
  \Phi^{\ep_{d+1}^-,\dots,\ep_1^-}_\Sigma\colon\hom(\ep_d^-,\ep^-_{d+1})\otimes\dots\otimes\hom(\ep^-_1,\ep^-_2)\to\hom(\ep^+_1,\ep^+_{d+1})
 \end{alignat*}
defined by
\begin{alignat*}{1}
 \Phi^{\ep_{d+1}^-,...,\ep_1^-}_\Sigma(\gamma_d^-,\dots,\gamma_1^-)=\sum_{\substack{\gamma^+\in\Rc(\La^+)\\ \bs{\delta}_1^-,...,\bs{\delta}_{d+1}^-}}\#\cM^0_\Sigma(\gamma^+;\bs{\delta}_1^-,\gamma_1^-,\bs{\delta}_2^-,\dots,\gamma_d^-,\bs{\delta}_{d+1}^-)\cdot\ep_1^-(\bs{\delta}_1^-)\\
 \dots\ep^-_{d+1}(\bs{\delta}_{d+1}^-)\cdot\gamma^+
\end{alignat*}
\end{itemize}
The induced functor on cohomology level gives a map on bilinearized Legendrian contact cohomology:
\begin{alignat*}{1}
 \Phi_\Sigma^{\ep_2^-,\ep_1^-}\colon LCH_{\ep_1^-,\ep_2^-}^*(\La^-)\to LCH_{\ep_1^+,\ep_2^+}^*(\La^+)
\end{alignat*}
which was shown to be an isomorphism if $\Sigma$ is a concordance, in \cite{CDGG1}.

In the case of the augmentation category $\Aug_+(\La)$ (defined in \cite{NRSSZ}), an exact Lagrangian cobordism from $\La^-$ to $\La^+$ also induces a functor $F:\Aug_+(\La^-)\to\Aug_+(\La^+)$. In particular, this functor was shown to be injective on equivalence classes of augmentations and cohomologically faithful by Yu Pan \cite{Yu}.

\section{Floer theory for Lagrangian cobordisms}\label{CTH}
\subsection{The Cthulhu complex}\label{cthulhu}
\begin{figure}[ht]     
      \begin{center}\includegraphics[width=14.5cm]{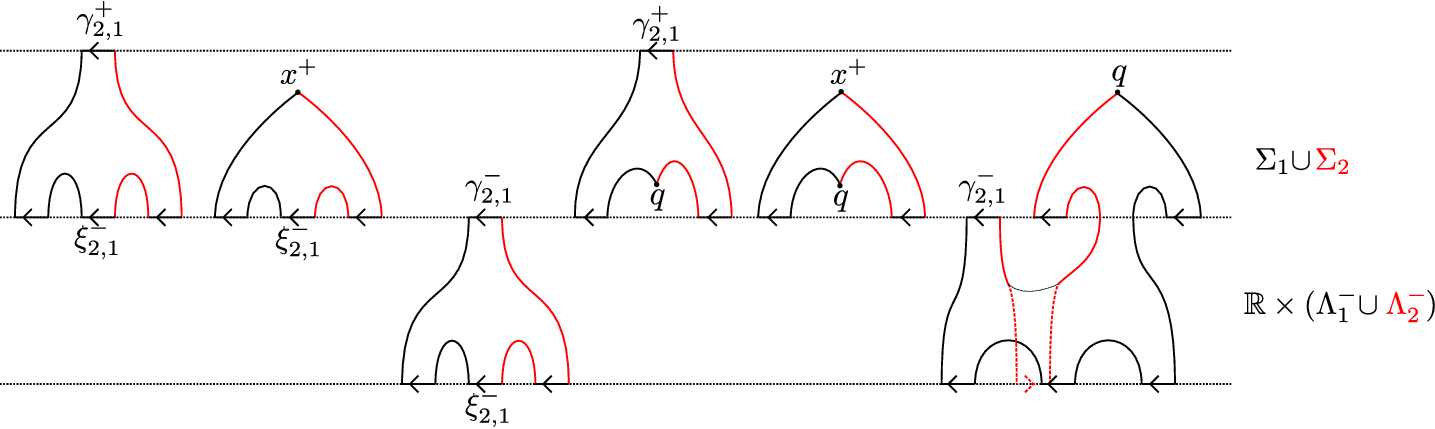}\end{center}
      \caption{From left to right: schematic picture of pseudo-holomorphic disks contributing to $d_{+-}(\xi_{2,1}^-)$, $d_{0-}(\xi_{2,1}^-)$, $d_{--}(\xi_{2,1}^-)$, $d_{+0}(q)$, $d_{00}(q)$ and $d_{-0}(q)$.
      }
      \label{Cthulhu}
\end{figure}

In this section we recall the definition of the Cthulhu complex, a Floer-type complex for Lagrangian cobordisms, defined by Chantraine, Dimitroglou-Rizell, Ghiggini and Golovko in \cite{CDGG2}. Let $\La_1^-\prec_{\Sigma_1}\La_1^+$ and $\La_2^-\prec_{\Sigma_2}\La_2^+$ be two transverse exact Lagrangian cobordisms in $\R\times Y$ with $\La_1^-, \La_1^+,\La_2^-,\La_2^+$ Legendrian submanifolds in $Y$. We assume that the Chekanov-Eliashberg algebras $\Ac(\La_1^-)$ and $\Ac(\La_2^-)$ admit augmentations $\ep_1^-$ and $\ep_2^-$ respectively, which induce augmentations $\ep_1^+$ and $\ep_2^+$ of $\Ac(\La_1^+)$ and $\Ac(\La_2^+)$ as we saw previously. Cthulhu homology is the homology of a graded complex associated to the pair $(\Sigma_1,\Sigma_2)$, denoted $(\Cth(\Sigma_1,\Sigma_2),\mathfrak{d}_{\ep_1^-,\ep_2^-})$, generated by intersection points in $\Sigma_1\cap\Sigma_2$, Reeb chords from $\La_2^+$ to $\La_1^+$ and Reeb chords from $\La_2^-$ to $\La_1^-$, with shifts in grading:
\begin{alignat*}{1}
 \Cth(\Sigma_1,\Sigma_2)=C(\La_1^+,\La_2^+)[2]\oplus CF(\Sigma_1,\Sigma_2)\oplus C(\La_1^-,\La_2^-)[1]
\end{alignat*}
where $CF(\Sigma_1,\Sigma_2)$ denotes the $\Z_2$-vector space generated by intersection points in $\Sigma_1\cap\Sigma_2$.
This is a graded complex. For the grading to be in $\Z$, we assume that $2c_1(P)$ as well as the Maslov classes of $\Sigma_1$ and $\Sigma_2$ vanish (which implies that the Maslov classes of $\La_1^\pm$ and $\La_2^\pm$ also vanish). The grading for Reeb chords is the same as the Legendrian contact homology grading (Section \ref{LCHdef}). For an intersection point $p\in\Sigma_1\cap\Sigma_2$ the grading is defined to be the Maslov index of a path of graded Lagrangians from $(T_p\Sigma_1)^\#$ to $(T_p\Sigma_2)^\#$ in $Gr^\#(T_p(\R\times Y),d(e^t\alpha)_p)$, the universal cover of the Grassmaniann of Lagrangian subspaces of $(T_p(\R\times Y),d(e^t\alpha)_p)$, see \cite[Section 11.j]{S} and \cite{CDGG2}. The differential on $\Cth(\Sigma_1,\Sigma_2)$ is a matrix 
\begin{alignat*}{1}
 \mathfrak{d}_{\ep_1^-,\ep_2^-}=\begin{pmatrix} d_{++} & d_{+0} & d_{+-}\\ 0 & d_{00} & d_{0-} \\ 0 & d_{-0} & d_{--}\end{pmatrix}
\end{alignat*}
where each component is defined by a count of rigid pseudo-holomorphic disks with boundary on $\Sigma_1$ and $\Sigma_2$ and asymptotic to intersection points and Reeb chords from $\La_2^\pm$ to $\La_1^\pm$ as follows:
\begin{enumerate}
 \item for $\xi^+_{2,1}\in\Rc(\La_2^+,\La_1^+)$: 
  \begin{alignat*}{1}
   d_{++}(\xi^+_{2,1})=\sum_{\gamma^+_{2,1}}\sum_{\bs{\beta}_1,\bs{\beta}_2}\#\widetilde{\cM^1}_{\R\times\La^+_{1,2}}(\gamma_{2,1}^+;\bs{\beta}_1,\xi_{2,1}^+,\bs{\beta}_2)\ep_1^+(\bs{\beta}_1)\ep_2^+(\bs{\beta}_2)\cdot\gamma_{2,1}^+
  \end{alignat*}
  where the sum is for $\gamma^+_{2,1}\in\Rc(\La_2^+,\La_1^+)$ and $\bs{\beta}_i$ words of Reeb chords of $\La^+_i$, for $i=1,2$. The map $d_{++}$ is the restriction to $C(\La_1^+,\La_2^+)$ of the bilinearized differential $\mu^1_{\ep_2^+,\ep_1^+}$ of the Legendrian contact cohomology of $\La_1^+\cup\La_2^+$.
  \item for $\xi^-_{2,1}\in\Rc(\La_2^-,\La_1^-)$:
  \begin{alignat*}{1}
   d_{+-}(\xi^-_{2,1})&=\sum_{\gamma^+_{2,1}}\sum_{\bs{\delta}_1,\bs{\delta}_2}\#\cM^0_{\Sigma_{1,2}}(\gamma^+_{2,1};\bs{\delta}_1,\xi_{2,1}^-,\bs{\delta}_2)\ep_1^-(\bs{\delta}_1)\ep_2^-(\bs{\delta}_2)\cdot\gamma_{2,1}^+\\
   d_{0-}(\xi_{2,1}^-)&=\sum_{x^+\in\Sigma_1\cap\Sigma_2}\sum_{\bs{\delta}_1,\bs{\delta}_2}\#\cM^0_{\Sigma_{1,2}}(x^+;\bs{\delta}_1,\xi_{2,1}^-,\bs{\delta}_2)\ep_1^-(\bs{\delta}_1)\ep_2^-(\bs{\delta}_2)\cdot x^+\\
   d_{--}(\xi_{2,1}^-)&=\sum_{\gamma^-_{2,1}}\sum_{\bs{\delta}_1,\bs{\delta}_2}\#\widetilde{\cM^1}_{\R\times\La^-_{1,2}}(\gamma^-_{2,1};\bs{\delta}_1,\xi_{2,1}^-,\bs{\delta}_2)\ep_1^-(\bs{\delta}_1)\ep_2^-(\bs{\delta}_2)\cdot\gamma_{2,1}^-\\
  \end{alignat*}
  and as for $d_{++}$, the map $d_{--}$ is the restriction of $\mu^1_{\ep_2^-,\ep_1^-}$ to $C(\La_1^-,\La_2^-)$.
 \item for $q\in\Sigma_1\cap\Sigma_2$ which is a jump from $\Sigma_1$ to $\Sigma_2$:
  \begin{alignat*}{2}
   d_{+0}(q)&=\sum_{\gamma^+_{2,1}}\sum_{\bs{\delta}_1,\bs{\delta}_2}\#\cM^0_{\Sigma_{1,2}}(\gamma^+_{2,1};\bs{\delta}_1,q,\bs{\delta}_2)\ep_1^-(\bs{\delta}_1)\ep_2^-(\bs{\delta}_2)\cdot\gamma_{2,1}^+\\
   d_{00}(q)&=\sum_{x^+}\sum_{\bs{\delta}_1,\bs{\delta}_2}\#\cM^0_{\Sigma_{1,2}}(x^+;\bs{\delta}_1,q,\bs{\delta}_2)\ep_1^-(\bs{\delta}_1)\ep_2^-(\bs{\delta}_2)\cdot x^+\\
    d_{-0}(q)&=b\circ\delta_{-0}^{\Sigma_{2,1}}(q)\\
   &=\sum_{\gamma^-_{1,2},\gamma^-_{2,1}}\sum_{\bs{\delta}_1,\bs{\delta}_2}\#\big(\widetilde{\cM^1}_{\R\times\La_{1,2}^-}(\gamma^-_{2,1};\bs{\delta}_1',\gamma^-_{1,2},\bs{\delta}_2'')\times\cM^0_{\Sigma_{1,2}}(\gamma_{1,2}^-;\bs{\delta}''_1,q,\bs{\delta}_2')\big)\\
   & \hspace{67mm}\cdot\ep_1^-(\bs{\delta}_1'\bs{\delta}_1'')\ep_2^-(\bs{\delta}_2'\bs{\delta}_2'')\cdot \gamma^-_{2,1}
  \end{alignat*}
where
\begin{itemize}
 \item the last sum is for $\bs{\delta}_i,\bs{\delta}_i',\bs{\delta}_i''$ words of Reeb chords of $\La_i$ such that $\bs{\delta}_i=\bs{\delta}_i'\bs{\delta}_i''$,
 \item $\delta_{-0}^{\Sigma_{2,1}}$ is the dual of $d_{0-}^{\Sigma_{2,1}}\colon C(\La_2^-,\La_1^-)\to CF(\Sigma_2,\Sigma_1)$ with Lagrangian label $(\Sigma_2,\Sigma_1)$,
 \item $b\colon C(\La_2^-,\La_1^-)\to C(\La_1^-,\La_2^-)$ is the map defined by the count of bananas:
 \begin{alignat*}{1}
 b(\gamma^-_{1,2})=\sum_{\gamma^-_{2,1}}\sum_{\bs{\delta}_1,\bs{\delta}_2}\#\widetilde{\cM^1}_{\R\times\La_{1,2}^-}(\gamma^-_{2,1};\bs{\delta}_1,\gamma^-_{1,2},\bs{\delta}_2)\cdot\ep_1^-(\bs{\delta}_1)\ep_2^-(\bs{\delta}_2)\cdot \gamma^-_{2,1}
\end{alignat*}
\end{itemize}
\end{enumerate}
See Figure \ref{Cthulhu} for examples of pseudo-holomorphic disks which contribute to the components of $\mathfrak{d}_{\ep_1^-,\ep_2^-}$, except for $d_{++}$ whose contributing disks are of the same type of those contributing to $d_{--}$, but with boundary on $\R\times(\La_1^+\cup\La_2^+)$.
\begin{rem}
 The other components of the differential vanish for energy reasons.
\end{rem}
When transversality holds, it is again possible to express the dimension of the moduli spaces above by the degree of the asymptotics. In particular, from \cite[Proposition 3.2]{CDGG2} we have:
\begin{alignat*}{1}
 & \dim\widetilde{\cM}_{\R\times\La^+_{1,2}}(\gamma_{2,1}^+;\bs{\beta}_1,\xi_{2,1}^+,\bs{\beta}_2)=|\gamma_{2,1}^+|+|\xi_{2,1}^+|-|\bs{\beta}_1|-|\bs{\beta}_2|-1\\
 & \dim\widetilde{\cM}_{\R\times\La^-_{1,2}}(\gamma_{2,1}^-;\bs{\delta}_1,\zeta_{2,1}^-,\bs{\delta}_2)=|\gamma_{2,1}^-|-|\zeta_{2,1}^-|-|\bs{\delta}_1|-|\bs{\delta}_2|-1 \\
 & \dim\widetilde{\cM}_{\R\times\La^-_{1,2}}(\gamma_{2,1}^-;\bs{\delta}_1,\gamma_{1,2}^-,\bs{\delta}_2)= |\gamma_{2,1}^-|+|\gamma_{1,2}^-|-|\bs{\delta}_1|-|\bs{\delta}_2|+1-n\\
 & \dim\cM_{\Sigma_{1,2}}(\gamma_{2,1}^+;\bs{\delta}_1,q,\bs{\delta}_2)=|\gamma_{2,1}^+|-|q|-|\bs{\delta}_1|-|\bs{\delta}_2|+1\\
 & \dim\cM_{\Sigma_{1,2}}(\gamma_{2,1}^+;\bs{\delta}_1,\zeta_{2,1}^-,\bs{\delta}_2)=|\gamma_{2,1}^+|-|\zeta_{2,1}^-|-|\bs{\delta}_1|-|\bs{\delta}_2|\\
 & \dim\cM_{\Sigma_{1,2}}(x^+;\bs{\delta}_1,q,\bs{\delta}_2)=|x^+|-|q|-|\bs{\delta}_1|-|\bs{\delta}_2|-1 \\
 & \dim\cM_{\Sigma_{1,2}}(x^+;\bs{\delta}_1,\zeta_{2,1}^-,\bs{\delta}_2)=|x^+|-|\zeta_{2,1}^-|-|\bs{\delta}_1|-|\bs{\delta}_2|-2 \\
 &
 \dim\cM_{\Sigma_{1,2}}(\gamma_{1,2}^-;\bs{\delta}_1,q,\bs{\delta}_2)=n+1-|q|-|\gamma_{2,1}^-|-|\bs{\delta}_1|-|\bs{\delta}_2|-1\\ 
 \end{alignat*}
This gives that the map $\mathfrak{d}_{\ep_1^-,\ep_2^-}$ is of degree $1$. Without the shifts in grading, we obtain that $d_{+0}$ is of degree $-1$, $d_{+-}$ and $d_{-0}$ are of degree $0$, $d_{++}$, $d_{00}$ and $d_{--}$ are degree $1$ maps and $d_{0-}$ is of degree $2$.

The necessary transversality results in order to make the above moduli spaces transversely cut out are given in \cite{CDGG2}. Briefly, as we already saw in the previous section, for Legendrian contact homology-type moduli spaces, cylindrical almost complex structures on $\R\times Y$ are generically regular. This is also the case for moduli spaces of bananas, since that even if the curves in those spaces have two positive Reeb chords asymptotics, these Reeb chords are distinct, and so the curve is always somewhere injective.

Now, if $J^\pm$ are regular for Legendrian contact homology type moduli spaces and banana moduli spaces, then moduli spaces $\cM_\Sigma(\gamma^+;\gamma_1^-,\dots,\gamma_d^-)$ are transversely cut out for a generic almost complex structure in $\Jc^{adm}_{J^-,J^+,T}(\R\times Y)$, using results of \cite[Chapter 3]{DS}. The regularity assumption on $J^\pm$ permits in particular to achieve transversality for pseudo-holomorphic curves coming from the gluing of a curve in $\cM_\Sigma(\gamma^+;\gamma_1^-,\dots,\gamma_d^-)$ and a curve in $\cM_{\R\times\La^\pm}(\gamma;\gamma_1,\dots,\gamma_d)$.

Finally, moduli spaces of the types $\cM_{\Sigma_{1,2}}(x^+;\bs{\delta}_1,q,\bs{\delta}_2)$ and $\cM_{\Sigma_{1,2}}(\gamma^+_{2,1};\bs{\delta}_1,\gamma^-_{2,1},\bs{\delta}_2)$ are transversely cut out for a generic domain dependent almost complex structure $$J:[0,1]\to\Jc^{adm}_{J^-,J^+,T}(\R\times Y)$$
generalizing results of \cite{AD}. The domain dependence here is just a time-dependence because the domain of a curve is biholomorphic to a strip $\R\times[0,1]$ with marked points on the boundary (asymptotic to pure Reeb chords), and we want invariance of the almost complex structure by translation of the $\R$-coordinate.

\begin{teo}\cite{CDGG2} Given $\Sigma_1,\Sigma_2\subset\R\times Y$ exact Lagrangian cobordisms as above,
 \begin{enumerate}
  \item $\mathfrak{d}_{\ep_1^-,\ep_2^-}^2=0$, and
  \item The complex $(\Cth(\Sigma_1,\Sigma_2),\mathfrak{d}_{\ep_1^-,\ep_2^-})$ is acyclic.
 \end{enumerate}
\end{teo}
The first point of this theorem is proven by studying breakings of pseudo-holomorphic curves of index $1$ with boundary on $\Sigma_1\cup\Sigma_2$, or of index $2$ with boundary on $\R\times\La_1^-\cup\R\times\La_2^-$, and two mixed asymptotics. In Section \ref{proof}, we will use the same ideas to prove Theorem \ref{teo_prod}. The second point of the theorem comes from the fact that it is possible to displace the cobordisms in $\R\times Y$ such that $\Sigma_1$ and $\Sigma_2$ no longer have intersection points and such that there are no more Reeb chords from $\La_2^\pm$ to $\La_1^\pm$. Briefly, this is done by first wrapping the ends of one of the two cobordisms by a Hamiltonian isotopy in such a way that the complex we get has only intersection points generators (no more Reeb chords) and is canonically isomorphic to the original Cthulhu complex. Then, the invariance of the Cthulhu complex by a compactly supported Hamiltonian isotopy permits to separate the two cobordisms so that there are no more generators, which implies that the complex vanishes, as well as its homology.

Let us denote $\partial_{-\infty}=\begin{pmatrix} d_{00} & d_{0-} \\ d_{-0} & d_{--}\end{pmatrix}$ the submatrix of $\mathfrak{d}_{\ep_1^-,\ep_2^-}$, then 
\begin{alignat}{1}0=(\mathfrak{d}_{\ep_1^-,\ep_2^-})^2=\left(\begin{array}{cc}
         d_{++}^2 & \begin{array}{cc}
                     *_{+0} & *_{+-}\\
                    \end{array}\\
	 \begin{array}{c}
	      0\\ 0\\
	 \end{array} 
	 & \partial_{-\infty}^2
        \end{array}\right) \label{del_carre}
\end{alignat}
where $*_{+0}=d_{++}d_{+0}+d_{+0}d_{00}+d_{+-}d_{-0}$ and $*_{+-}=d_{++}d_{+-}+d_{+0}d_{0-}+d_{+-}d_{--}$. So in particular, $(C(\La_1^+,\La_2^+),d_{++})$ is a subcomplex of the Cthulhu complex and
\begin{alignat*}{1}
\big(CF_{-\infty}(\Sigma_1,\Sigma_2):=CF(\Sigma_1,\Sigma_2)\oplus C(\La_1^-,\La_2^-)[1],\partial_{-\infty}\big)
\end{alignat*}
is a quotient complex. Relation \eqref{del_carre} implies also that $d_{+0}+d_{+-}\colon CF_{-\infty}(\Sigma_1,\Sigma_2)\to C(\La_1^+,\La_2^+)$ is a chain map, i.e. the Cthulhu complex is the cone of $d_{+0}+d_{+-}$. This map, that we denote now $\Fc_{21}^1$ and sometimes just $\Fc^1$ when the pair of cobordisms in clear from the context, is in fact a quasi-isomorphism due to the acyclicity of the Cthulhu complex.

\subsection{Hamiltonian perturbations}\label{perturbations}
Given a cobordism $\La^-\prec_{\Sigma}\La^+$ in $(\R\times P\times\R,d(e^t\alpha))$, we consider a special type of Hamiltonian isotopies by which we deform $\Sigma$, in order to extract some properties of the Cthulhu complex. More precisely, we use a Hamiltonian $H\colon\R\times P\times\R\to\R$ that depends only on the real coordinate in the symplectization, which means that $H(t,p,z)=h(t)$, where $h\colon\R\to\R$ is a smooth function. The associated Hamiltonian flow is by definition the flow of the Hamiltonian vector field $X_H$ defined by $\iota_{X_H}d(e^t\alpha)=-dH$. We can compute that $X_H(t,p,z)=e^{-t}h'(t)\partial_z$ and so the flow $\Phi_H$ is given by:
\begin{alignat*}{1}
  &\Phi_H^s\colon\R\times P\times\R\to\R\times P\times\R\\
  &\Phi_H^s(t,p,z)=(t,p,z+se^{-t}h'(t))
\end{alignat*}
Now, since $\Sigma$ is an exact Lagrangian cobordism, $\Phi_H^s(\Sigma)$ is also an exact Lagrangian cobordism. Indeed, if $f_{\Sigma}:\Sigma\to\R$ is the primitive of $e^t\alpha$ restricted to $\Sigma$, then we have: 
\begin{alignat*}{1}
 e^t\alpha_{|\Phi_H^s(\Sigma)}&=(\Phi_H^s)^*(e^t(dz+\beta))\\
 &=e^t(d(z+se^{-t}h'(t))+\beta)_{|\Sigma}\\
 &=e^t(dz+se^{-t}(h''-h')dt+\beta)_{|\Sigma}\\
 &=e^t\alpha_{|\Sigma}+s(h''-h')dt_{|\Sigma}\\
\end{alignat*}
So, a primitive of $e^t\alpha_{|\Phi_H^s(\Sigma)}$ is given by
\begin{alignat}{1}
 f_{\Phi_H^s(\Sigma)}(t,p,z)=f_{\Sigma}(t,p,z)+s(h'-h)(t)\label{nouveauf}
\end{alignat}
In particular, when the function $h$ is for example the function $h_D$ below, the primitive $f_{\Phi_H^s(\Sigma)}$ given by \eqref{nouveauf} vanishes on the negative end of $\Phi_H^s(\Sigma)$.
This type of Hamiltonian isotopy is useful to wrap the cylindrical ends of the cobordisms, and the way to wrap depends on the choice of the function $h\colon\R\to\R$ to define the Hamiltonian. Let us describe here one type of perturbation (see \cite{CDGG2} for other perturbations). Given $T>0$, we define a function $h_D\colon\R\to\R$ by:
\begin{alignat*}{1}
          &h_D(t)=e^t\mbox{ for }t\leq-T-1\\
          &h_D(t)=A\mbox{ for }t\in[-T,T]\\
          &h_D(t)=e^t-B\mbox{ for }t\geq T+1\\
          &h_D'(t)\geq0\mbox{ for }t\in[-T-1,-T]\cup[T,T+1]
          \end{alignat*}
where $A$ and $B$ are positive constants. Then we denote $H_D$ the corresponding Hamiltonian on $\R\times P\times\R$. Now we look how this Hamiltonian wraps the cylindrical ends of a cobordism. Let $\Sigma_1$ be an exact Lagrangian cobordism and consider its image by the flow at time $\epsilon$, $\Phi^\epsilon_{H_D}(\Sigma_1)$, for a small $\epsilon>0$, and denote it $\widetilde{\Sigma}_2$. The cobordisms $\Sigma_1,\widetilde{\Sigma}_2$ are not transverse, indeed:
\begin{enumerate}
	\item on $[-T,T]\times Y$ they coincide,
	\item $\widetilde{\Sigma}_2$ has cylindrical ends $(-\infty,-T-1]\times\widetilde{\Lambda}_2^-$ and $[T+1,\infty)\times\widetilde{\Lambda}_2^+$, where $\widetilde{\La}_2^\pm=\La_1^\pm+\epsilon\frac{\partial}{\partial z}$, so the Legendrian links $\La_1^-\cup\widetilde{\La}_2^-$ and $\La_1^+\cup\widetilde{\La}_2^+$ are degenerate, i.e. they have an infinite number of Reeb chords.
\end{enumerate}
In order to get a pair of transverse cobordisms, we perturb $\widetilde{\Sigma}_2$ as follows. We explain briefly the perturbation we need and refer to \cite[Section 6]{DR} for some more detailed construction. By the Weinstein Lagrangian neighborhood theorem, there is a neighborhood of $\Sigma_1$ symplectomorphic to a neighborhood $U_0$ of the $0$-section of $T^*\Sigma_1$, such that $\Sigma_1$ is identified with the $0$-section. If $\epsilon$ is sufficiently small, then $\widetilde{\Sigma}_2=\Phi_{H_D}^\epsilon(\Sigma_1)$ is identified with a Lagrangian in $U_0$ which can be seen as the graph of $d(\epsilon H_D)$ for the function $\epsilon H_D$ restricted to $\Sigma_1$. In particular, on $\Sigma_1\cap\big([-T,T]\times Y\big)$, the graph of $d(\epsilon H_D)$ coincides with the $0$-section. Take $\Sigma_2\subset\R\times Y$ to be the exact Lagrangian cobordism identified with the graph of $df$, for $f:\Sigma_1\to\R$ a Morse function such that:
\begin{enumerate}
	\item $f$ is a small perturbation of $\epsilon H_D$ on $\Sigma_1$,
	\item the critical points of $f$ are all contained in $\Sigma_1\cap\big([-T,T]\times Y\big)$,
	\item the cylindrical ends $\La_2^\pm$ of $\Sigma_2$ are small Morse perturbations of $\widetilde{\La}_2^\pm$: there are Morse functions $f^\pm$ on $\widetilde{\La}_2^\pm$ such that on a neighborhood of $\widetilde{\La}_2^\pm$ contactomorphic to a neighborhood of the $0$-section of $J^1(\widetilde{\La}_2^\pm)$, $\La_2^\pm$ is identified with $j^1(f^\pm)$.
\end{enumerate}
The pair $(\Sigma_1,\Sigma_2)$ is now a pair of transverse exact Lagrangian cobordisms. For a small enough Morse perturbation as above, Formula \eqref{nouveauf} gives that every intersection point in $CF(\Sigma_1,\Sigma_2)$ has negative action, we say then that the pair $(\Sigma_1,\Sigma_2)$ is \textit{directed}. 
Such a pair of cobordisms satisfy some properties listed in the following proposition:

\begin{prop}\cite{CDGG2}\label{prop_iso_can}
 Let $(\Sigma_1,\Sigma_2)$ be a directed pair of Lagrangian cobordisms such that $\Sigma_2$ is a small perturbation of $\Phi_{H_D}^\epsilon(\Sigma_1)$ as above by a Morse function $f$ on $\Sigma_1$. Let $T>0$ be such that $\Sigma_i\backslash([-T,T]\times Y\cap\Sigma_i)$ are cylindrical, and consider a domain dependent almost complex structure $J_t$ in $\Jc^{adm}_{J^-,J^+,T}(\R\times Y)$ such that $J^\pm$ are in $\Jc^{cyl}_\pi(\R\times Y)$. Assume moreover that $\Ac(\La_1^-)$ admits augmentations $\ep_1^-, \ep_2^-$ which induce augmentations $\ep_1^+$ and $\ep_2^+$ of $\Ac(\La_1^+)$, then:
 \begin{enumerate}
  \item there are canonical isomorphisms of the Chekanov-Eliashberg algebras $(\Ac(\La_1^-),\partial_1^-)\simeq(\Ac(\La_2^-),\partial_2^-)$ and $(\Ac(\La_1^+),\partial_1^+)\simeq(\Ac(\La_2^+),\partial_2^+)$, and so in particular $\ep_1^\pm$ and $\ep^\pm_2$ can also be considered as augmentations of $\Ac(\La_2^\pm)$ under this identification,
  \item $H^*\big(C(\La_1^-,\La_2^-),d_{--}\big)\simeq LCH_{\ep_1^-,\ep_2^-}^*(\La_1^-)$, where $H^*\big(C(\La_1^-,\La_2^-),d_{--}\big)$ denotes the homology of the complex $(C(\La_1^-,\La_2^-),d_{--}\big)$ using augmentations $\ep_1^-,\ep_2^-$ to compute $d_{--}$,
  \item $H^*\big(C(\La_1^+,\La_2^+),d_{++}\big)\simeq LCH_{\ep_1^+,\ep_2^+}^*(\La_1^+)$,
  \item if $J_t$ is regular and induced by a Riemanniann metric $g$ such that $(f,g)$ is a Morse-Smale pair in a neighborhood of $\overline{\Sigma}_1$, then $HF_*(\Sigma_1,\Sigma_2)\simeq H^{Morse}_{n+1-*}(f)\simeq H_{n+1-*}(\overline{\Sigma}_1,\partial_-\overline{\Sigma}_1;\Z_2)\simeq H^*(\overline{\Sigma}_1,\partial_+\overline{\Sigma}_1;\Z_2)$.
 \end{enumerate}
\end{prop}

\section{Product structure}\label{PROD}
\subsection{Definition of the product}\label{defprod}

Let $\La_i^-\prec_{\Sigma_i}\La_i^+$, $i=1,2,3$, be three transverse exact Lagrangian cobordisms, and $T>0$ such that $\Sigma_i\backslash([-T,T]\times Y\cap\Sigma_i)$ are cylindrical. Recall that the moduli spaces which are useful to define the product are of different types. First, we need moduli spaces of pseudo-holomorphic curves with boundary on the negative cylindrical ends of the cobordisms and with two or three mixed asymptotics:
  \begin{alignat*}{1}
    &\cM_{\R\times\La^-_{i,j}}(\gamma_{j,i};\bs{\delta}_i,\gamma_{i,j},\bs{\delta}_i),\, i,j\in\{1,2,3\}\\
    &\cM_{\R\times\La^-_{1,2,3}}(\gamma_{3,1};\bs{\delta}_1,\gamma_{1,2},\bs{\delta}_2,\gamma_{2,3},\bs{\delta_3})\\
    &\cM_{\R\times\La^-_{1,2,3}}(\gamma_{3,1};\bs{\delta}_1,\gamma_{2,1},\bs{\delta}_2,\gamma_{2,3},\bs{\delta_3})\\
    &\cM_{\R\times\La^-_{1,2,3}}(\gamma_{3,1};\bs{\delta}_1,\gamma_{1,2},\bs{\delta}_2,\gamma_{3,2},\bs{\delta_3})\\
    &\cM_{\R\times\La^-_{1,2,3}}(\gamma_{3,1};\bs{\delta}_1,\gamma_{2,1},\bs{\delta}_{2},\gamma_{3,2},\bs{\delta}_3)
  \end{alignat*}
with $\gamma_{i,j}\in\Rc(\La_i^-,\La_j^-)$ for $1\leq i,j\leq3$.
Remark that the first one is a moduli space of bananas which already appeared in the definition of the Cthulhu differential. The four others are moduli spaces of pseudo-holomorphic curves having $\gamma_{3,1}$ as a positive Reeb chord asymptotic and the other mixed chords are positive or negative asymptotics depending on their direction.

Then we also need moduli spaces of pseudo-holomorphic curves with boundary in the non-cylindrical parts of the cobordisms, and again with two or three mixed asymptotics:
 \begin{alignat*}{1}
    &\cM_{\Sigma_{1,2}}(x^+;\bs{\delta}_1,\gamma_{2,1},\bs{\delta}_2)
 \end{alignat*}
for $x^+\in\Sigma_1\cap\Sigma_2$, and
 \begin{alignat*}{1}
    &\cM_{\Sigma_{1,2,3}}(x^+;\bs{\delta}_1,x_1,\bs{\delta}_2,x_2,\bs{\delta}_3)\\
    &\cM_{\Sigma_{1,2,3}}(x^+;\bs{\delta}_1,x_1,\bs{\delta}_2,\gamma_{3,2},\bs{\delta}_3)\\
    &\cM_{\Sigma_{1,2,3}}(x^+;\bs{\delta}_1,\gamma_{2,1},\bs{\delta}_2,x_2,\bs{\delta}_3)\\
    &\cM_{\Sigma_{1,2,3}}(x^+;\bs{\delta}_1,\gamma_{2,1},\bs{\delta}_2,\gamma_{3,2},\bs{\delta}_3)
 \end{alignat*}
for $x^+\in\Sigma_1\cap\Sigma_3$, and also
 \begin{alignat*}{1}
    &\cM_{\Sigma_{1,2,3}}(\gamma_{1,3};\bs{\delta}_1,x_1,\bs{\delta}_2,x_2,\bs{\delta}_3)\\
    &\cM_{\Sigma_{1,2,3}}(\gamma_{1,3};\bs{\delta}_1,\gamma_{2,1},\bs{\delta}_2,x_2,\bs{\delta}_3)\\
    &\cM_{\Sigma_{1,2,3}}(\gamma_{1,3};\bs{\delta}_1,x_1,\bs{\delta}_2,\gamma_{3,2},\bs{\delta}_3)
  \end{alignat*}
where $\gamma_{1,3}$ is a negative Reeb chord asymptotic.
  
We achieve transversality for these moduli spaces using domain dependent almost complex structures. First, remark that moduli spaces of curves with boundary on the negative cylindrical ends above are transversely cut out for a generic almost complex structure in $\Jc^{cyl}(\R\times Y)$. Indeed, even if some curves have several positive asymptotics, these are all distinct so it is always possible to find an injective point (same argument as for Legendrian contact homology type moduli spaces).

Now, consider $J^\pm\in\Jc^{cyl}(\R\times Y)$ regular almost complex structures for the moduli spaces of curves with boundary on the positive and negative cylindrical ends respectively, with two or three mixed Reeb chords asymptotics and negative pure Reeb chords asymptotics (in fact, we do not need here regularity of $J^+$ but in any case it will be useful in Section \ref{teo2}). We know that Cthulhu moduli spaces of curves with boundary on non-cylindrical parts of the cobordisms are transversely cut out for a generic time-dependent almost complex structure $J_t:[0,1]\to\Jc^{adm}_{J^-,J^+,T}(\R\times Y)$. So, let us denote by $J_t^{\Sigma_i,\Sigma_j}$, for $i,j\in\{1,2,3\}$, a regular time-dependent almost complex structure for Cthulhu moduli spaces associated to the pair of cobordisms $(\Sigma_i,\Sigma_j)$, with the convention that $J_t^{\Sigma_i,\Sigma_i}$ is a constant path. Then, given a consistent universal choice of strip-like ends, we use Seidel's result \cite[Section (9k)]{S} to deduce that a universal domain dependent almost complex structure $\Jc_{2,\underline{\Sigma}}:\Sc^3\to\Jc_{J^-,J^+,T}^{adm}(\R\times Y)$ can be perturb to a regular one for the moduli spaces above with boundary on non cylindrical parts and three mixed asymptotics.

\begin{figure}[ht]     
 	\begin{center}\includegraphics[width=7cm]{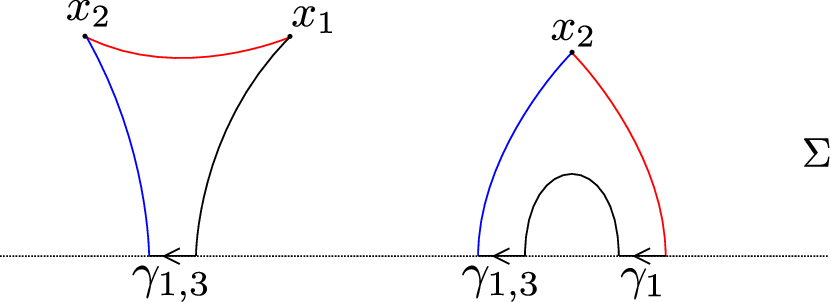}\end{center}
 	\caption{Curves contributing to $f^{(2)}(x_2,x_1)$ and $f^{(2)}(x_2,\gamma_1)$.}
 	\label{f2}
\end{figure}

\begin{figure}[ht]     
      \begin{center}\includegraphics[width=10cm]{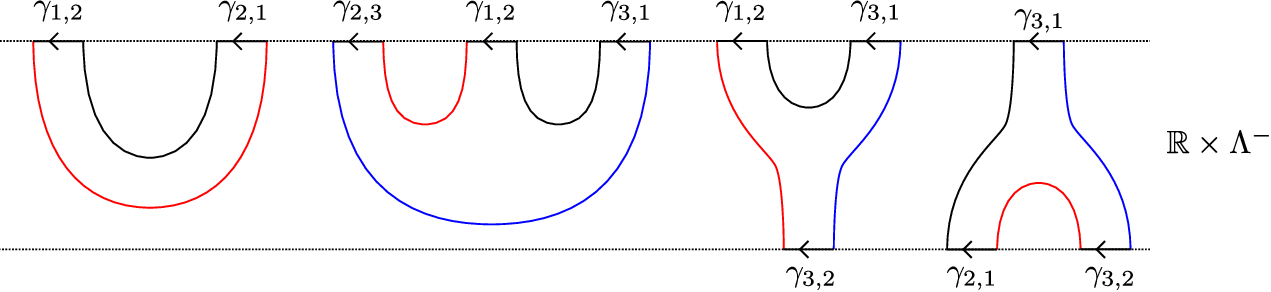}\end{center}
      \caption{Curves contributing to $b(\gamma_{1,2})$, $b^{(2)}(\gamma_{2,3},\gamma_{1,2})$, $b^{(2)}(\gamma_{3,2},\gamma_{1,2})$ and $b^{(2)}(\gamma_{3,2},\gamma_{2,1})$.}
      \label{les_bananes}
\end{figure}

This means that we can find a regular domain dependent almost complex structure with values in $\Jc^{adm}_{J^-,J^+,T}(\R\times Y)$ such that all the moduli spaces we have encountered until now are simultaneously smooth manifolds.
\begin{rem}
 In all the section, as before, we define maps by a count of rigid pseudo-holomorphic curves. This count will always be modulo $2$.
\end{rem}
Let us assume that the Chekanov-Eliashberg algebras $\Ac(\La^-_i)$, $i=1,2,3$, admit augmentations $\ep_i^-$. We want to define a map:
  \begin{alignat*}{1}
    & \mfm_2\colon CF_{-\infty}^*(\Sigma_2,\Sigma_3)\otimes CF_{-\infty}^*(\Sigma_1,\Sigma_2)\to CF_{-\infty}^*(\Sigma_1,\Sigma_3)\\
  \end{alignat*}
linear in each variable. This map decomposes as $\mfm_2=m^0+m^-$, where $m^0$ takes values in $CF^*(\Sigma_1,\Sigma_3)$ and $m^-$ takes values in $C^*(\La_1^-,\La_3^-)$. In order to do this, we define these maps on pairs of generators, which means that we must define the eight following components:
  \begin{alignat*}{3}
    & m_{00}^0 & \colon & CF^*(\Sigma_2,\Sigma_3)\otimes CF^*(\Sigma_1,\Sigma_2) \to CF^*(\Sigma_1,\Sigma_3)\\
    & m_{00}^- & \colon & CF^*(\Sigma_2,\Sigma_3)\otimes CF^*(\Sigma_1,\Sigma_2) \to C^*(\La_1^-,\La_3^-)\\
    & m_{0-}^0 & \colon & CF^*(\Sigma_2,\Sigma_3)\otimes C^*(\La_1^-,\La_2^-) \to CF^*(\Sigma_1,\Sigma_3)\\
    & m_{0-}^- & \colon & CF^*(\Sigma_2,\Sigma_3)\otimes C^*(\La_1^-,\La_2^-) \to C^*(\La_1^-,\La_3^-)\\
    & m_{-0}^0 & \colon & C^*(\La_2^-,\La_3^-)\otimes CF^*(\Sigma_1,\Sigma_2) \to CF^*(\Sigma_1,\Sigma_3)\\
    & m_{-0}^- & \colon & C^*(\La_2^-,\La_3^-)\otimes CF^*(\Sigma_1,\Sigma_2) \to C^*(\La_1^-,\La_3^-)\\
    & m_{--}^0 & \colon & C^*(\La_2^-,\La_3^-)\otimes C^*(\La_1^-,\La_2^-) \to CF^*(\Sigma_1,\Sigma_3)\\
    & m_{--}^- & \colon & C^*(\La_2^-,\La_3^-)\otimes C^*(\La_1^-,\La_2^-) \to C^*(\La_1^-,\La_3^-)\\
  \end{alignat*}
Let us begin by $m^0$. We set:
 \begin{alignat*}{1}
    & m^0_{00}(x_2,x_1)=\sum\limits_{x^+,\bs{\delta}_i}\#\cM^0_{\Sigma_{123}}(x^+;\bs{\delta}_1,x_1,\bs{\delta}_2,x_2,\bs{\delta}_3)\ep_1^-(\bs{\delta}_1)\ep_2^-(\bs{\delta}_2)\ep_3^-(\bs{\delta}_3)\cdot x^+\\
    & m^0_{0-}(x_2,\gamma_1)=\sum\limits_{x^+,\bs{\delta}_i}\#\cM^0_{\Sigma_{123}}(x^+;\bs{\delta}_1,\gamma_1,\bs{\delta}_2,x_2,\bs{\delta}_3)\ep_1^-(\bs{\delta}_1)\ep_2^-(\bs{\delta}_2)\ep_3^-(\bs{\delta}_3)\cdot x^+\\
    & m^0_{-0}(\gamma_2,x_1)=\sum\limits_{x^+,\bs{\delta}_i}\#\cM^0_{\Sigma_{123}}(x^+;\bs{\delta}_1,x_1,\bs{\delta}_2,\gamma_2,\bs{\delta}_3)\ep_1^-(\bs{\delta}_1)\ep_2^-(\bs{\delta}_2)\ep_3^-(\bs{\delta}_3)\cdot x^+\\
    & m^0_{--}(\gamma_2,\gamma_1)=\sum\limits_{x^+,\bs{\delta}_i}\#\cM^0_{\Sigma_{123}}(x^+;\bs{\delta}_1,\gamma_1,\bs{\delta}_2,\gamma_2,\bs{\delta}_3)\ep_1^-(\bs{\delta}_1)\ep_2^-(\bs{\delta}_2)\ep_3^-(\bs{\delta}_3)\cdot x^+
  \end{alignat*}
where the sums are for $x^+\in\Sigma_1\cap\Sigma_3$ and for each $i\in\{1,2,3\}$, $\bs{\delta}_i$ is a word of Reeb chords of $\La_i^-$. Then, to define $m^-$ we first introduce intermediate maps. We recall that there is a canonical identification of complexes $CF_{n+1-*}(\Sigma_b,\Sigma_a)=CF^*(\Sigma_a,\Sigma_b)$ and we denote $C_*(\La_a^-,\La_b^-)$ the dual of $C^*(\La_a^-,\La_b^-)$. We consider a map:
  \begin{alignat*}{1}
  f^{(2)}\colon CF_{-\infty}(\Sigma_2,\Sigma_3)\otimes CF_{-\infty}(\Sigma_1,\Sigma_2)\to C_{n-1-*}(\La^-_3,\La^-_1)
  \end{alignat*} 
defined on each pair of generators by: 
  \begin{alignat*}{1}
    & f^{(2)}(x_2,x_1)=\sum_{\gamma_{1,3},\bs{\delta}_i}\#\cM^0_{\Sigma_{1,2,3}}(\gamma_{1,3};\bs{\delta}_1,x_1,\bs{\delta}_2,x_2,\bs{\delta}_3)\ep_1^-(\bs{\delta}_1)\ep_2^-(\bs{\delta}_2)\ep_3^-(\bs{\delta}_3)\cdot\gamma_{1,3}\\
    & f^{(2)}(x_2,\gamma_1)=\sum_{\gamma_{1,3},\bs{\delta}_i}\#\cM^0_{\Sigma_{1,2,3}}(\gamma_{1,3};\bs{\delta}_1,\gamma_1,\bs{\delta}_2,x_2,\bs{\delta}_3)\ep_1^-(\bs{\delta}_1)\ep_2^-(\bs{\delta}_2)\ep_3^-(\bs{\delta}_3)\cdot\gamma_{1,3}\\
    & f^{(2)}(\gamma_2,x_1)=\sum_{\gamma_{1,3},\bs{\delta}_i}\#\cM^0_{\Sigma_{1,2,3}}(\gamma_{1,3};\bs{\delta}_1,x_1,\bs{\delta}_2,\gamma_2,\bs{\delta}_3)\ep_1^-(\bs{\delta}_1)\ep_2^-(\bs{\delta}_2)\ep_3^-(\bs{\delta}_3)\cdot\gamma_{1,3}\\
    & f^{(2)}(\gamma_2,\gamma_1)=0
  \end{alignat*}
where $x_2\in CF(\Sigma_2,\Sigma_3)$, $x_1\in CF(\Sigma_1,\Sigma_2)$, $\gamma_2\in \Rc(\La_3^-,\La_2^-)$ and $\gamma_1\in\Rc(\La_2^-,\La_1^-)$ (see Figure \ref{f2}). This map is the analogue of the map
  \begin{alignat*}{1}
     & \delta_{-0}^{\Sigma_{21}}:CF_{n+1-*}(\Sigma_2,\Sigma_1)=CF^*(\Sigma_1,\Sigma_2)\to C_{n-1-*}(\La_2^-,\La_1^-)\\
  \end{alignat*}
with three mixed asymptotics instead of two, where $\delta_{-0}^{\Sigma_{21}}$ is the dual of $d_{0-}^{\Sigma_{21}}:C^{*-2}(\La_2^-,\La_1^-)\to CF^*(\Sigma_2,\Sigma_1)$ (see Section \ref{cthulhu}). The Lagrangian label being given by the asymptotics, we will now denote by $f^{(1)}$ the maps $\delta_{-0}^{\Sigma_{21}}$ and $\delta_{-0}^{\Sigma_{32}}$, which we extend to the whole complex $CF_{-\infty}(\Sigma_i,\Sigma_j)$ by setting $f^{(1)}(\gamma_{j,i})=\gamma_{j,i}$ for a mixed Reeb chord.
Then we generalize the banana map $b$ with a map $b^{(2)}$ defined by a count of pseudo-holomorphic disks with three mixed asymptotics. Let us denote $\mathfrak{C}(\La_i^-,\La_j^-)=C^*(\La_i,\La_j)\oplus C_{n-1-*}(\La_j^-,\La_i^-)$, we define:
\begin{alignat*}{1}
 b^{(2)}\colon\mathfrak{C}(\La_2^-,\La_3^-)\otimes\mathfrak{C}(\La_1^-,\La_2^-)\to C^*(\La_1^-,\La_3^-)
\end{alignat*}
by
  \begin{alignat*}{1}
    & b^{(2)}(\gamma_{2,3},\gamma_{1,2})=\sum_{\gamma_{3,1},\bs{\delta}_i}\#\widetilde{\cM^1}_{\R\times\La^-_{123}}(\gamma_{3,1};\bs{\delta}_1,\gamma_{1,2},\bs{\delta}_2,\gamma_{2,3},\bs{\delta}_3)\ep_1^-(\bs{\delta}_1)\ep_2^-(\bs{\delta}_2)\ep_3^-(\bs{\delta}_3)\cdot\gamma_{3,1}\\
    & b^{(2)}(\gamma_{2,3},\gamma_{2,1})=\sum_{\gamma_{3,1},\bs{\delta}_i}\#\widetilde{\cM^1}_{\R\times\La^-_{123}}(\gamma_{3,1};\bs{\delta}_1,\gamma_{2,1},\bs{\delta}_2,\gamma_{2,3},\bs{\delta}_3)\ep_1^-(\bs{\delta}_1)\ep_2^-(\bs{\delta}_2)\ep_3^-(\bs{\delta}_3)\cdot\gamma_{3,1}\\
    & b^{(2)}(\gamma_{3,2},\gamma_{1,2})=\sum_{\gamma_{3,1},\bs{\delta}_i}\#\widetilde{\cM^1}_{\R\times\La^-_{123}}(\gamma_{3,1};\bs{\delta}_1,\gamma_{1,2},\bs{\delta}_2,\gamma_{3,2},\bs{\delta}_3)\ep_1^-(\bs{\delta}_1)\ep_2^-(\bs{\delta}_2)\ep_3^-(\bs{\delta}_3)\cdot\gamma_{3,1}\\
    & b^{(2)}(\gamma_{3,2},\gamma_{2,1})=\sum\limits_{\gamma_{3,1},\bs{\delta}_i}\#\widetilde{\cM^1}_{\R\times\La^-_{123}}(\gamma_{3,1};\bs{\delta}_1,\gamma_{2,1},\bs{\delta}_2,\gamma_{3,2},\bs{\delta}_3)\ep_1^-(\bs{\delta}_1)\ep_2^-(\bs{\delta}_2)\ep_3^-(\bs{\delta}_3)\cdot\gamma_{3,1}
  \end{alignat*}
with $\gamma_{i,j}\in\Rc(\La_i^-,\La_j^-)$ and words of Reeb chords $\bs{\delta}_i$ of $\La_i^-$. Figure \ref{les_bananes} shows examples of curves counted by $b^{(2)}$.
\begin{rem}
The map $b^{(2)}$ restricted to $C^*(\La_2^-,\La_3^-)\otimes C^*(\La_1^-,\La_2^-)$ is in fact equal to the product $\mu^2_{\ep^-_{3,2,1}}$ in the augmentation category $\Aug_-(\La_1^-\cup\La_2^-\cup\La_3^-)$ restricted to this sub-complex, where $\ep_{3,2,1}^-$ is the diagonal augmentation of $\Ac(\La_1^-\cup\La_2^-\cup\La_3^-)$ built from $\ep_1^-,\ep_2^-$ and $\ep_3^-$ (Section \ref{linearisation}).
\end{rem}
Now we define the map $m^-$ by the following formula. For a pair $(a_2,a_1)\in CF_{\infty}(\Sigma_2,\Sigma_3)\otimes CF_{-\infty}(\Sigma_1,\Sigma_2)$, we set:
  \begin{alignat*}{1}
   m^-(a_2,a_1)=b\circ f^{(2)}(a_2,a_1)+b^{(2)}(f^{(1)}(a_2),f^{(1)}(a_1))
  \end{alignat*}
More precisely, for each pair of asymptotics, we have:  \begin{alignat*}{1}
    & m^-_{00}(x_2,x_1)=b\circ f^{(2)}(x_2,x_1)+b^{(2)}(f^{(1)}(x_2),f^{(1)}(x_1))\\
    & m^-_{0-}(x_2,\gamma_1)=b\circ f^{(2)}(x_2,\gamma_1)+b^{(2)}(f^{(1)}(x_2),\gamma_1)\\
    & m^-_{-0}(\gamma_2,x_1)=b\circ f^{(2)}(\gamma_2,x_1)+b^{(2)}(\gamma_2,f^{(1)}(x_1))\\
    & m^-_{--}(\gamma_2,\gamma_1)=b^{(2)}(\gamma_2,\gamma_1)
  \end{alignat*}
Contrary to the definition of $m^0$, when at least one input is an intersection point we need to count broken curves instead of just one type of pseudo-holomorphic disk, in order to associate a positive Reeb chord in $C^*(\La_1^-,\La_3^-)$ to the two inputs. These broken curves have two levels, one level contains curves with boundary on the non cylindrical parts of the cobordisms, and the other level contains a curve with boundary on the negative cylindrical ends of the cobordisms. These configurations look like pseudo-holomorphic buildings but they are not because their components cannot be glued, so in particular, they are rigid. These configurations are part of what we will call \textit{unfinished pseudo-holomorphic buildings}, for which we now give a definition in a general setting. This definition is very close to the definition of pseudo-holomorphic building given in Section \ref{comp}. The main difference is that after boundary connected sum of the components at nodes, the remaining asymptotics contain a positive Reeb chord in the bottom level and so the gluing results do not apply.
As in Section \ref{comp}, consider $d+1$ transverse exact Lagrangian cobordisms $\La_i^-\prec_{\Sigma_i}\La_i^+$, and the following Lagrangian labels $\underline{\Sigma}=(\Sigma_1,\dots,\Sigma_{d+1})$ and $\underline{\R\times\Lambda}^\pm=(\R\times\La_1^\pm,\dots,\R\times\La_{d+1}^\pm)$.
Recall that given a planar rooted tree $T$ with $d-1$ leaves, $d\geq2$, we associate to each interior vertex $v$ a triple $(S_{r_v},I_v,L_v)$. Let us distinguish this time one vertex that we denote $v_0$ which is the unique interior vertex connected to the root of $T$.
Again, the set of all nodes $\bigcup\limits_{v} I_v$ contains an even number of elements organized in pair and we denote
\begin{alignat*}{1}
\bigcup\limits_{v} I_v=\{p_1,\bar{p}_1\}\cup\{p_2,\bar{p}_2\}\cup\dots\cup\{p_k,\bar{p}_k\}
\end{alignat*}
the partition of $\bigcup I_v$ in such pairs, where $k$ is the number of interior edges of $T$.

\begin{defi}\label{unfb}
	Given integers $k^-,k^+\geq0$, an \textit{unfinished pseudo-holomorphic building of height $k^-|1|k^+$} in $\R\times Y$ with boundary on $\underline{\Sigma}$ is the data of:
	\begin{enumeratec}
		\item a planar rooted tree $T$ and the corresponding union of triples $\bigcup\limits_{v}(S_{r_v},I_v,L_v)$ and the distinguished vertex $v_0$,
		\item a choice of asymptotic in $A(\underline{\Sigma})$ for each node in $\bigcup I_v$. We require that for each pair $\{p_j,\bar{p}_j\}$, the same asymptotic is assigned to $p_j$ and $\bar{p}_j$,
		\item a choice of asymptotic for each marked point in $\bigcup L_v$,
		\item a pair $(u_v,\rho_v)$ for each interior vertex $v\neq v_0$ of $T$, where $u_v:S_{r_v}\to\R\times Y$ is a pseudo-holomorphic disk asymptotic to the given asymptotics assigned to elements in $I_v\cup L_v$, and $\rho_v$ is the floor of $v$, satisfying $-k^-\leq\rho_v\leq k^+$,
		\item a pair $(u_{v_0},\rho_{v_0})$ where $u_{v_0}:S_{r_{v_0}}\to\R\times Y$ is a pseudo-holomorphic disk asymptotic to the given asymptotics assigned to elements in $I_{v_0}\cup L_{v_0}$, in the bottom level (i.e. $\rho_{v_0}<0$), so in particular $u_{v_0}$ is only asymptotic to Reeb chords (no intersection points).
	
		These data are required to satisfy some conditions. First, let us denote $v_1,\dots,v_j$ the interior vertices of $T$ connected to $v_0$ respectively by edges that we denote $e_1,\dots, e_j$. Add a vertex $l_{e_i}$ on each edge $e_i$, in particular we have $|l_{e_i}|=2$. Denote $\widehat{T}$ the planar rooted tree obtained from $T$ by adding these vertices. Consider now the subtrees $T_1,\dots,T_j$ of $\widehat{T}$ rooted at $l_{e_i}$ and containing all the descendants of $l_{e_i}$ (under the canonical orientation from the root to the leaves given by a rooted tree). Remark that then $T_0=\widehat{T}\backslash\bigcup_jT_j$ is a planar rooted tree with one interior vertex, $v_0$, leaves $l_{e_1},\dots,l_{e_j}$, a root $l_0$ (the original root of $T$) and potentially other leaves of the original tree $T$ that we denote $l_1,\dots,l_s$. The numbering is not significant here. Remark that each $l_{e_i}$ corresponds to a pair of nodes of $T$ that we denote without loss of generality $\{p_i,\bar{p}_i\}$, and we assume $p_i\in I_{v_i}$ and $\bar{p}_i\in I_{v_0}$. The following conditions are required:
		\begin{enumeratec}
			\item each floor $-k^-\leq\rho\leq k^+$ admits at least one non trivial disk,
			\item each tree $T_i$ associated with the given asymptotic data and pseudo-holomorphic disks chosen above gives a pseudo-holomorphic building,
			\item $S_{r_{v_0}}$ has boundary punctures at $\bar{p}_i$, $1\leq i\leq j$, and boundary punctures corresponding to $l_0$ and $l_1,\dots,l_s$ that we still denote $l_0,\dots,l_s\in L_{v_0}$ by abuse of notation, such that
			\begin{enumeratec}
			\item $u_{v_0}$ has a positive Reeb chord asymptotic at $l_0$ and negative Reeb chord asymptotics at $l_i$ for $1\leq  i\leq s$,
			\item if $u_{v_i}$ has a positive (resp. negative) asymptotic to a Reeb chord $\gamma\in A(\underline{\Sigma})$ at the node $p_i$, then $u_{v_0}$ must have a negative (resp. positive) asymptotic to $\gamma$ at $\bar{p}_i$, and we have $\rho_{v_i}=\rho_{v_0}-1$ (resp. $\rho_{v_i}=\rho_{v_0}+1$).
			\end{enumeratec}		
		\end{enumeratec}
	\end{enumeratec}
\end{defi}

\begin{rem}
	Although the previous definition is quite long and complicated, remark that an unfinished pseudo-holomorphic building is an object which satisfies the conditions of Definition \ref{phb} except the condition (g), and has the following additional properties:
	\begin{itemize}
		\item the middle level is non-empty,
		\item the component $u_{v_0}$ lives in the bottom level and has exactly one positive Reeb chord asymptotic which is not a node,
		\item the other pseudo-holomorphic disks defining the building (but not $u_{v_0}$) satisfy the condition (g) of Definition \ref{phb}.
	\end{itemize}
\end{rem}

Next, we define the notion of equivalence of unfinished pseudo-holomorphic buildings. The definition is actually the same as Definition \ref{def:eqphb} for pseudo-holomorphic buildings.
\begin{defi}[Equivalence of unfinished pseudo-holomorphic buildings]\label{defequi}
	Two unfinished pseudo-holomorphic buildings are \textit{equivalent} if they become the same after the removal of an appropriate number of trivial cylinders together with the obvious deformation of the underlying planar rooted tree, to each of them in the bottom and top levels. 
	In other words, two unfinished pseudo-holomorphic buildings are equivalent if they become the same after:
		\begin{itemize}
			\item removing all possible trivial cylinders in the bottom and top levels attached to asymptotics assigned to $\cup L_v$,
			\item removing simultaneously trivial cylinders attached to all the positive ends and/or all the negative ends corresponding to nodes in $\cup I_v$ of a component in the bottom or top level.
		\end{itemize}
\end{defi}

\begin{figure}[ht]
	\begin{center}\includegraphics[width=13cm]{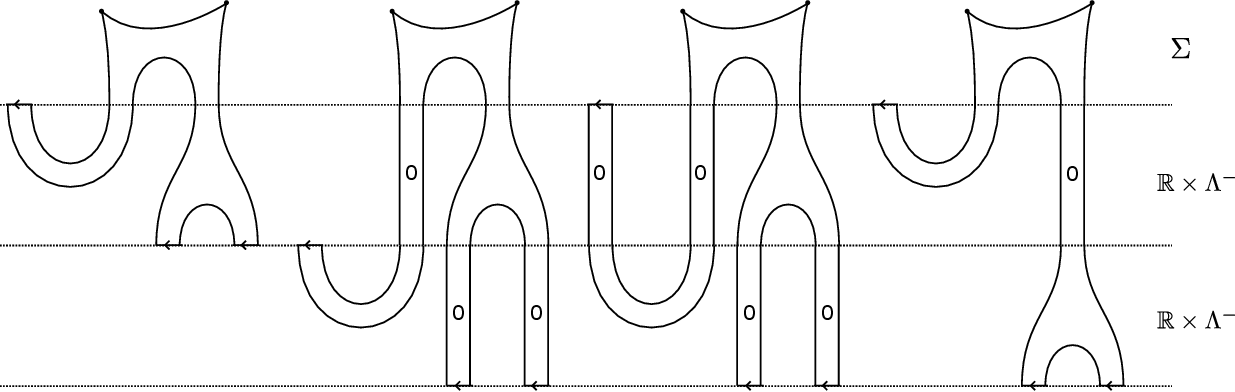}\end{center}
	\caption{Four equivalent unfinished pseudo-holomorphic buildings, where the "0" indicates a trivial strip.}
	\label{equivalence}
\end{figure}

See Figure \ref{equivalence} for an example of $4$ equivalent unfinished pseudo-holomorphic buildings.
\begin{rem}
	The trees associated to equivalent unfinished pseudo-holomorphic buildings are the same up to the addition/removal of vertices of valency $2$.
\end{rem}
\begin{rem}
	Observe that condition (5)(a) in Definition \ref{unfb} implies that we do not consider any unfinished building with boundary only on the negative ends $\underline{\R\times\La^-}$, because such an unfinished building would actually be equivalent to a pseudo-holomorphic building as defined in Section \ref{comp}. 
\end{rem}

\begin{figure}    
      \begin{center}\includegraphics[width=8cm]{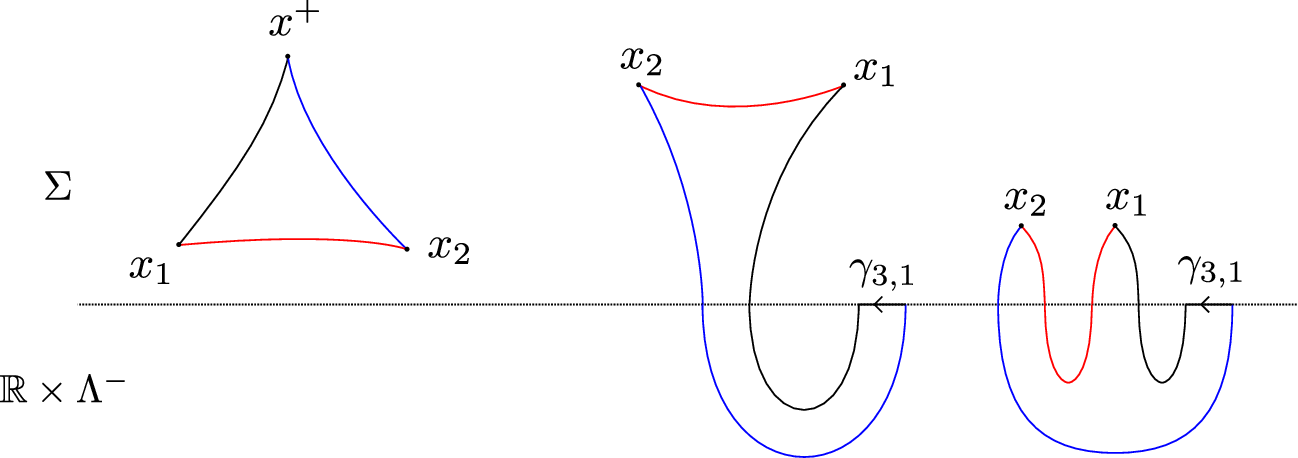}\end{center}
      \caption{Left: curve contributing to $m^0_{00}(x_2,x_1)$; right: curves contributing to $m^-_{00}(x_2,x_1)$.}
      \label{produitxx}
\end{figure}

\begin{figure}[ht]     
      \begin{center}\includegraphics[width=8cm]{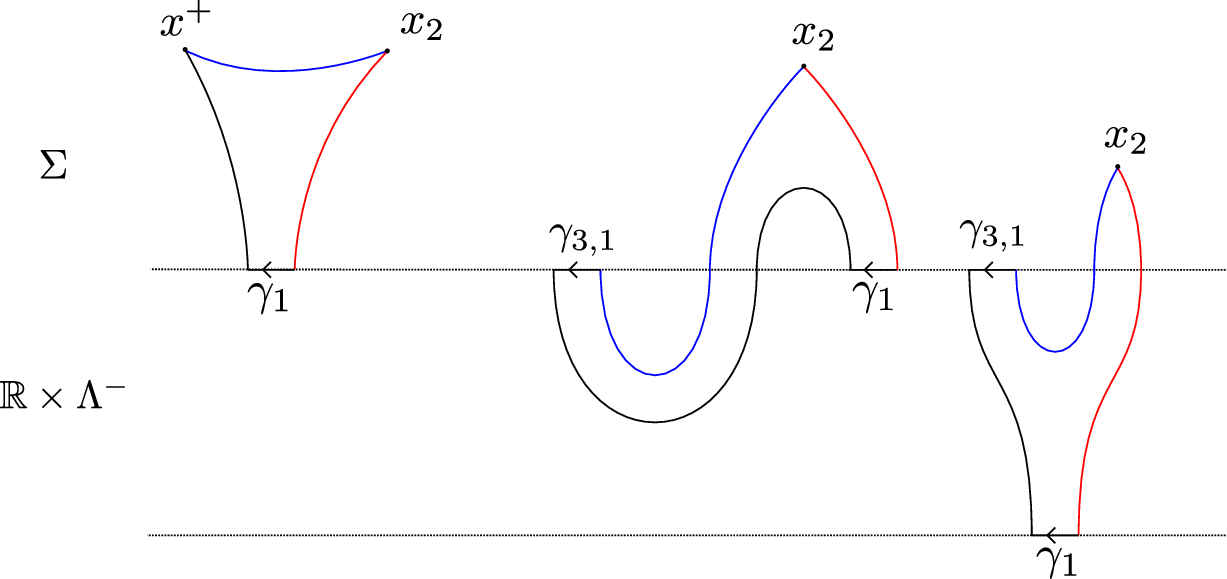}\end{center}
      \caption{Left: curve contributing to $m^0_{0-}(x_2,\gamma_1)$; right: curves contributing to $m^-_{0-}(x_2,\gamma_1)$}
      \label{produitxc}
\end{figure}

\begin{figure}[ht]     
      \begin{center}\includegraphics[width=8cm]{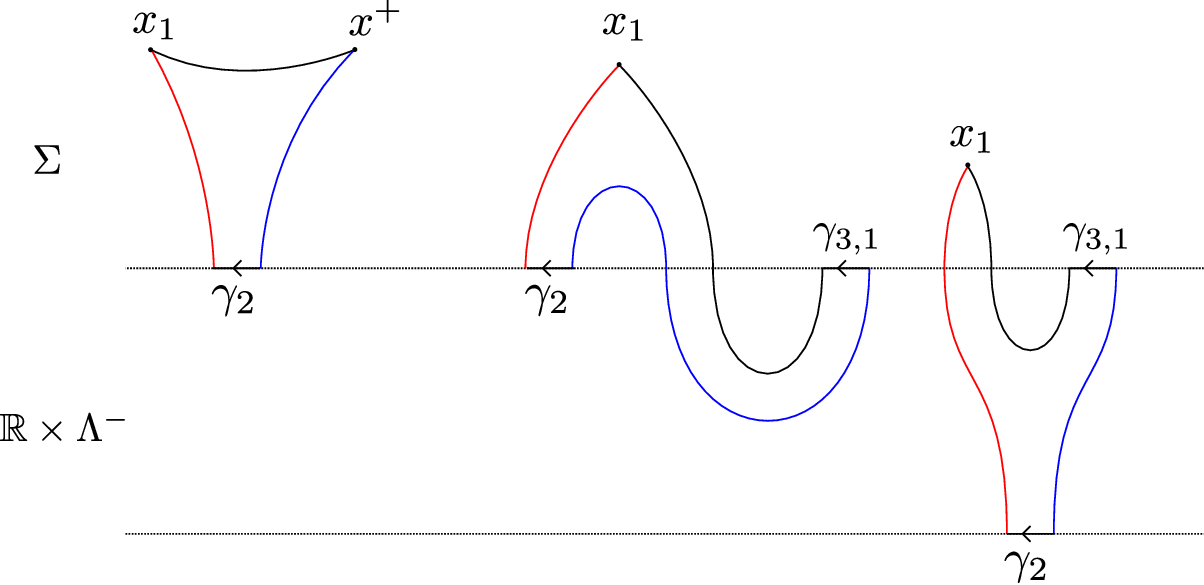}\end{center}
      \caption{Left: curve contributing to $m^0_{-0}(\gamma_2,x_1)$; right: curves contributing to $m^-_{-0}(\gamma_2,x_1)$}
      \label{produitcx}
\end{figure}

\begin{figure}[ht]     
      \begin{center}\includegraphics[width=7cm]{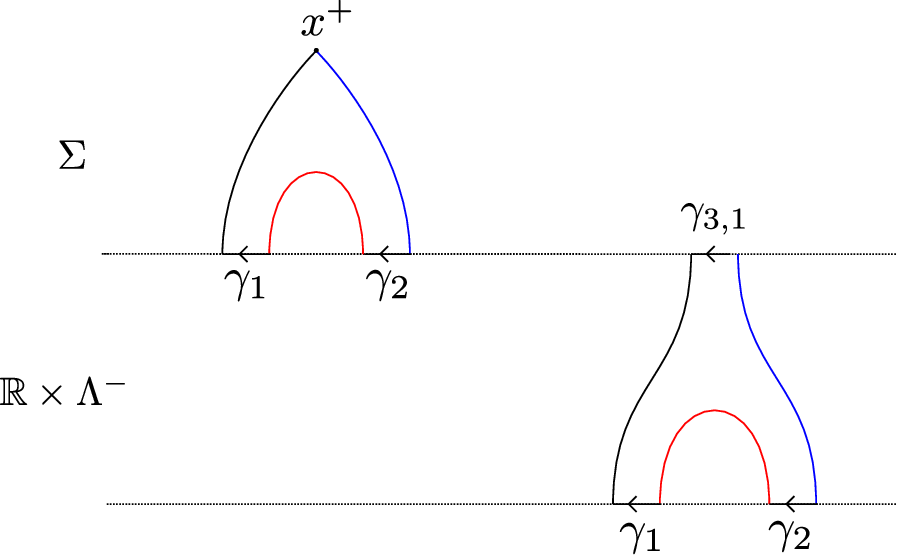}\end{center}
      \caption{Left: curve contributing to $m^0_{--}(\gamma_2,\gamma_1)$; right: curve contributing to $m^-_{--}(\gamma_2,\gamma_1)$.}
      \label{produitcc}
\end{figure}

Roughly speaking, one can imagine unfinished buildings as being buildings where we have removed some component in the middle level in such a way that it is not a building anymore.
Thus the map $m^-$ counts unfinished pseudo-holomorphic buildings. On Figures \ref{produitxx}, \ref{produitxc}, \ref{produitcx} and \ref{produitcc} are schematized the different types of curves and unfinished buildings that contribute to $\mfm_2$.
\begin{rem}
 By \cite[Proposition 3.2]{CDGG2} for curves with boundary on three transverse exact Lagrangian cobordisms instead of two, we can express the dimension of the moduli spaces involved in the definition of $\mfm_2$ by the degree of the asymptotics. Then it is not hard to check that $\mfm_2$ is a degree $0$ map, with the shift in grading for Reeb chords (see Section \ref{CTH}).
\end{rem}

\subsection{Proof of Theorem \ref{teo_prod}}\label{proof}
In this section, we prove that $\mfm_2$ satisfies the Leibniz rule:
\begin{alignat*}{1}
 \mfm_2(-,\partial_{-\infty})+\mfm_2(\partial_{-\infty},-)+\partial_{-\infty}\circ\mfm_2(-,-)=0
\end{alignat*}
In order to do this, we show that the above relation is satisfied for each pair of generators in $CF_{-\infty}(\Sigma_2,\Sigma_3)\otimes CF_{-\infty}(\Sigma_1,\Sigma_2)$. For example, for $(x_2,x_1)\in CF(\Sigma_2,\Sigma_3)\otimes CF(\Sigma_1,\Sigma_2)$, this gives:
\begin{alignat*}{1}
 &\mfm_2\big(x_2,\partial_{-\infty}(x_1)\big)+\mfm_2\big(\partial_{-\infty}(x_2),x_1\big)+\partial_{-\infty}\circ\mfm_2(x_2,x_1)=0\\
 &\Leftrightarrow\,\mfm_2\big(x_2,(d_{00}+d_{-0})(x_1)\big)+\mfm_2\big((d_{00}+d_{-0})(x_2),x_1\big)\\
 &\hspace{43mm}+(d_{00}+d_{-0})\circ m^0(x_2,x_1)+(d_{0-}+d_{--})\circ m^-(x_2,x_1)=0\\
 &\Leftrightarrow\,\Big(m^0\big(x_2,(d_{00}+d_{-0})(x_1)\big)+m^0\big((d_{00}+d_{-0})(x_2),x_1\big)+d_{00}\circ m^0(x_2,x_1)+d_{0-}\circ m^-(x_2,x_1)\Big)\\
 &\hspace{45mm}+\Big(m^-\big(x_2,(d_{00}+d_{-0})(x_1)\big)+m^-\big((d_{00}+d_{-0})(x_2),x_1\big)\\
 &\hspace{45mm}+d_{-0}\circ m^0(x_2,x_1)+d_{--}\circ m^-(x_2,x_1)\Big)=0\\
\end{alignat*}
and in the last equality the two terms into big brackets must vanish because the first one is an element in $CF(\Sigma_1,\Sigma_3)$ and the second one is an element in $C(\La_1^-,\La_3^-)$. Thus, considering each pair of generators we obtain in total eight relations to prove which are the following.
\vspace{2mm}	

\noindent1.\textit{ For a pair $(x_2,x_1)\in CF(\Sigma_2,\Sigma_3)\otimes CF(\Sigma_1,\Sigma_2)$}:
       \begin{alignat}{3}
        &\bullet\,\,m_{00}^0(x_2,d_{00}(x_1))+m_{00}^0(d_{00}(x_2),x_1)+d_{00}\circ m_{00}^0(x_2,x_1)\label{rel0xx}  \\
	&\hspace{3cm}+m_{0-}^0(x_2,d_{-0}(x_1))+m_{-0}^0(d_{-0}(x_2),x_1)+d_{0-}\circ m_{00}^-(x_2,x_1)=0\nonumber\\
	&\bullet\,\, m_{00}^-(x_2,d_{00}(x_1))+m_{00}^-(d_{00}(x_2),x_1)+m_{0-}^-(x_2,d_{-0}(x_1))\label{rel-xx} \\
	&\hspace{3cm}+m_{-0}^-(d_{-0}(x_2),x_1)+d_{-0}\circ m_{00}^0(x_2,x_1)+d_{--}\circ m_{00}^-(x_2,x_1)=0\nonumber
       \end{alignat}
\noindent2.\textit{ For a pair $(x_2,\gamma_1)\in CF(\Sigma_2,\Sigma_3)\otimes C(\La_1^-\La_2^-)$}:
      \begin{alignat}{1}
       &\bullet\,\,m_{00}^0(x_2,d_{0-}(\gamma_1))+m_{0-}^0(d_{00}(x_2),\gamma_1)+d_{00}\circ m_{0-}^0(x_2,\gamma_1)\label{rel0xc}\\
       &\hspace{3cm}+m_{0-}^0(x_2,d_{--}(\gamma_1))+m_{--}^0(d_{-0}(x_2),\gamma_1)+d_{0-}\circ m_{0-}^-(x_2,\gamma_1)=0\nonumber \\
       &\bullet\,\, m_{00}^-(x_2,d_{0-}(\gamma_1))+m_{0-}^-(d_{00}(x_2),\gamma_1)+m_{0-}^-(x_2,d_{--}(\gamma_1))\label{rel-xc}\\
       &\hspace{3cm}+m_{--}^-(d_{-0}(x_2),\gamma_1)+d_{-0}\circ m_{0-}^0(x_2,\gamma_1)+d_{--}\circ m_{0-}^-(x_2,\gamma_1)=0\nonumber
      \end{alignat}
\noindent3.\textit{ For a pair $(\gamma_2,x_1)\in C(\La_2^-,\La_3^-)\otimes CF(\Sigma_1,\Sigma_2)$}:
      \begin{alignat}{1}
       &\bullet\,\, m_{-0}^0(\gamma_2,d_{00}(x_1))+m_{00}^0(d_{0-}(\gamma_2),x_1)+d_{00}\circ m_{-0}^0(\gamma_2,x_1)\label{rel0cx} \\
       &\hspace{3cm}+m_{--}^0(\gamma_2,d_{-0}(x_1))+m_{-0}^0(d_{--}(\gamma_2),x_1)+d_{0-}\circ m_{-0}^-(\gamma_2,x_1)=0\nonumber \\
       &\bullet\,\, m_{-0}^-(\gamma_2,d_{00}(x_1))+m_{00}^-(d_{0-}(\gamma_2),x_1)+m_{--}^-(\gamma_2,d_{-0}(x_1)) \label{rel-cx}\\
       &\hspace{3cm}+m_{-0}^-(d_{--}(\gamma_2),x_1)+d_{-0}\circ m_{-0}^0(\gamma_2,x_1)+d_{--}\circ m_{-0}^-(\gamma_2,x_1)=0\nonumber
      \end{alignat}
\noindent4.\textit{ For a pair $(\gamma_2,\gamma_1)\in C(\La_2^-,\La_3^-)\otimes C(\La_1^-,\La_2^-)$}:
      \begin{alignat}{1}
       &\bullet\,\, m_{-0}^0(\gamma_2,d_{0-}(\gamma_1))+m_{0-}^0(d_{0-}(\gamma_2),\gamma_1)+d_{00}\circ m_{--}^0(\gamma_2,\gamma_1)\label{rel0cc}\\
       &\hspace{31mm}+m_{--}^0(\gamma_2,d_{--}(\gamma_1))+m_{--}^0(d_{--}(\gamma_2),\gamma_1)+d_{0-}\circ m_{--}^-(\gamma_2,\gamma_1)=0\nonumber\\
       &\bullet\,\, m_{--}^-(\gamma_2,d_{--}(\gamma_1))+m_{--}^-(d_{--}(\gamma_2),\gamma_1)+d_{--}\circ m_{--}^-(\gamma_2,\gamma_1)=0 \label{rel-cc} 
      \end{alignat}
\vspace{1mm}

To obtain these relations, we study the different types of pseudo-holomorphic buildings involved in the definition of each term appearing in the relations. Each curve in these buildings are rigid because the Cthulhu differential and the map $\mfm_2$ are defined by a count of rigid configurations. This means that the curves are of index $0$ if their boundary is on non-cylindrical Lagrangians, and of index $1$ if their boundary is on the negative cylindrical ends of the cobordisms. Compactness and gluing results imply that these broken curves are in bijection with elements in the boundary of the compactification of some moduli spaces. We recall below some properties that must be satisfied by the pseudo-holomorphic buildings we will consider here:
\begin{enumerate}
 \item each curve in a pseudo-holomorphic building must have positive energy, so for example each component with only Reeb chords asymptotics must have at least one positive Reeb chord asymptotic. For curves with also intersection points asymptotics, as the action is independent of the label, it is possible to have curves with only negative action asymptotics, but in any case the energy must be positive (see Section \ref{action_energie} and Subsection \ref{Rel-xc} below),
 \item each curve has a non negative Fredholm index because of the regularity of the almost complex structure, 
 \item the following relation on indices must be satisfied: if $u_1,\dots,u_k$ are curves forming a pseudo-holomorphic building, the glued solution $u$ has index given by $\ind(u)=\nu+\sum_i\ind(u_i)$, where $\nu$ is the number of pair of nodes asymptotic to intersection points (Section \ref{comp}).
\end{enumerate} 

\begin{figure}[ht]    
     \begin{center}\includegraphics[width=2cm]{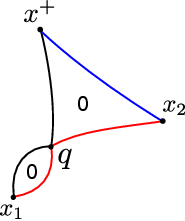}\end{center}
      \caption{Pseudo-holomorphic building contributing to $m_{00}^0(x_2,d_{00}(x_1))$.}
     \label{recolle1}
\end{figure}

\subsubsection{Relation \eqref{rel0xx}}
The first term appearing in this relation is $m_{00}^0(x_2,d_{00}(x_1))$. For every intersection point $x^+\in\Sigma_1\cap\Sigma_3$, the coefficient $\langle m_{00}^0(x_2,d_{00}(x_1)),x^+\rangle$ is defined by a count of pseudo-holomorphic buildings whose components are two index-$0$ curves with boundary on $\Sigma_1\cup\Sigma_2\cup\Sigma_3$, which have a common asymptotic at an intersection point $q\in CF(\Sigma_1,\Sigma_2)$. One curve contributes to $\langle d_{00}(x_1),q\rangle$ and the other contributes to $\langle m_{00}^0(x_2,q),x^+\rangle$ (see Figure \ref{recolle1}, where the numbers in the curves indicate the Fredholm index).

\begin{figure}[ht]    
     \begin{center}\includegraphics[width=145mm]{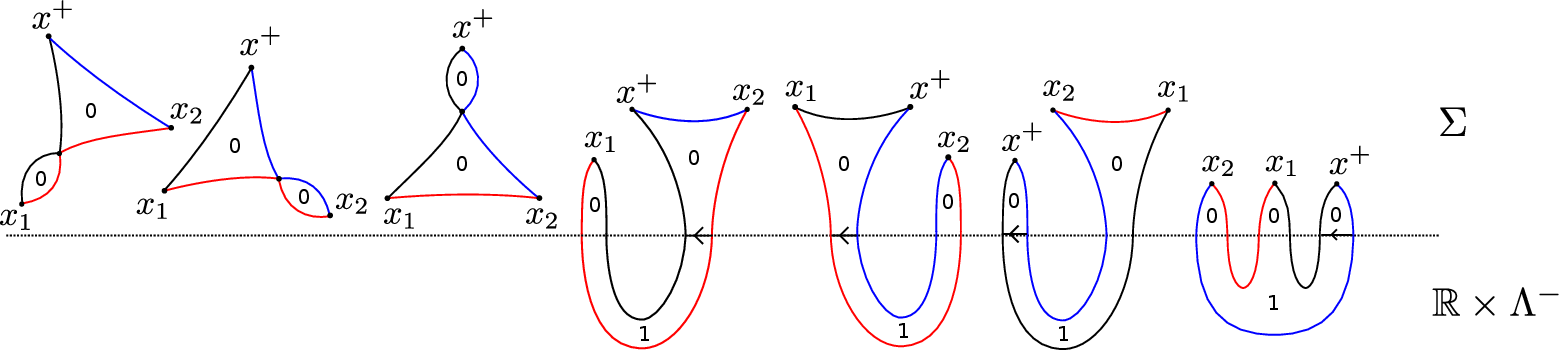}\end{center}
      \caption{Pseudo-holomorphic buildings in the boundary of the compactification of $\cM^1(x^+;\bs{\delta}_1^-,x_1,\bs{\delta}_2^-,x_2,\bs{\delta}_3^-)$.}
     \label{degexx1}
\end{figure}

The two curves can be glued together along $q$ and the resulting curve is an index-$1$ curve in the moduli space $\cM^1_{\Sigma_{123}}(x^+;\bs{\delta}_1,x_1,\bs{\delta}_2,x_2,\bs{\delta}_3)$. This implies that the holomorphic buildings contributing to $m_{00}^0(x_2,d_{00}(x_1))$ are in the boundary of the compactification of this moduli space. In fact, each term of the Relation \eqref{rel0xx} is defined by a count of holomorphic buildings whose components can be glued to give a curve in $\cM^1_{\Sigma_{123}}(x^+;\bs{\delta}_1,x_1,\bs{\delta}_2,x_2,\bs{\delta}_3)$. Thus now we look at all the possible breakings that can occur for a one parameter family of curves in this dimension $1$ moduli space. The curve can break on:
\begin{enumerate}
 \item an intersection point in $\Sigma_1\cap\Sigma_2$, $\Sigma_2\cap\Sigma_3$, or $\Sigma_3\cap\Sigma_1$, giving a pseudo-holomorphic building with one level containing two curves with a common asymptotic at this intersection point,
 \item a Reeb chord, giving a building of height $1|1|0$, the middle level containing index $0$ curves with boundary on $\overline{\Sigma}_1\cup\overline{\Sigma}_2\cup\overline{\Sigma}_3$, the bottom level containing an index-$1$ curve with boundary on $\R\times(\La^-_1\cup\La_2^-\cup\La_3^-)$.
\end{enumerate}
\begin{figure}[ht]    
     \begin{center}\includegraphics[width=11cm]{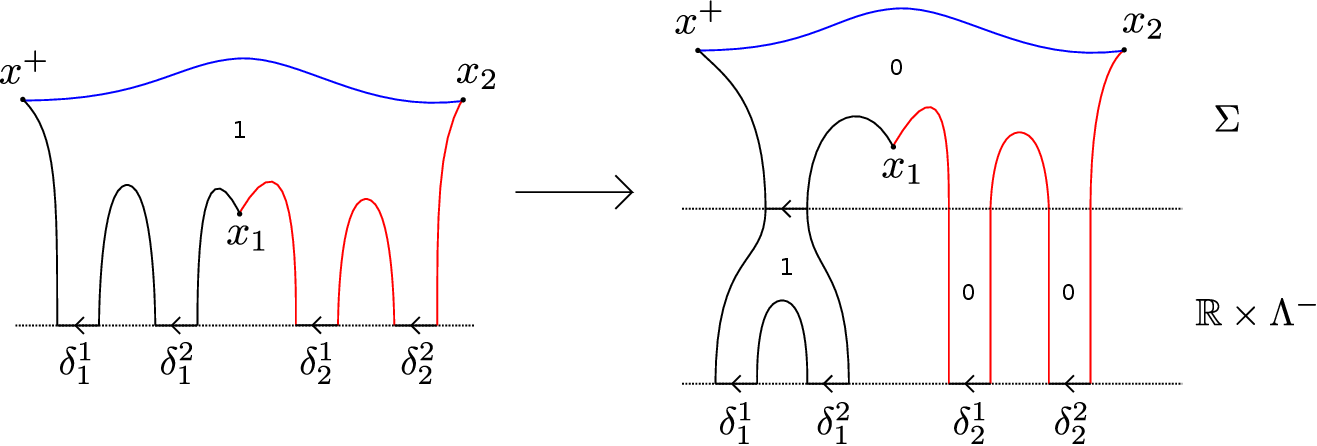}\end{center}
      \caption{$\partial$-breaking of a curve in $\cM^1(x^+;\bs{\delta}_1,x_1,\bs{\delta}_2,x_2)$.}
     \label{deltabreak}
\end{figure}
\begin{rem}
 In the second case, if the curve breaks on a pure Reeb chord $\gamma\in\Rc(\La_i^-)$ for $i\in\{1,2,3\}$, this is called a \textit{$\partial$-breaking} (see Figure \ref{deltabreak}). One component of such a broken curve contributes to $\partial^i(\gamma)$, where $\partial^i$ is the differential of the Legendrian contact homology of $\La_i^-$. We denote by $\overline{\cM}^\partial(x^+;\bs{\delta}_1,x_1,\bs{\delta}_2,x_2,\bs{\delta}_3)$ the union of all the $\partial$-breakings obtained as degeneration of curves in $\cM^1_{\Sigma_{123}}(x^+;\bs{\delta}_1,x_1,\bs{\delta}_2,x_2,\bs{\delta}_3)$. Now, the Cthulhu differential and the maps involved in the definition of the product $\mfm_2$ are defined by a count of elements in some moduli spaces of curves with two or three mixed asymptotics and every possible words of pure Reeb chords asymptotics $\bs{\delta}_i$. Thus, the $\partial$-breakings on a chord $\gamma$ for every possible words of pure chords $\bs{\delta}_i$ will contain all the curves contributing to $\partial_i(\gamma)$. Then, in the definition of the Cthulhu differential and the product, pure chords are augmented by $\ep_i^-$ and by definition $\ep_i^-\circ\partial^i=0$, so this means that the total contribution of $\partial$-breakings vanishes. 
\end{rem}
The boundary of the compactification of $\cM^1_{\Sigma_{123}}(x^+;\bs{\delta}_1,x_1,\bs{\delta}_2,x_2,\bs{\delta}_3)$ can be decomposed as follows:
\begin{alignat*}{2}
  & \partial\overline{\cM^1}_{\Sigma_{123}}(x^+;\bs{\delta}_1,x_1,\bs{\delta}_2,x_2,\bs{\delta}_3)=\overline{\cM}^\partial(x^+,\bs{\delta}_1,x_1,\bs{\delta}_2,x_2,\bs{\delta}_3)\\
  & \bigcup_{\substack{p\in\Sigma_1\cap\Sigma_2\\ \bs{\delta}_1'\bs{\delta}_1''=\bs{\delta}_1,\bs{\delta}_2'\bs{\delta}_2''=\bs{\delta}_2}}\cM(x^+;\bs{\delta}_1',p,\bs{\delta}_2'',x_2,\bs{\delta}_3)\times\cM_{\Sigma_{12}}(p;\bs{\delta}_1'',x_1,\bs{\delta}_2')\\
  & \bigcup_{\substack{q\in\Sigma_2\cap\Sigma_3\\ \bs{\delta}_2'\bs{\delta}_2''=\bs{\delta}_2,\bs{\delta}_3'\bs{\delta}_3''=\bs{\delta}_3}}\cM(x^+;\bs{\delta}_1,x_1,\bs{\delta}_2',q,\bs{\delta}_3'')\times\cM_{\Sigma_{23}}(q;\bs{\delta}_2'',x_2,\bs{\delta}_3')\\
  & \bigcup_{\substack{r\in\Sigma_1\cap\Sigma_3\\ \bs{\delta}_1'\bs{\delta}_1''=\bs{\delta}_1,\bs{\delta}_3'\bs{\delta}_3''=\bs{\delta}_3}}\cM_{\Sigma_{13}}(x^+;\bs{\delta}_1',r,\bs{\delta}_3'')\times\cM(r;\bs{\delta}_1'',x_1,\bs{\delta}_2,x_2,\bs{\delta}_3')\\
  & \bigcup\cM(x^+;\bs{\delta}_1',\xi_{2,1},\bs{\delta}_2''',x_2,\bs{\delta}_3)\times\widetilde{\cM}(\xi_{2,1};\bs{\delta}_1'',\xi_{1,2},\bs{\delta}_2'')\times\cM_{\Sigma_{12}}(\xi_{1,2};\bs{\delta}_1''',x_1,\bs{\delta}_2')\\
  & \bigcup\cM(x^+;\bs{\delta}_1,x_1,\bs{\delta}_2',\xi_{3,2},\bs{\delta}_3''')\times\widetilde{\cM}(\xi_{3,2};\bs{\delta}_2'',\xi_{2,3},\bs{\delta}_3'')\times\cM_{\Sigma_{23}}(\xi_{2,3};\bs{\delta}_2''',x_2,\bs{\delta}_3')\\
  & \bigcup\cM_{\Sigma_{13}}(x^+;\bs{\delta}_1',\xi_{3,1},\bs{\delta}_3''')\times\widetilde{\cM}(\xi_{3,1};\bs{\delta}_1'',\xi_{1,3},\bs{\delta}_3'')\times\cM(\xi_{1,3};\bs{\delta}_1''',x_1,\bs{\delta}_2,x_2,\bs{\delta}'_3)\\
  & \bigcup\cM_{\Sigma_{13}}(x^+;\bs{\delta}_1',\xi_{3,1},\bs{\delta}_3''')\times\widetilde{\cM}(\xi_{3,1};\bs{\delta}_1'',\xi_{1,2},\bs{\delta}_2'',\xi_{2,3},\bs{\delta}_3'')\times\cM_{\Sigma_{12}}(\xi_{1,2};\bs{\delta}_1''',x_1,\bs{\delta}_2')\\
  &\hspace{10cm}\times\cM_{\Sigma_{23}}(\xi_{2,3};\bs{\delta}_2''',x_2,\bs{\delta}_3')
\end{alignat*}
where the $\bs{\delta}_i',\bs{\delta}_i'',\bs{\delta}_i'''$ are words of Reeb chords of $\La_i^-$ such that $\bs{\delta}_i'\bs{\delta}_i''=\bs{\delta}_i$ for the three first unions, and $\bs{\delta}_i'\bs{\delta}_i''\bs{\delta}_i'''=\bs{\delta}_i$ in the four last unions where we sum also respectively for:
 \begin{alignat*}{3}
  & \bullet \quad\xi_{1,2}\in \Rc(\La_1^-,\La_2^-),\,\xi_{2,1}\in\Rc(\La_2^-,\La_1^-),\\
  & \bullet\quad\xi_{3,1}\in \Rc(\La_3^-,\La_1^-),\,\xi_{1,3}\in \Rc(\La_1^-,\La_3^-),\\
  & \bullet\quad\xi_{3,2}\in\Rc(\La_3^-,\La_2^-),\,\xi_{2,3}\in\Rc(\La_2^-,\La_3^-),\\
  & \bullet\quad\xi_{3,1}\in\Rc(\La_3^-,\La_1^-),\,\xi_{2,3}\in\Rc(\La_2^-,\La_3^-),\,\xi_{1,2}\in\Rc(\La_1^-,\La_2^-).
 \end{alignat*}
See Figure \ref{degexx1} for a schematic picture of the above pseudo-holomorphic buildings (except the $\partial$-breakings because we do not draw the pure Reeb chords).
From this, we can deduce Relation \eqref{rel0xx}. Indeed, there is a one-to-one correspondence between buildings involved in the definition of each term in Relation \eqref{rel0xx} (from left to right) and buildings in $\partial\overline{\cM^1}_{\Sigma_{123}}(x^+;\bs{\delta}_1,x_1,\bs{\delta}_2,x_2,\bs{\delta}_3)$ (from left to right on Figure \ref{degexx1}), except that the last term of Relation \eqref{rel0xx} is defined by a count of the two last types of buildings at the right of the figure. Moreover, $\overline{\cM^1}_{\Sigma_{123}}(x^+;\bs{\delta}_1,x_1,\bs{\delta}_2,x_2,\bs{\delta}_3)$ is a compact $1$-dimensional manifold so its boundary consists of an even number of points, hence the count of such points vanishes over $\Z_2$ and we get:
 \begin{alignat*}{1}
  m_{00}^0(x_2,d_{00}(x_1))&+m_{00}^0(d_{00}(x_2),x_1)+d_{00}\circ m_{00}^0(x_2,x_1)\\
  &+m_{0-}^0(x_2,d_{-0}(x_1))+m_{-0}^0(d_{-0}(x_2),x_1)+d_{0-}\circ m_{00}^-(x_2,x_1)=0\nonumber
 \end{alignat*}

\subsubsection{Relation \eqref{rel-xx}}\label{Rel-xx}
\begin{figure}[ht]
      \begin{center}\includegraphics[width=7cm]{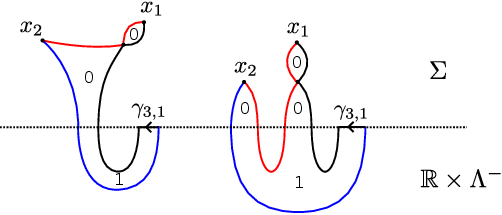}\end{center}
      \caption{Unfinished buildings contributing to $\langle m_{00}^-(x_2,d_{00}(x_1)),\gamma_{3,1}\rangle$.}
      \label{recolle5}
\end{figure}

\begin{figure}[ht]
      \begin{center}\includegraphics[width=8cm]{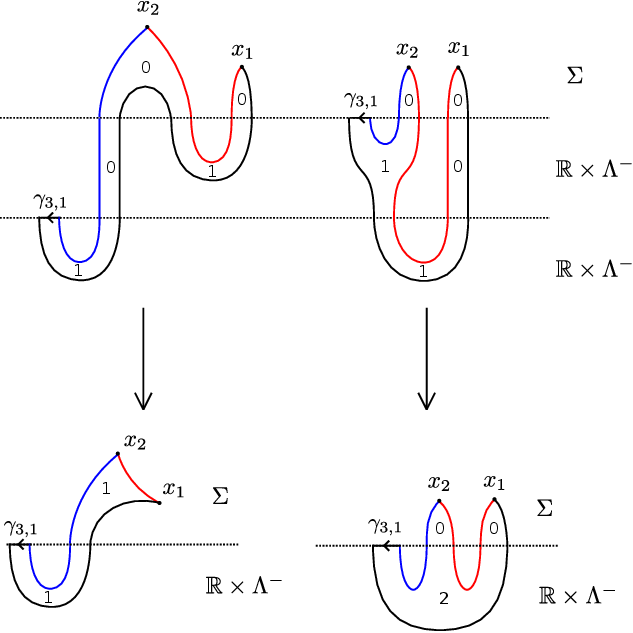}\end{center}
      \caption{Unfinished pseudo-holomorphic buildings contributing to the coefficient $\langle m_{-0}^-(x_2,d_{-0}(x_1)),\gamma_{3,1}\rangle$ and gluing of some levels.}
      \label{recolle6}
\end{figure}

The first term of Relation \eqref{rel-xx} is $m_{00}^-(x_2,d_{00}(x_1))$. For each Reeb chord $\gamma_{3,1}\in\Rc(\La_3^-,\La_1^-)$, the coefficient $\langle m_{00}^-(d_{00}(x_2),x_1),\gamma_{3,1}\rangle$ is defined by a count of unfinished buildings of two types, as we saw in Section \ref{defprod} for the definition of $m^-$. These unfinished buildings are of height $1|1|0$ and the components in the middle level form a pseudo-holomorphic building so its components can be glued (see Figure \ref{recolle5}). Indeed, the curves in the middle level glue together at an intersection point to produce unfinished pseudo-holomorphic buildings of height $1|1|0$  of two types. These live in the following products of moduli spaces:
\begin{alignat*}{1}
   & \widetilde{\cM^1}_{\R\times\La^-_{13}}(\gamma_{3,1};\bs{\xi}_1,\gamma_{1,3},\bs{\xi}_3)\times\cM^1_{\Sigma_{123}}(\gamma_{1,3};\bs{\delta}_1,x_1,\bs{\delta}_2,x_2,\bs{\delta}_3)\\
   & \widetilde{\cM^1}_{\R\times\La^-_{123}}(\gamma_{3,1};\bs{\zeta}_1,\gamma_{1,2},\bs{\zeta}_2,\gamma_{2,3},\bs{\zeta}_3)\times\cM^1_{\Sigma_{12}}(\gamma_{1,2};\bs{\beta}_1,x_1,\bs{\beta}_2)\times\cM^0_{\Sigma_{23}}(\gamma_{2,3};\bs{\delta}_2,x_2,\bs{\delta}_3)\\
\end{alignat*}
Similarly to the previous relation, we will study degeneration of curves in these products of moduli spaces. However, these are not the only one we have to consider. Indeed, the second term of Relation \eqref{rel-xx} is $m_{00}^-(d_{00}(x_2),x_1)$, and analogously to the first term, unfinished pseudo-holomorphic buildings contributing to $\langle m_{00}^-(x_2,d_{00}(x_1)),\gamma_{3,1}\rangle$ for a Reeb chord $\gamma_{3,1}\in\Rc(\La_3^-,\La_1^-)$ have a pseudo-holomorphic building middle level whose components can be glued. After gluing, we get unfinished buildings in the following products:
\begin{alignat*}{1}
   & \widetilde{\cM^1}_{\R\times\La^-_{13}}(\gamma_{3,1};\bs{\xi}_1,\gamma_{1,3},\bs{\xi}_3)\times\cM^1_{\Sigma_{123}}(\gamma_{1,3};\bs{\delta}_1,x_1,\bs{\delta}_2,x_2,\bs{\delta}_3)\\
	& \widetilde{\cM^1}_{\R\times\La^-_{123}}(\gamma_{3,1};\bs{\zeta}_1,\gamma_{1,2},\bs{\zeta}_2,\gamma_{2,3},\bs{\zeta}_3)\times\cM^0_{\Sigma_{12}}(\gamma_{1,2};\bs{\beta}_1,x_1,\bs{\beta}_2)\times\cM^1_{\Sigma_{23}}(\gamma_{2,3};\bs{\delta}_2,x_2,\bs{\delta}_3)\\
\end{alignat*}
Observe that the first product type is the same as one we already obtained above, but the second is different.
Let us consider now the third term of Relation \eqref{rel-xx}, which is by definition a sum
\begin{alignat*}{1}
 m_{-0}^-(x_2,d_{-0}(x_1))=b\circ f^{(2)}(x_2,d_{-0}(x_1))+b^{(2)}(f^{(1)}(x_2),d_{-0}(x_1))
\end{alignat*}
The first term of this sum counts unfinished buildings of height $2|1|0$ such that the components in the middle level and on floor $-1$ are curves contributing to $f^{(2)}(x_2,d_{-0}(x_1))$ and form a pseudo-holomorphic building of height $1|1|0$. On floor $-2$ there is one curve contributing to the map $b$. The second term, $b^{(2)}(f^{(1)}(x_2),d_{-0}(x_1))$, also counts unfinished holomorphic buildings of height $2|1|0$ but this time the components in the bottom level on floors $-2$ and $-1$ give a holomorphic building of height $2|0|0$. Gluing the components of these buildings, we get unfinished buildings in the following products (see Figure \ref{recolle6}):
\begin{alignat*}{1}
   & \widetilde{\cM^1}_{\R\times\La^-_{13}}(\gamma_{3,1};\bs{\xi}_1,\gamma_{1,3},\bs{\xi}_3)\times\cM^1_{\Sigma_{123}}(\gamma_{1,3};\bs{\delta}_1,x_1,\bs{\delta}_2,x_2,\bs{\delta}_3)\\
  & \widetilde{\cM^2}_{\R\times\La^-_{123}}(\gamma_{3,1};\bs{\zeta}_1,\gamma_{1,2},\bs{\zeta}_2,\gamma_{2,3},\bs{\zeta}_3)\times\cM^0_{\Sigma_{12}}(\gamma_{1,2};\bs{\beta}_1,x_1,\bs{\beta}_2)\times\cM^0_{\Sigma_{23}}(\gamma_{2,3};\bs{\delta}_2,x_2,\bs{\delta}_3)\\
\end{alignat*}
Again, we already got the first type of product but the second product is a new one we will have to study.
Then, the fourth term of Relation \eqref{rel-xx}, $m_{0-}^-(d_{-0}(x_2),x_1))$, is symmetric  to the third and so counts unfinished buildings such that some levels can be glued to give unfinished buildings in the same products of moduli spaces as above (for study of the third term). The fifth term is $d_{-0}\circ m_{00}^0(x_2,x_1)$. The middle level of unfinished buildings contributing to this term is a pseudo-holomorphic building with two curves which glue together at an intersection point in $\Sigma_1\cap\Sigma_3$. After gluing, we get unfinished holomorphic buildings in the product:
\begin{alignat*}{1}
 & \widetilde{\cM^1}_{\R\times\La^-_{13}}(\gamma_{3,1};\bs{\xi}_1,\gamma_{1,3},\bs{\xi}_3)\times\cM^1_{\Sigma_{123}}(\gamma_{1,3};\bs{\delta}_1,x_1,\bs{\delta}_2,x_2,\bs{\delta}_3)\\
\end{alignat*}
Finally, the last term of the relation, $d_{--}\circ m_{00}^-(x_2,x_1)$, counts unfinished buildings of height $2|1|0$ of two types, such that the components in the bottom level on floors $-2$ and $-1$ form a holomorphic buildings of height $2|0|0$. Gluing these components, we get unfinished buildings in the products:
 \begin{alignat*}{1}
   & \widetilde{\cM^2}_{\R\times\La^-_{13}}(\gamma_{3,1};\bs{\xi}_1,\gamma_{1,3},\bs{\xi}_3)\times\cM^0_{\Sigma_{123}}(\gamma_{1,3};\bs{\delta}_1,x_1,\bs{\delta}_2,x_2,\bs{\delta}_3)\\
   & \widetilde{\cM^2}_{\R\times\La^-_{123}}(\gamma_{3,1};\bs{\zeta}_1,\gamma_{1,2},\bs{\zeta}_2,\gamma_{2,3},\bs{\zeta}_3)\times\cM^0_{\Sigma_{12}}(\gamma_{1,2};\bs{\beta}_1,x_1,\bs{\beta}_2)\times\cM^0_{\Sigma_{23}}(\gamma_{2,3};\bs{\delta}_2,x_2,\bs{\delta}_3)\\
 \end{alignat*}
Now, in order to obtain Relation \eqref{rel-xx}, we need to study the boundary of the compactification of each product of moduli spaces appearing above, where all the broken curves that are involved in the definition of each term of the relation live. So, to sum up, we must study the boundary of the compactification of the following products:   
\begin{alignat}{1}
   &\widetilde{\cM^1}_{\R\times\La^-_{13}}(\gamma_{3,1};\bs{\xi}_1,\gamma_{1,3},\bs{\xi}_3)\times\cM^1_{\Sigma_{123}}(\gamma_{1,3};\bs{\delta}_1,x_1,\bs{\delta}_2,x_2,\bs{\delta}_3)\label{prodcomp1} \\
   &\widetilde{\cM^2}_{\R\times\La^-_{13}}(\gamma_{3,1};\bs{\xi}_1,\gamma_{1,3},\bs{\xi}_3)\times\cM^0_{\Sigma_{123}}(\gamma_{1,3};\bs{\delta}_1,x_1,\bs{\delta}_2,x_2,\bs{\delta}_3) \label{prodcomp2}\\
   &\widetilde{\cM^1}_{\R\times\La^-_{123}}(\gamma_{3,1};\bs{\zeta}_1,\gamma_{1,2},\bs{\zeta}_2,\gamma_{2,3},\bs{\zeta}_3)\times\cM^0_{\Sigma_{12}}(\gamma_{1,2};\bs{\beta}_1,x_1,\bs{\beta}_2)\times\cM^1_{\Sigma_{23}}(\gamma_{2,3};\bs{\delta}_2,x_2,\bs{\delta}_3) \label{prodcomp3} \\ 
   &\widetilde{\cM^1}_{\R\times\La^-_{123}}(\gamma_{3,1};\bs{\zeta}_1,\gamma_{1,2},\bs{\zeta}_2,\gamma_{2,3},\bs{\zeta}_3)\times\cM^1_{\Sigma_{12}}(\gamma_{1,2};\bs{\beta}_1,x_1,\bs{\beta}_2)\times\cM^0_{\Sigma_{23}}(\gamma_{2,3};\bs{\delta}_2,x_2,\bs{\delta}_3) \label{prodcomp4}\\ 
   &\widetilde{\cM^2}_{\R\times\La^-_{123}}(\gamma_{3,1};\bs{\zeta}_1,\gamma_{1,2},\bs{\zeta}_2,\gamma_{2,3},\bs{\zeta}_3)\times\cM^0_{\Sigma_{12}}(\gamma_{1,2};\bs{\beta}_1,x_1,\bs{\beta}_2)\times\cM^0_{\Sigma_{23}}(\gamma_{2,3};\bs{\delta}_2,x_2,\bs{\delta}_3)\label{prodcomp5}
\end{alignat}
In these products, moduli spaces of index-$0$ curves with boundary on non-cylindrical Lagrangians are compact $0$-dimensional manifolds, as well as the quotient of moduli spaces of index-$1$ curves with boundary on the negative cylindrical ends of the Lagrangian cobordisms. On the other hand, moduli spaces of index-$1$ curves with boundary on non-cylindrical Lagrangians are non compact $1$-dimensional manifolds, as well as the quotient of moduli spaces of index-$2$ curves with boundary on cylindrical Lagrangians. By compactness results, these $1$-dimensional moduli spaces can be compactified and the boundary of the compactification consists of pseudo-holomorphic buildings with rigid components. Thus, we need to describe the followings spaces:
\begin{enumerate}
 \item[1.] $\partial\overline{\cM^1}(\gamma_{1,3};\bs{\delta}_1,x_1,\bs{\delta}_2,x_2,\bs{\delta}_3),$
 \item[2.] $\partial\overline{\cM^2}(\gamma_{3,1};\bs{\xi}_1,\gamma_{1,3},\bs{\xi}_3),$
 \item[3.] $\partial\overline{\cM^1}(\gamma_{2,3};\bs{\delta}_2,x_2,\bs{\delta}_3),$
 \item[4.] $\partial\overline{\cM^1}(\gamma_{1,2};\bs{\beta}_1,x_1,\bs{\beta}_2),$
 \item[5.] $\partial\overline{\cM^2}(\gamma_{3,1};\bs{\zeta}_1,\gamma_{1,2},\bs{\zeta}_2,\gamma_{2,3},\bs{\zeta_3})$
\end{enumerate}
where we write $\overline{\cM^2}$ instead of $\overline{\widetilde{\cM^2}}=\overline{(\cM^2/\R)}$ to simplify notation for the compactification of the quotient of a moduli space of index-$2$ curves with boundary on cylindrical Lagrangians.
Once we understand the boundaries of the compactified moduli spaces above, we understand all the broken curves appearing as degeneration of unfinished buildings in the products \eqref{prodcomp1}, \eqref{prodcomp2}, \eqref{prodcomp3}, \eqref{prodcomp4} and \eqref{prodcomp5}.

\vspace{2mm}

\begin{figure}[ht]
      \begin{center}\includegraphics[width=14cm]{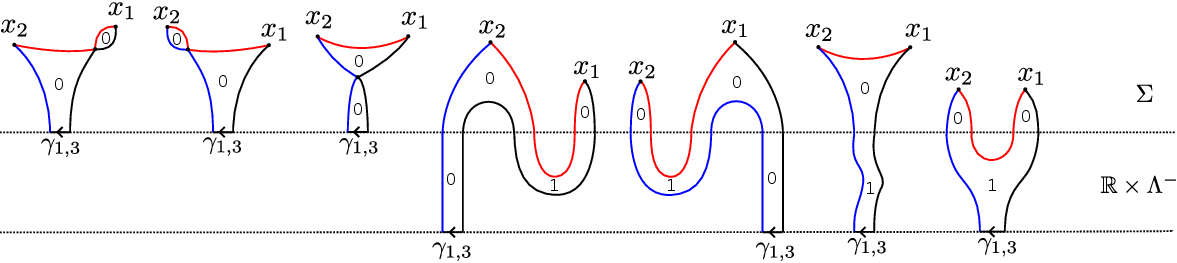}\end{center}
      \caption{Pseudo-holomorphic buildings in the boundary of $\overline{\cM^1}(\gamma_{1,3};\bs{\delta}_1,x_1,\bs{\delta}_2,x_2,\bs{\delta}_3)$.}
      \label{degene_f_1}
\end{figure}

\noindent1. $\partial\overline{\cM^1}(\gamma_{1,3};\bs{\delta}_1,x_1,\bs{\delta}_2,x_2,\bs{\delta}_3)$: the different pseudo-holomorphic buildings in this space are listed below, where the unions are, depending on cases, over intersection points $p\in\Sigma_1\cap\Sigma_2$, $q\in\Sigma_2\cap\Sigma_3$, $r\in\Sigma_1\cap\Sigma_3$, Reeb chords $\xi_{i,j}\in\Rc(\La_i^-,\La_j^-)$ for $1\leq i\neq j\leq3$, and words of pure chords $\bs{\delta}_i', \bs{\delta}_i'', \bs{\delta}_i'''$ of $\La_i^-$ for $i=1,2,3$ satisfying $\bs{\delta}_i'\bs{\delta}_i''=\bs{\delta}_i$ or $\bs{\delta}_i'\bs{\delta}_i''\bs{\delta}_i'''=\bs{\delta}_i$.
\begin{alignat*}{1}
   & \partial\overline{\cM^1}(\gamma_{1,3};\bs{\delta}_1,x_1,\bs{\delta}_2,x_2,\bs{\delta}_3)=\overline{\cM}^\partial(\gamma_{1,3};\bs{\delta}_1,x_1,\bs{\delta}_2,x_2,\bs{\delta}_3)\\
   &\bigcup_{p,\bs{\delta}_i',\bs{\delta}_i''}\cM(\gamma_{1,3};\bs{\delta}_1',p,\bs{\delta}_2'',x_2,\bs{\delta}_3)\times\cM(p;\bs{\delta}_1'',x_1,\bs{\delta}_2')\\
   &\bigcup_{q,\bs{\delta}_i'\bs{\delta}_i''}\cM(\gamma_{1,3};\bs{\delta}_1,x_1,\bs{\delta}_2',q,\bs{\delta}_3'')\times\cM(q;\bs{\delta}_2'',x_2,\bs{\delta}_3')\\
   &\bigcup_{r,\bs{\delta}_i',\bs{\delta}_i''}\cM(\gamma_{1,3};\bs{\delta}_1',r,\bs{\delta}_3'')\times\cM(r;\bs{\delta}_1'',x_1,\bs{\delta}_2,x_2,\bs{\delta}_3')\\
   & \bigcup_{\substack{\xi_{2,1},\xi_{1,2}\\ \bs{\delta}_i',\bs{\delta}_i'',\bs{\delta}_i'''}}\cM(\gamma_{1,3};\bs{\delta}_1',\xi_{2,1},\bs{\delta}_2''',x_2,\bs{\delta}_3)\times\widetilde{\cM}(\xi_{2,1};\bs{\delta}_1'',\xi_{1,2},\bs{\delta}_2'')\times\cM(\xi_{1,2};\bs{\delta}_1''',x_1,\bs{\delta}_2')\\
   & \bigcup_{\substack{\xi_{2,3},\xi_{3,2}\\ \bs{\delta}_i',\bs{\delta}_i'',\bs{\delta}_i'''}}\cM(\gamma_{1,3};\bs{\delta}_1,x_1,\bs{\delta}_2',\xi_{3,2},\bs{\delta}_3''')\times\widetilde{\cM}(\xi_{3,2};\bs{\delta}_2'',\xi_{2,3},\bs{\delta}_3'')\times\cM(\xi_{2,3};\bs{\delta}_2''',x_2,\bs{\delta}_3')\\
   & \bigcup_{\substack{\xi_{1,3}\\ \bs{\delta}_i',\bs{\delta}_i''}}\widetilde{\cM}(\gamma_{1,3};\bs{\delta}_1',\xi_{1,3},\bs{\delta}_3'')\times\cM(\xi_{1,3};\bs{\delta}_1'',x_1,\bs{\delta}_2,x_2,\bs{\delta}_3')\\
   & \bigcup_{\substack{\xi_{2,3},\xi_{1,2}\\ \bs{\delta}_i',\bs{\delta}_i'',\bs{\delta}_i'''}}\widetilde{\cM}(\gamma_{1,3};\bs{\delta}_1',\xi_{1,2},\bs{\delta}_2'',\xi_{2,3},\bs{\delta}_3'')\times\cM(\xi_{2,3};\bs{\delta}_2''',x_2,\bs{\delta}_3')\times\cM(\xi_{1,2};\bs{\delta}_1'',x_1,\bs{\delta}_2')
\end{alignat*}
See Figure \ref{degene_f_1} for a schematic picture of the different types of pseudo-holomorphic buildings in $\partial\overline{\cM^1}(x_2;\bs{\delta}_3,\gamma_{1,3},\bs{\delta}_1,x_1,\bs{\delta}_2)$.

\begin{figure}[ht]
      \begin{center}\includegraphics[width=10cm]{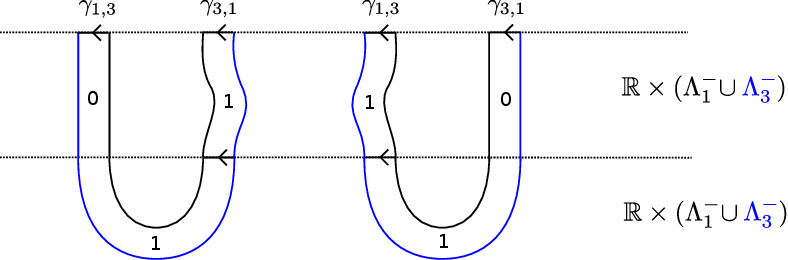}\end{center}
      \caption{Pseudo-holomorphic buildings in the boundary of $\overline{\cM^2}(\gamma_{3,1};\bs{\delta}_1,\gamma_{1,3},\bs{\delta}_3)$.}
      \label{brisure_banane}
\end{figure}

\noindent 2. $\partial\overline{\cM^2}(\gamma_{3,1};\bs{\xi}_1,\gamma_{1,3},\bs{\xi}_3)$: pseudo-holomorphic buildings appearing as degeneration of index-$2$ bananas are of two types. We have:
\begin{alignat*}{2}
   & \partial\overline{\cM^2}(\gamma_{3,1};\bs{\xi}_1,\gamma_{1,3},\bs{\xi}_3)=\overline{\cM}^\partial(\gamma_{3,1};\bs{\xi}_1,\gamma_{1,3},\bs{\xi}_3)\\
   & \bigcup_{\xi_{3,1},\bs{\xi}_i',\bs{\xi}_i''}\widetilde{\cM}(\gamma_{3,1};\bs{\xi}_1',\xi_{3,1},\bs{\xi}_3'')\times\widetilde{\cM}(\xi_{3,1};\bs{\xi}_1'',\gamma_{1,3},\bs{\xi}_3')\\
   & \bigcup_{\xi_{1,3},\bs{\xi}_i',\bs{\xi}_i''}\widetilde{\cM}(\gamma_{3,1};\bs{\xi}_1',\xi_{1,3},\bs{\xi}_3'')\times\widetilde{\cM}(\xi_{1,3};\bs{\xi}_1'',\gamma_{1,3},\bs{\xi}_3')\\
\end{alignat*}
with again $\xi_{i,j}\in\Rc(\La_i^-,\La_j^-)$, and $\bs{\xi}_i',\bs{\xi}_i''$ words of Reeb chords of $\La_i^-$, with $\bs{\xi}_i'\bs{\xi}_i''=\bs{\xi}_i$ (see Figure \ref{brisure_banane}).

\begin{figure}[ht]
      \begin{center}\includegraphics[width=145mm]{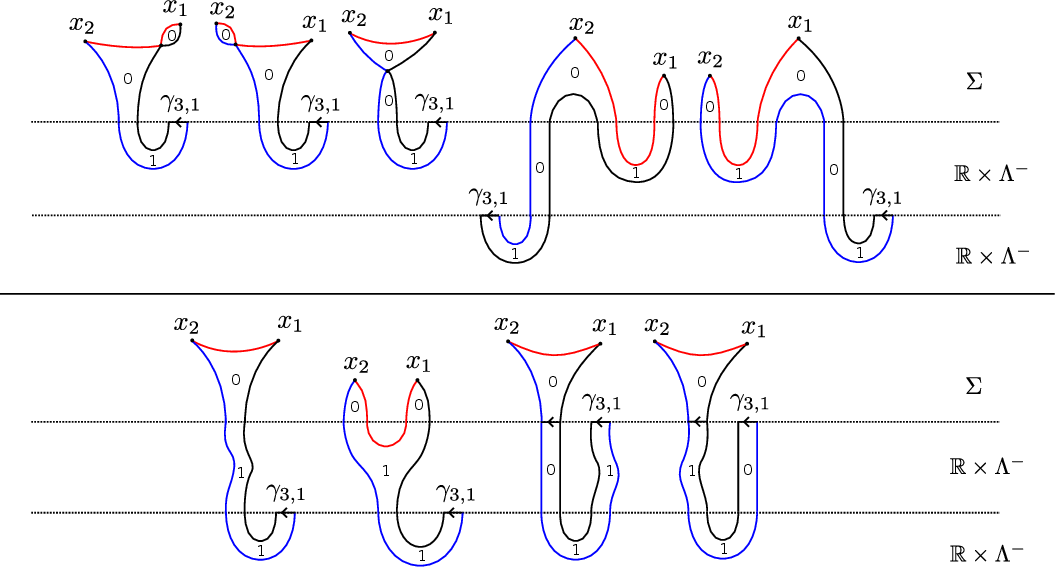}\end{center}
      \caption{Unfinished pseudo-holomorphic buildings in \eqref{prodcomp6} and \eqref{prodcomp7}.}
      \label{degene_f1_banane}
\end{figure}

With $1.$ and $2.$ above, we can describe all the types of broken curves in the boundary of the compactification of the products \eqref{prodcomp1} and \eqref{prodcomp2}. Instead of writing again huge unions of moduli spaces, in Figure \ref{degene_f1_banane} we drew schematic pictures of the corresponding unfinished holomorphic buildings. The first seven (from left to right and top to bottom) are in:
\begin{alignat}{1}
  \widetilde{\cM^1}(\gamma_{3,1};\bs{\xi}_1,\gamma_{1,3},\bs{\xi}_3)\times\partial\overline{\cM^1}(\gamma_{1,3};\bs{\delta}_1,x_1,\bs{\delta}_2,x_2,\bs{\delta}_3)\label{prodcomp6}
\end{alignat}
and the last two are in:
\begin{alignat}{1}
  \partial\overline{\cM^2}(\gamma_{3,1};\bs{\xi}_1,\gamma_{1,3},\bs{\xi}_3)\times\cM^0(\gamma_{1,3};\bs{\delta}_1,x_1,\bs{\delta}_2,x_2,\bs{\delta}_3)\label{prodcomp7}
\end{alignat}
\begin{rem}\label{rem_bord}
 In the bottom of Figure \ref{degene_f1_banane}, two of the unfinished holomorphic buildings are equivalent (Definition \ref{defequi}): the leftmost compensate with the rightmost. The leftmost is in \eqref{prodcomp6}: the component in the middle level together with the component in the bottom level in floor $-1$ form a pseudo-holomorphic building which lives in $\partial\overline{\cM^1}(\gamma_{1,3};\bs{\delta}_1,x_1,\bs{\delta}_2,x_2,\bs{\delta}_3)$. The rightmost unfinished building in the bottom of the figure is in \eqref{prodcomp7}: the components of the bottom level form a pseudo-holomorphic building living in $\partial\overline{\cM^2}(\gamma_{3,1};\bs{\xi}_1,\gamma_{1,3},\bs{\xi}_3)$. These two unfinished buildings correspond thus to different geometric configurations because they live in the boundary of the compactification of two different products of moduli spaces (one is in \eqref{prodcomp6}, the other in \eqref{prodcomp7}). However, these buildings differ only by a trivial strip $\R\times\gamma_{3,1}$ so they contribute algebraically to the same map which is in this case $b\circ\delta_{--}\circ f^{(2)}(x_2,x_1)$, where $\delta_{--}$ is the dual of $d_{--}$.
\end{rem}
In order deduce the algebraic relation that these boundary elements give, we introduce a new map:
\begin{alignat*}{1}
 \Delta^{(2)}\colon\mathfrak{C}_{n-1-*}(\La_3^-,\La_2^-)\times\mathfrak{C}_{n-1-*}(\La_2^-,\La_1^-)\to C_{n-1-*}(\La_3^-,\La_1^-)
\end{alignat*}
defined on pairs of generators by:
     \begin{alignat*}{1}
      &\Delta^{(2)}(\gamma_{2,3},\gamma_{1,2})=\sum_{\gamma_{1,3}}\#\cM^0(\gamma_{1,3};\bs{\delta}_1,\gamma_{1,2},\bs{\delta}_2,\gamma_{2,3},\bs{\delta}_3)\ep_1^-(\bs{\delta}_1)\ep_2^-(\bs{\delta}_2)\ep_3^-(\bs{\delta}_3)\cdot\gamma_{1,3}\\
      &\Delta^{(2)}(\gamma_{2,3},\gamma_{2,1})=\sum_{\gamma_{1,3}}\#\cM^0(\gamma_{1,3};\bs{\delta}_1,\gamma_{2,1},\bs{\delta}_2,\gamma_{2,3},\bs{\delta}_3)\ep_1^-(\bs{\delta}_1)\ep_2^-(\bs{\delta}_2)\ep_3^-(\bs{\delta}_3)\cdot\gamma_{1,3}\\
      &\Delta^{(2)}(\gamma_{3,2},\gamma_{1,2})=\sum_{\gamma_{1,3}}\#\cM^0(\gamma_{1,3};\bs{\delta}_1,\gamma_{1,2},\bs{\delta}_2,\gamma_{3,2},\bs{\delta}_3)\ep_1^-(\bs{\delta}_1)\ep_2^-(\bs{\delta}_2)\ep_3^-(\bs{\delta}_3)\cdot\gamma_{1,3}\\
      &\Delta^{(2)}(\gamma_{3,2},\gamma_{2,1})=0
     \end{alignat*}
Now, as for Relation \eqref{rel0xx}, the mod-$2$ count of unfinished pseudo-holomorphic buildings in the products \eqref{prodcomp6} and \eqref{prodcomp7} equals $0$. On the other hand, these broken curves contribute to some composition of maps we have defined earlier. This implies that the following relation is satisfied:
    \begin{alignat}{1}
      & b\circ f^{(2)}(x_2,d_{00}(x_1))+b\circ f^{(2)}(d_{00}(x_2),x_1)+d_{-0}\circ m_{00}^0(x_2,x_1)\nonumber\\ 
      & +b\circ f^{(2)}(x_2,d_{-0}(x_1))+b\circ f^{(2)}(d_{-0}(x_2),x_1)+b\circ \Delta^{(2)}(f^{(1)}(x_2),f^{(1)}(x_1))\label{rel-xx1}\\
      &\hspace{71mm}+d_{--}\circ b\circ f^{(2)}(x_2,x_1)=0\nonumber
    \end{alignat}
where we did not write the term $b\circ\delta_{--}\circ f^{(2)}(x_2,x_1)$ as it would appear twice so this vanishes over $\Z_2$ (see Remark \ref{rem_bord}).
\vspace{2mm}

\begin{figure}[ht]
      \begin{center}\includegraphics[width=3cm]{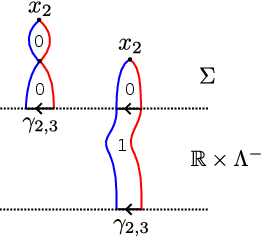}\end{center}
      \caption{Pseudo-holomorphic buildings in $\partial\overline{\cM^1}(\gamma_{2,3};\bs{\delta}_2,x_2,\bs{\delta}_3)$.}
      \label{degene_delta_-0}
\end{figure}

\noindent 3. $\partial\overline{\cM^1}(\gamma_{2,3};\bs{\delta}_2,x_2,\bs{\delta}_3)$: pseudo-holomorphic buildings in this space are of the following type (see Figure \ref{degene_delta_-0}).
\begin{alignat*}{2}
 &\partial\overline{\cM^1}(\gamma_{2,3};\bs{\delta}_2,x_2,\bs{\delta}_3)=\overline{\cM^1}^\partial(\gamma_{2,3};\bs{\delta}_2,x_2,\bs{\delta}_3)\\
 &\bigcup_{\substack{q\in\Sigma_2\cap\Sigma_3\\ \bs{\delta_i}',\bs{\delta}_i''}}\cM(\gamma_{2,3};\bs{\delta}_2',q,\bs{\delta}_3'',)\times\cM(q;\bs{\delta}_2'',x_2,\bs{\delta}_3')\\
 &\bigcup_{\substack{\xi_{2,3},\bs{\delta_i}',\bs{\delta}_i''}}\widetilde{\cM}(\gamma_{2,3};\bs{\delta}_2',\xi_{2,3},\bs{\delta}_3'')\times\cM(\xi_{2,3};\bs{\delta}_2'',x_2,\bs{\delta}_3')\\
\end{alignat*}

\noindent 4. $\partial\overline{\cM^1}(\gamma_{1,2};\bs{\beta}_1,x_1,\bs{\beta}_2)$: same types of degenerations as above (case $3.$).

\vspace{2mm}

\begin{figure}[ht]
      \begin{center}\includegraphics[width=14cm]{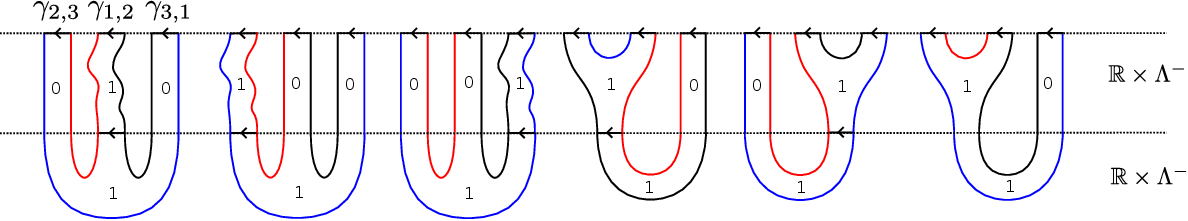}\end{center}
      \caption{Pseudo-holomorphic buildings in the boundary of the compactification of $\widetilde{\cM^2}(\gamma_{3,1};\bs{\delta}_1,\gamma_{1,2},\bs{\delta}_2,\gamma_{2,3},\bs{\delta}_3)$.}
      \label{brisure_banane1}
\end{figure}

\noindent 5. $\partial\overline{\cM^2}(\gamma_{3,1};\bs{\zeta}_1,\gamma_{1,2},\bs{\zeta}_2,\gamma_{2,3},\bs{\zeta}_3)$: we describe here degenerations of index-$2$ bananas with three positive Reeb chords asymptotics (see Figure \ref{brisure_banane1} for a schematic picture of the corresponding broken curves).
\begin{alignat*}{2}
   & \partial\overline{\cM^2}(\gamma_{3,1};\bs{\zeta}_1,\gamma_{1,2},\bs{\zeta}_2,\gamma_{2,3},\bs{\zeta}_3)=\overline{\cM}^\partial(\gamma_{3,1};\bs{\zeta}_1,\gamma_{1,2},\bs{\zeta}_2,\gamma_{2,3},\bs{\zeta}_3)\\
   & \bigcup_{\xi_{1,2},\bs{\zeta}_i'\bs{\zeta}_i''}\widetilde{\cM}(\gamma_{3,1};\bs{\zeta}_1',\xi_{1,2},\bs{\zeta}_2'',\gamma_{2,3},\bs{\zeta}_3)\times\widetilde{\cM}(\xi_{1,2};\bs{\zeta}_1'',\gamma_{1,2},\bs{\zeta}_2')\\
   & \bigcup_{\xi_{2,3},\bs{\zeta}_i'\bs{\zeta}_i''}\widetilde{\cM}(\gamma_{3,1};\bs{\zeta}_1,\gamma_{1,2},\bs{\zeta}_2',\xi_{2,3},\bs{\zeta}_3'')\times\widetilde{\cM}(\xi_{2,3};\bs{\zeta}_2'',\gamma_{2,3};\bs{\zeta}_3')\\
   & \bigcup_{\xi_{3,1},\bs{\zeta}_i'\bs{\zeta}_i''}\widetilde{\cM}(\gamma_{3,1};\bs{\zeta}_1',\xi_{3,1},\bs{\zeta}_3'')\times\widetilde{\cM}(\xi_{3,1};\bs{\zeta}_1'',\gamma_{1,2},\bs{\zeta}_2,\gamma_{2,3},\bs{\zeta}_3')\\
   & \bigcup_{\xi_{2,1},\bs{\zeta}_i'\bs{\zeta}_i''}\widetilde{\cM}(\gamma_{3,1};\bs{\zeta}_1',\xi_{2,1},\bs{\zeta}_2'',\gamma_{2,3},\bs{\zeta}_3)\times\widetilde{\cM}(\xi_{2,1};\bs{\zeta}_1'',\gamma_{1,2},\bs{\zeta}_2')\\
   & \bigcup_{\xi_{3,2},\bs{\zeta}_i'\bs{\zeta}_i''}\widetilde{\cM}(\gamma_{3,1};\bs{\zeta}_1,\gamma_{2,1},\bs{\zeta}_2',\xi_{3,2},\bs{\zeta}_3'')\times\widetilde{\cM}(\xi_{3,2};\bs{\zeta}_2'',\gamma_{2,3},\bs{\zeta}_3')\\
   & \bigcup_{\xi_{1,3},\bs{\zeta}_i'\bs{\zeta}_i''}\widetilde{\cM}(\gamma_{3,1};\bs{\zeta}_1',\xi_{1,3},\bs{\zeta}_3'')\times\widetilde{\cM}(\xi_{1,3};\bs{\zeta}_1'',\gamma_{1,2},\bs{\zeta}_2,\gamma_{2,3},\bs{\zeta}_3')\\
\end{alignat*}
where $\bs{\zeta}_i',\bs{\zeta}_i''$ are words of pure chords of $\La_i^-$ such that $\bs{\zeta}_i'\bs{\zeta}_i''=\bs{\zeta}_i$.
\begin{figure}[ht]
      \begin{center}\includegraphics[width=13cm]{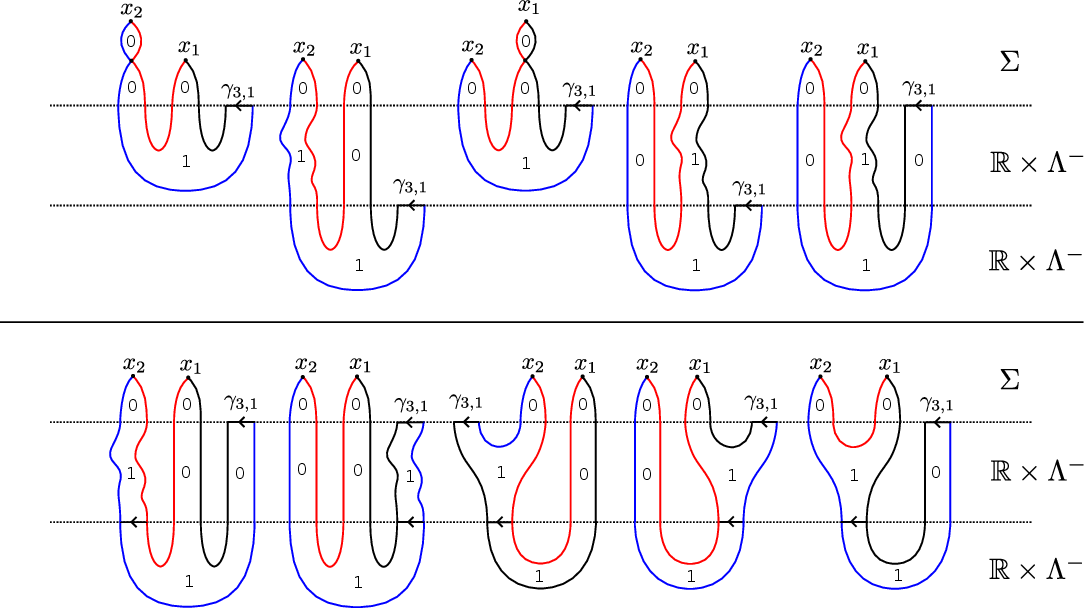}\end{center}
      \caption{Unfinished buildings in the boundary of the compactification of the products \eqref{prodcomp3}, \eqref{prodcomp4} and \eqref{prodcomp5}.}
      \label{degene_delta_b1}
\end{figure}

By $3., 4.$ and $5.$, we can describe all the types of unfinished pseudo-holomorphic buildings in the boundary of the compactification of the products \eqref{prodcomp3}, \eqref{prodcomp4} and \eqref{prodcomp5}, that is to say, unfinished buildings in the spaces:
  \begin{alignat*}{2}
   & \widetilde{\cM^1}(\gamma_{3,1};\bs{\zeta}_1,\gamma_{1,2},\bs{\zeta}_2,\gamma_{2,3},\bs{\zeta}_3)\times\partial\overline{\cM^1}(\gamma_{2,3};\bs{\delta}_2,x_2,\bs{\delta}_3)\times\cM^0(\gamma_{1,2};\bs{\beta}_1,x_1,\bs{\beta}_2)\\
   & \bigcup\widetilde{\cM^1}(\gamma_{3,1};\bs{\zeta}_1,\gamma_{1,2},\bs{\zeta}_2,\gamma_{2,3},\bs{\zeta}_3)\times\cM^0(\gamma_{2,3};\bs{\delta}_2,x_2,\bs{\delta}_3)\times\partial\overline{\cM^1}(\gamma_{1,2};\bs{\beta}_1,x_1,\bs{\beta}_2)\\
   & \bigcup\partial\overline{\cM^2}(\gamma_{3,1};\bs{\zeta}_1,\gamma_{1,2},\bs{\zeta}_2,\gamma_{2,3},\bs{\zeta}_3)\times\cM^0(\gamma_{2,3};\bs{\delta}_2,x_2,\bs{\delta}_3)\times\cM^0(\gamma_{1,2};\bs{\beta}_1,x_1,\bs{\beta}_2)
  \end{alignat*}
The corresponding buildings are schematized on Figure \ref{degene_delta_b1}: the first two (from left to right and top to bottom) are in the boundary of the compactification of \eqref{prodcomp3}, the following two are in the boundary of the compactification of \eqref{prodcomp4}, and finally the six others are in the boundary of the compactification of \eqref{prodcomp5}. As in the previous case (see Remark \ref{rem_bord}), several unfinished holomorphic buildings are equivalent (so their algebraic contributions are the same). Indeed, the second and the sixth one, differing by a trivial strip $\R\times\gamma_{3,1}$ contribute to $b^{(2)}(\delta_{--}\circ f^{(1)}(x_2),f^{(1)}(x_1))$, and the fourth and the fifth, differing also by the same type of trivial strip contribute to $b^{(2)}(f^{(1)}(x_2),\delta_{--}\circ f^{(1)}(x_1))$. We obtain this time the relation:
   \begin{alignat}{2}
    & b^{(2)}(f^{(1)}\circ d_{00}(x_2),f^{(1)}(x_1))+b^{(2)}(f^{(1)}(x_2),f^{(1)}\circ d_{00}(x_1))\nonumber\\
    & +d_{--}\circ b^{(2)}(f^{(1)}(x_2),f^{(1)}(x_1))+b^{(2)}(f^{(1)}(x_2),d_{-0}(x_1))+b^{(2)}(d_{-0}(x_2),f^{(1)}(x_1))\label{rel-xx2}\\
    &\hspace{73mm}+b\circ \Delta^{(2)}(f^{(1)}(x_2),f^{(1)}(x_1))=0 \nonumber
   \end{alignat}
Combining Relations \eqref{rel-xx1} and \eqref{rel-xx2}, we get Relation \eqref{rel-xx}:
  \begin{alignat*}{1}
   & m_{00}^-(d_{00}(x_2),x_1)+m_{00}^-(x_2,d_{00}(x_1))+m_{-0}^-(d_{-0}(x_2),x_1)\\
   &\hspace{2cm}+m_{0-}^-(x_2,d_{-0}(x_1))+d_{-0}\circ m_{00}^0(x_2,x_1)+d_{--}\circ m_{00}^-(x_2,x_1)=0
  \end{alignat*}
where the term $b\circ \Delta^{(2)}(f^{(1)}(x_2),f^{(1)}(x_1))$ disappeared because it is at the same time in \eqref{rel-xx1} and \eqref{rel-xx2} so vanishes over $\Z_2$.

\subsubsection{Relation \eqref{rel0xc}}
This relation is really analogous to Relation \eqref{rel0xx} except that one of the three mixed asymptotics is a Reeb chord. Each term in \eqref{rel0xc} counts pseudo-holomorphic buildings of height $0|1|0$ or $1|1|0$ whose components can be glued on index-$1$ disks in the moduli space $\cM^1_{\Sigma_{123}}(x^+;\bs{\delta}_1,\gamma_1,\bs{\delta}_2,x_2,\bs{\delta}_3)$. To determine Relation \eqref{rel0xc}, we have thus to study the broken curves in the boundary of the compactification of this moduli space. This gives:
\begin{alignat*}{2}
  & \partial\overline{\cM^1}_{\Sigma_{123}}(x^+;\bs{\delta}_1,\gamma_1,\bs{\delta}_2,x_2,\bs{\delta}_3)=\overline{\cM}^\partial(x^+,\bs{\delta}_1,\gamma_1,\bs{\delta}_2,x_2,\bs{\delta}_3)\\
  & \bigcup_{p,\bs{\delta}_i',\bs{\delta}_i''}\cM(x^+;\bs{\delta}_1',p,\bs{\delta}_2'',x_2,\bs{\delta}_3)\times\cM_{\Sigma_{12}}(p;\bs{\delta}_1'',\gamma_1,\bs{\delta}_2')\\
  & \bigcup_{q,\bs{\delta}_i',\bs{\delta}_i''}\cM(x^+;\bs{\delta}_1,\gamma_1,\bs{\delta}_2',q,\bs{\delta}_3'')\times\cM_{\Sigma_{23}}(q;\bs{\delta}_2'',x_2,\bs{\delta}_3')\\
  & \bigcup_{r,\bs{\delta}_i',\bs{\delta}_i''}\cM_{\Sigma_{13}}(x^+;\bs{\delta}_1',r,\bs{\delta}_3'')\times\cM(r;\bs{\delta}_1'',\gamma_1,\bs{\delta}_2,x_2,\bs{\delta}_3')\\
  & \bigcup_{\xi_{2,1},\bs{\delta}_i',\bs{\delta}_i''}\cM(x^+;\bs{\delta}_1',\xi_{2,1},\bs{\delta}_2'',x_2,\bs{\delta}_3)\times\widetilde{\cM}(\xi_{2,1};\bs{\delta}_1'',\gamma_1,\bs{\delta}_2')\\
  & \bigcup_{\substack{\xi_{3,2},\xi_{2,3}\\ \bs{\delta}_i',\bs{\delta}_i'',\bs{\delta}_i'''}}\cM(x^+;\bs{\delta}_1,\gamma_1,\bs{\delta}_2',\xi_{3,2},\bs{\delta}_3''')\times\widetilde{\cM}(\xi_{3,2};\bs{\delta}_2'',\xi_{2,3},\bs{\delta}_3'')\times\cM_{\Sigma_{23}}(\xi_{2,3};\bs{\delta}_2''',x_2,\bs{\delta}_3')\\
  & \bigcup_{\substack{\xi_{3,1},\xi_{1,3}\\ \bs{\delta}_i',\bs{\delta}_i'',\bs{\delta}_i'''}}\cM(x^+;\bs{\delta}_1',\xi_{3,1},\bs{\delta}_3''')\times\widetilde{\cM}(\xi_{3,1};\bs{\delta}_1'',\xi_{1,3},\bs{\delta}_3'')\times\cM(\xi_{1,3};\bs{\delta}_1''',\gamma_1,\bs{\delta}_2,x_2,\bs{\delta}_3')\\
  & \bigcup_{\substack{\xi_{3,1},\xi_{2,3}\\ \bs{\delta}_i',\bs{\delta}_i'',\bs{\delta}_i'''}}\cM(x^+;\bs{\delta}_1',\xi_{3,1},\bs{\delta}_3''')\times\widetilde{\cM}(\xi_{3,1};\bs{\delta}_1'',\gamma_1,\bs{\delta}_2',\xi_{2,3},\bs{\delta}_3'')\times\cM_{\Sigma_{23}}(\xi_{2,3};\bs{\delta}_2'',x_2,\bs{\delta}_3')\\
\end{alignat*}
where the three first unions are respectively for $p\in\Sigma_1\cap\Sigma_2$, $q\in\Sigma_2\cap\Sigma_3$ and $r\in\Sigma_1\cap\Sigma_3$.
\begin{figure}[ht]  
      \begin{center}\includegraphics[width=13cm]{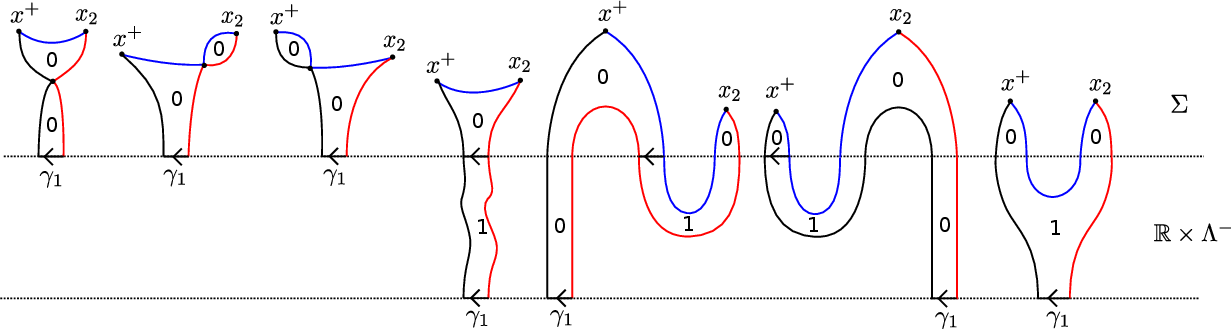}\end{center}
      \caption{Pseudo-holomorphic buildings in the boundary of the compactification of $\cM^1_{\Sigma_{123}}(x^+;\bs{\delta}_1,\gamma_1,\bs{\delta}_2,x_2,\bs{\delta}_3)$.}
      \label{degexc1}
\end{figure}
On Figure \ref{degexc1}, each broken configuration contributes from left to right to the terms of Relation \eqref{rel0xc}, and so we get:
\begin{alignat*}{1}
 &m_{00}^0(x_2,d_{0-}(\gamma_1))+m_{0-}^0(d_{00}(x_2),\gamma_1)+d_{00}\circ m_{0-}^0(x_2,\gamma_1)\\
 &\hspace{2cm}+m_{0-}^0(x_2,d_{--}(\gamma_1))+m_{--}^0(d_{-0}(x_2),\gamma_1)+d_{0-}\circ m_{0-}^-(x_2,\gamma_1)=0
\end{alignat*}

\subsubsection{Relation \eqref{rel-xc}}\label{Rel-xc}
\begin{figure}[ht]
      \begin{center}\includegraphics[width=5cm]{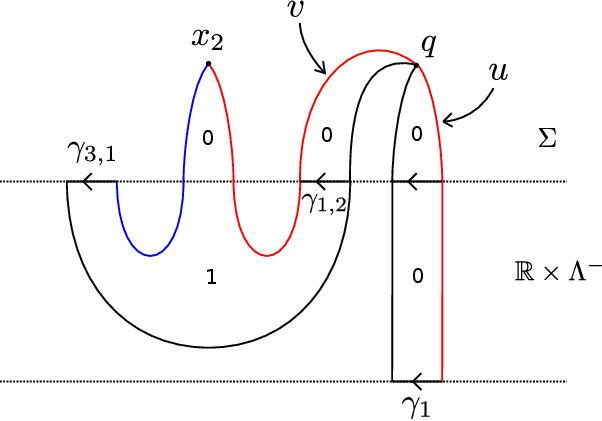}\end{center}
      \caption{Impossible breaking.}
      \label{brisure_imposs}
\end{figure}

Again, to find this relation we argue the same way as for Relation \eqref{rel-xx}. First, let us remark that one term in this relation already vanishes for energy reasons. More precisely, by definition we have:
\begin{alignat*}{1}
 & m_{00}^-(x_2,d_{0-}(\gamma_1))=b\circ f^{(2)}(x_2,d_{0-}(\gamma_1))+b^{(2)}(f^{(1)}(x_2),f^{(1)}\circ d_{0-}(\gamma_1))\\
\end{alignat*}
but $b^{(2)}(f^{(1)}(x_2),f^{(1)}\circ d_{0-}(\gamma_1))=0$ because such a term would count negative energy curves which is not possible, see Figure \ref{brisure_imposs}. Indeed, if there exist pseudo-holomorphic curves $u\in\cM^0(q;\bs{\delta}_1,\gamma_1,\bs{\delta}_2)$ and $v\in\cM^0(\gamma_{1,2};\bs{\zeta}_1,q,\bs{\zeta}_2)$, then the energies of these curves are given by (see Section \ref{action_energie}):
\begin{alignat*}{1}
 & E_{d(\chi\alpha)}(u)=\mathfrak{a}(q)-\mathfrak{a}(\gamma_1)-\mathfrak{a}(\bs{\delta}_1)-\mathfrak{a}(\bs{\delta}_2)\\
 & E_{d(\chi\alpha)}(v)=-\mathfrak{a}(q)-\mathfrak{a}(\gamma_{1,2})-\mathfrak{a}(\bs{\zeta}_1)-\mathfrak{a}(\bs{\zeta}_2)
\end{alignat*}
The energy of a non-constant pseudo-holomorphic curve is always strictly positive and the action of Reeb chords is also always positive, so the existence of $v$ implies that $q$ is an intersection point with a strictly negative action, which then contradicts the existence of $u$.
The other terms of Relation \eqref{rel-xc} are defined by a count of unfinished buildings of height $0|1|0$, $1|1|0$ or $2|1|0$. In each case, either the curves in the middle level form a pseudo-holomorphic building, or the curves in the middle level and on floor $-1$, or the curves on floors $-1$ and $-2$ form a building. In any case, after gluing, we get unfinished buildings in the following products of moduli spaces
\begin{alignat*}{1}
 & \widetilde{\cM^1}(\gamma_{3,1};\bs{\xi}_1,\gamma_{1,3},\bs{\xi}_3)\times\cM^1(\gamma_{1,3};\bs{\delta}_1,\gamma_1,\bs{\delta}_2,x_2,\bs{\delta}_3)\\
 & \widetilde{\cM^2}(\gamma_{3,1};\bs{\xi}_1,\gamma_{1,3},\bs{\xi}_3)\times\cM^0(\gamma_{1,3};\bs{\delta}_1,\gamma_1,\bs{\delta}_2,x_2,\bs{\delta}_3)\\
 & \widetilde{\cM^1}(\gamma_{3,1};\bs{\zeta}_1,\gamma_1,\bs{\zeta}_2,\gamma_{2,3},\bs{\zeta}_3)\times\cM^1(\gamma_{2,3};\bs{\delta}_2,x_2,\bs{\delta}_3)\\
 & \widetilde{\cM^2}(\gamma_{3,1};\bs{\zeta}_1,\gamma_1,\bs{\zeta}_2,\gamma_{2,3},\bs{\zeta}_3)\times\cM^0(\gamma_{2,3};\bs{\delta}_2,x_2,\bs{\delta}_3)
\end{alignat*}
Now, in order to get the relation, we have to find the broken curves in the boundary of the compactification of these products, i.e. broken curves in:
\begin{alignat}{2}
  & \widetilde{\cM^1}(\gamma_{3,1};\bs{\xi}_1,\gamma_{1,3},\bs{\xi}_3)\times\partial\overline{\cM^1}(\gamma_{1,3};\bs{\delta}_1,\gamma_1,\bs{\delta}_2,x_2,\bs{\delta}_3)\label{proddd1}\\
  & \partial\overline{\cM^2}(\gamma_{3,1};\bs{\xi}_1,\gamma_{1,3},\bs{\xi}_3)\times\cM^0(\gamma_{1,3};\bs{\delta}_1,\gamma_1,\bs{\delta}_2,x_2,\bs{\delta}_3)\label{proddd2}\\
  & \widetilde{\cM^1}(\gamma_{3,1};\bs{\zeta}_1,\gamma_1,\bs{\zeta}_2,\gamma_{2,3},\bs{\zeta}_3)\times\partial\overline{\cM^1}(\gamma_{2,3};\bs{\delta}_2,x_2,\bs{\delta}_3)\label{proddd3}\\
  & \partial\overline{\cM^2}(\gamma_{3,1};\bs{\zeta}_1,\gamma_1,\bs{\zeta}_2,\gamma_{2,3},\bs{\zeta}_3)\times\cM^0(\gamma_{2,3};\bs{\delta}_2,x_2,\bs{\delta}_3)\label{proddd4}
\end{alignat}
We already described $\partial\overline{\cM^2}(\gamma_{3,1};\bs{\xi}_1,\gamma_{1,3},\bs{\xi}_3)$ and $\partial\overline{\cM^1}(\gamma_{2,3};\bs{\delta}_2,x_2,\bs{\delta}_3)$ in Section \ref{Rel-xx}, so it remains to study $\partial\overline{\cM^1}(\gamma_{1,3};\bs{\delta}_1,\gamma_1,\bs{\delta}_2,x_2,\bs{\delta}_3)$ and $\partial\overline{\cM^2}(\gamma_{3,1};\bs{\zeta}_1,\gamma_1,\bs{\zeta}_2,\gamma_{2,3},\bs{\zeta}_3)$. First, we have the following decomposition:
\begin{alignat*}{2}
   &\partial\overline{\cM^1}(\gamma_{1,3};\bs{\delta}_1,\gamma_1,\bs{\delta}_2,x_2,\bs{\delta}_3)=\overline{\cM}^\partial(\gamma_{1,3};\bs{\delta}_1,\gamma_1,\bs{\delta}_2,x_2,\bs{\delta}_3)\\
   &\bigcup_{\substack{p\in\Sigma_1\cap\Sigma_2\\ \bs{\delta}_i',\bs{\delta}_i''}}\cM(\gamma_{1,3};\bs{\delta}_1',p,\bs{\delta}_2'',x_2,\bs{\delta}_3)\times\cM(p;\bs{\delta}_1'',\gamma_1,\bs{\delta}_2')\\
   &\bigcup_{\substack{q\in\Sigma_2\cap\Sigma_3\\ \bs{\delta}_i',\bs{\delta}_i''}}\cM(\gamma_{1,3};\bs{\delta}_1,\gamma_1,\bs{\delta}_2',q,\bs{\delta}_3'')\times\cM(q;\bs{\delta}_2'',x_2,\bs{\delta}_3')\\
   &\bigcup_{\substack{r\in\Sigma_1\cap\Sigma_3\\ \bs{\delta}_i',\bs{\delta}_i''}}\cM(\gamma_{1,3};\bs{\delta}_1',r,\bs{\delta}_3'')\times\cM(r;\bs{\delta}_1'',\gamma_1,\bs{\delta}_2,x_2,\bs{\delta}_3')\\
   &\bigcup_{\xi_{2,1},\bs{\delta}_i',\bs{\delta}_i''}\cM(\gamma_{1,3};\bs{\delta}_1',\xi_{2,1},\bs{\delta}_2'',x_2,\bs{\delta}_3)\times\widetilde{\cM}(\xi_{2,1};\bs{\delta}_1'',\gamma_1,\bs{\delta}_2')\\
   & \bigcup_{\xi_{2,3},\bs{\delta}_i',\bs{\delta}_i''}\widetilde{\cM}(\gamma_{1,3};\bs{\delta}_1,\gamma_1,\bs{\delta}_2',\xi_{2,3},\bs{\delta}_3'')\times\cM(\xi_{2,3};\bs{\delta}_2'',x_2,\bs{\delta}_3')\\
   & \bigcup_{\xi_{1,3},\bs{\delta}_i',\bs{\delta}_i''}\widetilde{\cM}(\gamma_{1,3};\bs{\delta}_1',\xi_{1,3},\bs{\delta}_3'')\times\cM(\xi_{1,3};\bs{\delta}_1'',\gamma_1,\bs{\delta}_2,x_2,\bs{\delta}_3')\\
\end{alignat*}
Finally, the buildings occurring as degeneration of index-$2$ bananas with two positive Reeb chord asymptotics and one negative one are of the following type:
\begin{alignat*}{2}
   & \partial\overline{\cM^2}(\gamma_{3,1};\bs{\zeta}_1,\gamma_1,\bs{\zeta}_2,\gamma_{2,3},\bs{\zeta}_3)=\overline{\cM}^\partial(\gamma_{3,1};\bs{\zeta}_1,\gamma_{1},\bs{\zeta}_2,\gamma_{2,3},\bs{\zeta}_3)\\
   & \bigcup_{\xi_{2,1},\bs{\zeta}_i',\bs{\zeta}_i''}\widetilde{\cM}(\gamma_{3,1};\bs{\zeta}_1',\xi_{2,1},\bs{\zeta}_2'',\gamma_{2,3},\bs{\zeta}_3)\times\widetilde{\cM}(\xi_{2,1};\bs{\zeta}_1'',\gamma_{1},\bs{\zeta}_2')\\
   & \bigcup_{\xi_{2,3},\bs{\zeta}_i',\bs{\zeta}_i''}\widetilde{\cM}(\gamma_{3,1};\bs{\zeta}_1,\gamma_{1},\bs{\zeta}_2',\xi_{2,3},\bs{\zeta}_3'')\times\widetilde{\cM}(\xi_{2,3};\bs{\zeta}_2'',\gamma_{2,3},\bs{\zeta}_3')\\
   & \bigcup_{\xi_{3,2},\bs{\zeta}_i',\bs{\zeta}_i''}\widetilde{\cM}(\gamma_{3,1};\bs{\zeta}_1,\gamma_{1},\bs{\zeta}_2',\xi_{3,2},\bs{\zeta}_3'')\times\widetilde{\cM}(\xi_{3,2};\bs{\zeta}_2'',\gamma_{2,3},\bs{\zeta}_3')\\
   & \bigcup_{\xi_{3,1},\bs{\zeta}_i',\bs{\zeta}_i''}\widetilde{\cM}(\gamma_{3,1};\bs{\zeta}_1',\xi_{3,1},\bs{\zeta}_3'')\times\widetilde{\cM}(\xi_{3,1};\bs{\zeta}_1'',\gamma_{1},\bs{\zeta}_2,\gamma_{2,3},\bs{\zeta}_3')\\
   & \bigcup_{\xi_{1,3},\bs{\zeta}_i',\bs{\zeta}_i''}\widetilde{\cM}(\gamma_{3,1};\bs{\zeta}_1',\xi_{1,3},\bs{\zeta}_3'')\times\widetilde{\cM}(\xi_{1,3};\bs{\zeta}_1'',\gamma_{1},\bs{\zeta}_2,\gamma_{2,3},\bs{\zeta}_3')
\end{alignat*}
\begin{figure}[ht] 
	\begin{center}\includegraphics[width=13cm]{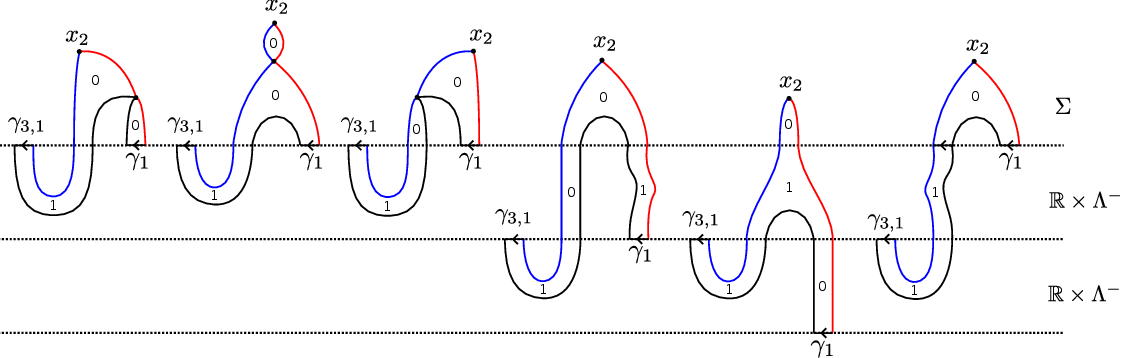}\end{center}
	\caption{Broken curves in \eqref{proddd1}.}
	\label{dege25}
\end{figure}
\begin{figure}[ht] 
	\begin{center}\includegraphics[width=12cm]{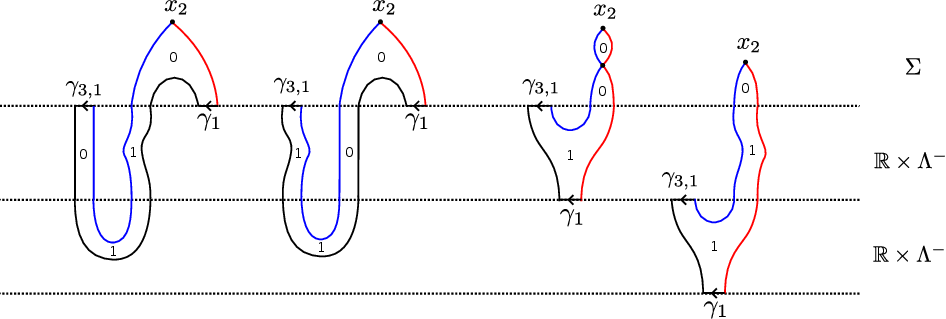}\end{center}
	\caption{Broken curves in \eqref{proddd2} and \eqref{proddd3}.}
	\label{dege2627}
\end{figure}
\begin{figure}[ht] 
	\begin{center}\includegraphics[width=12cm]{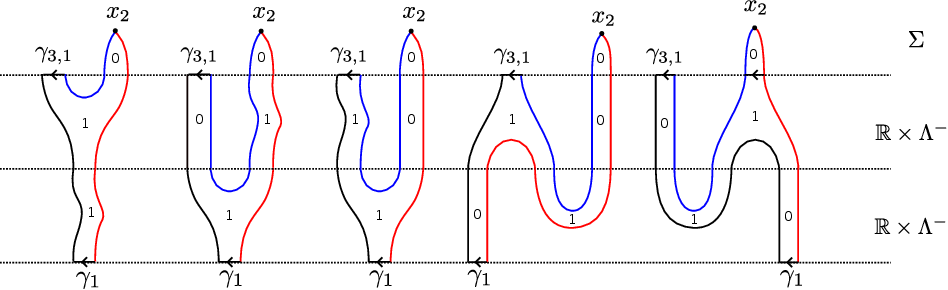}\end{center}
	\caption{Broken curves in \eqref{proddd4}.}
	\label{dege28}
\end{figure}

The different types of unfinished buildings corresponding to elements in the products \eqref{proddd1}, \eqref{proddd2}, \eqref{proddd3} and \eqref{proddd4} are schematized on Figures \ref{dege25}, \ref{dege2627} and \ref{dege28}. Observe that again some unfinished buildings are equivalent (see Remark \ref{rem_bord}). Indeed, the fifth configuration (from left to right) on Figure \ref{dege25} and the last configuration on Figure \ref{dege28} contribute both algebraically to $b\circ\Delta^{(2)}(f^{(1)}(x_2),\gamma_1)$. Then the last configuration on Figure \ref{dege25} and the first on Figure \ref{dege2627} contribute both algebraically to $b\circ\delta_{--}\circ f(x_2,\gamma_1)$. Finally, the last unfinished building on Figure \ref{dege2627} and the second one on Figure \ref{dege28} contribute to $b^{(2)}(\delta_{--}\circ f^{(1)}(x_2),\gamma_1)$.
Summing the contributions of the remaining unfinished buildings gives the Relation \eqref{rel-xc}:
\begin{alignat*}{2}
 & m_{00}^-(x_2,d_{0-}(\gamma_1))+m_{0-}^-(d_{00}(x_2),\gamma_1)+m_{0-}^-(x_2,d_{--}(\gamma_1))\\
 &\hspace{2cm} +m_{--}^-(d_{-0}(x_2),\gamma_1)+d_{-0}\circ m_{0-}^0(x_2,\gamma_1)+d_{--}\circ m_{0-}^-(x_2,\gamma_1)=0
\end{alignat*}

\subsubsection{Relations \eqref{rel0cx} and \eqref{rel-cx}}
By symmetry, Relations \eqref{rel0cx} and \eqref{rel-cx} for a pair $(\gamma_2,x_1)$ are obtained by studying same types of holomorphic curves as for Relations \eqref{rel0xc} and \eqref{rel-xc} corresponding to a pair of asymptotics $(x_2,\gamma_1)$. 
\begin{figure} 
      \begin{center}\includegraphics[width=125mm]{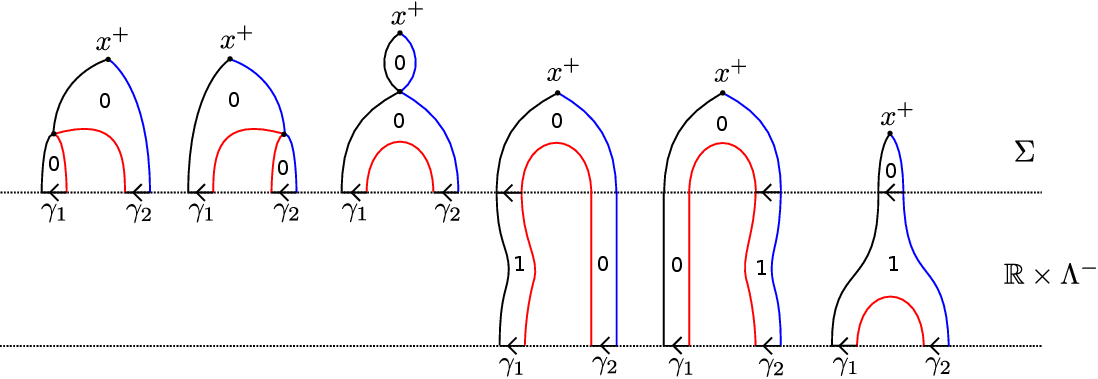}\end{center}
      \caption{Pseudo-holomorphic buildings in $\partial\overline{\cM^1}(x^+;\bs{\delta}_1,\gamma_1,\bs{\delta}_2,\gamma_2,\bs{\delta}_3)$.}
      \label{degecc1}
\end{figure}

\subsubsection{Relation \eqref{rel0cc}}
Each term of this relation corresponds to a count of broken curves in the boundary of the compactification of $\cM^1(x^+;\bs{\delta}_1,\gamma_1,\bs{\delta}_2,\gamma_2,\bs{\delta}_3)$, and we have (see Figure \ref{degecc1}):
\begin{alignat*}{1}
  &\partial\overline{\cM^1}(x^+;\bs{\delta}_1,\gamma_1,\bs{\delta}_2,\gamma_2,\bs{\delta}_3)=\overline{\cM}^\partial(x^+;\bs{\delta}_1,\gamma_1,\bs{\delta}_2,\gamma_2,\bs{\delta}_3)\\
  &\bigcup\limits_{\substack{q\in\Sigma_1\cap\Sigma_3\\ \bs{\delta}_i',\bs{\delta}_i''}}\cM^0(x^+;\bs{\delta}_1',q,\bs{\delta}_3'')\times\cM^0(q;\bs{\delta}_1'',\gamma_1,\bs{\delta}_2,\gamma_2,\bs{\delta}_3')\\
  &\bigcup\limits_{\substack{q\in\Sigma_1\cap\Sigma_2\\ \bs{\delta}_i',\bs{\delta}_i''}}\cM^0(x^+;\bs{\delta}_1',q,\bs{\delta}_2'',\gamma_2,\bs{\delta}_3)\times\cM^0(q;\bs{\delta}_1'',\gamma_1,\bs{\delta}_2')\\
  &\bigcup\limits_{\substack{q\in\Sigma_2\cap\Sigma_3\\ \bs{\delta}_i',\bs{\delta}_i''}}\cM^0(x^+;\bs{\delta}_1,\gamma_1,\bs{\delta}_2',q,\bs{\delta}_3'')\times\cM^0(q;\bs{\delta}_2'',\gamma_2,\bs{\delta}_3')\\
  &\bigcup\limits_{\xi_{3,1},\bs{\delta}_i',\bs{\delta}_i''}\cM^0(x^+;\bs{\delta}_1',\xi_{3,1},\bs{\delta}_3'')\times\widetilde{\cM^1}(\xi_{3,1};\bs{\delta}_1'',\gamma_1,\bs{\delta}_2,\gamma_2,\bs{\delta}_3')\\
  &\bigcup\limits_{\xi_{2,1},\bs{\delta}_i',\bs{\delta}_i''}\cM^0(x^+;\bs{\delta}_1',\xi_{2,1},\bs{\delta}_2'',\gamma_2,\bs{\delta}_3)\times\widetilde{\cM^1}(\xi_{2,1};\bs{\delta}_1'',\gamma_1,\bs{\delta}_2')\\
  &\bigcup\limits_{\xi_{3,2},\bs{\delta}_i',\bs{\delta}_i''}\cM^0(x^+;\bs{\delta}_1,\gamma_1,\bs{\delta}_2',\xi_{3,2},\bs{\delta}_3'')\times\widetilde{\cM^1}(\xi_{3,2};\bs{\delta}_2'',\gamma_2,\bs{\delta}_3')\\  
\end{alignat*}

\subsubsection{Relation \eqref{rel-cc}}
\begin{figure}[ht] 
      \begin{center}\includegraphics[width=75mm]{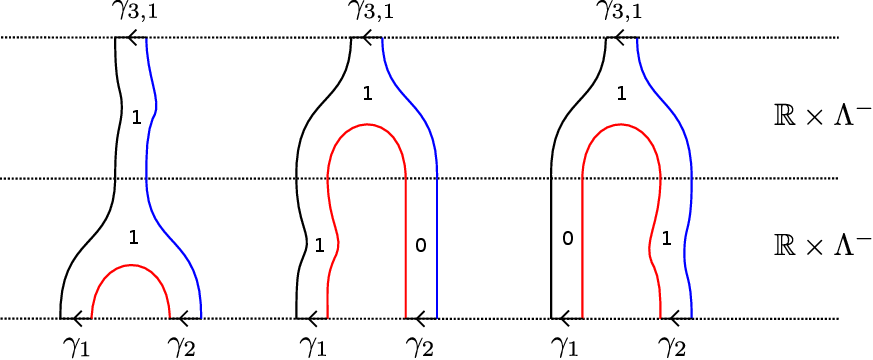}\end{center}
      \caption{Pseudo-holomorphic buildings in the boundary of $\overline{\cM^2}(\gamma_{3,1};\bs{\delta}_1,\gamma_1,\bs{\delta}_2,\gamma_2,\bs{\delta}_3)$}
      \label{degecc2}
\end{figure}

The product of two Reeb chords being given by the product $\mu^2_{\ep_{3,2,1}^+}$ in the augmentation category $\Aug_-(\La_1^-\cup\La_2^-\cup\La_2^-)$, Relation \eqref{rel-cc} is satisfied because it is the $A_\infty$-relation for $d=2$ satisfied by the maps $\{\mu^d\}_{d\geq 1}$ (see \eqref{relAug} in Section \ref{linearisation}). We recall the different kinds of degeneration of a curve in $\widetilde{\cM^2}(\gamma_{3,1};\bs{\delta}_1,\gamma_1,\bs{\delta}_2,\gamma_2,\bs{\delta}_3)$ (see Figure \ref{degecc2}):
\begin{alignat*}{2}
 & \partial\overline{\cM^2}(\gamma_{3,1};\bs{\delta}_1,\gamma_1,\bs{\delta}_2,\gamma_2,\bs{\delta}_3)=\overline{\cM}^\partial(\gamma_{3,1};\bs{\delta}_1,\gamma_1,\bs{\delta}_2,\gamma_2,\bs{\delta}_3)\\
 &\bigcup\widetilde{\cM}(\gamma_{3,1};\bs{\delta}_1',\xi_{3,1},\bs{\delta}_3'')\times\widetilde{\cM}(\xi_{3,1};\bs{\delta}_1'',\gamma_1,\bs{\delta}_2,\gamma_2,\bs{\delta}_3')\\
 & \bigcup\widetilde{\cM}(\gamma_{3,1};\bs{\delta}_1',\xi_{2,1},\bs{\delta}_2'',\gamma_2,\bs{\delta}_3)\times\widetilde{\cM}(\xi_{2,1};\bs{\delta}_1'',\gamma_1,\bs{\delta}_2')\\
 & \bigcup\widetilde{\cM}(\gamma_{3,1};\bs{\delta}_1,\gamma_1,\bs{\delta}_2',\xi_{3,2},\bs{\delta}_3'')\times\widetilde{\cM}(\xi_{3,2};\bs{\delta}_2'',\gamma_2,\bs{\delta}_3')\\
\end{alignat*}

\subsection{Proof of Theorem \ref{teofoncteur} and Corollary \ref{isoES}}\label{teo2}

\begin{figure}[ht]  
      \begin{center}\includegraphics[width=9cm]{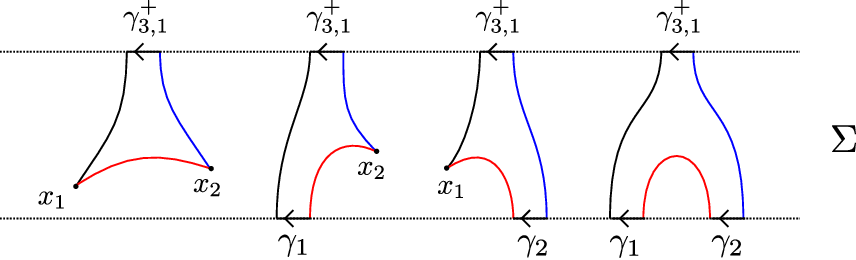}\end{center}
      \caption{Examples of curves in the moduli spaces \eqref{psixx}, \eqref{psixc}, \eqref{psicx}, and \eqref{psicc} respectively.}
      \label{phi2}
\end{figure}

Recall that $\Fc^1:=d_{+0}+d_{+-}$ and that the acyclicity of the Cthulhu complex implies that $\Fc^1$ is a quasi-isomorphism (Section \ref{cthulhu}). We consider as before three transverse exact Lagrangian cobordisms $\Sigma_1,\Sigma_2$ and $\Sigma_3$ such that the algebras $\Ac(\La_i^-)$ admit augmentations. As introduced in Section \ref{moduli_spaces}, we need to consider the following moduli spaces of curves with boundary on $\Sigma_1\cup\Sigma_2\cup\Sigma_3$ (see Figure \ref{phi2}):
  \begin{alignat}{1}
    & \cM_{\Sigma_{123}}(\gamma_{3,1}^+;\bs{\delta}_1,x_1,\bs{\delta}_2,x_2,\bs{\delta}_3) \label{psixx}\\
    & \cM_{\Sigma_{123}}(\gamma_{3,1}^+;\bs{\delta}_1,\gamma_1,\bs{\delta}_2,x_2,\bs{\delta}_3) \label{psixc} \\
    & \cM_{\Sigma_{123}}(\gamma_{3,1}^+;\bs{\delta}_1,x_1,\bs{\delta}_2,\gamma_2,\bs{\delta}_3) \label{psicx}\\
    & \cM_{\Sigma_{123}}(\gamma_{3,1}^+;\bs{\delta}_1,\gamma_1,\bs{\delta}_2,\gamma_2,\bs{\delta}_3) \label{psicc}
  \end{alignat}
with
\begin{itemize}
 \item[$\bullet$] $\gamma_{3,1}^+\in\Rc(\La_3^+,\La_1^+)$, $\gamma_1\in\Rc(\La_2^-,\La_1^-)$ and $\gamma_2\in\Rc(\La_3^-,\La_2^-)$,
 \item[$\bullet$] $x_1\in\Sigma_1\cap\Sigma_2$, $x_2\in\Sigma_2\cap\Sigma_3$,
 \item[$\bullet$] $\bs{\delta}_i$ are words of Reeb chords of $\La_i^-$, for $i=1,2,3$.
\end{itemize}
By a count of rigid pseudo-holomorphic disks in these moduli spaces, we introduce a map:
  \begin{alignat*}{1}
   \Fc^2\colon CF_{-\infty}(\Sigma_2,\Sigma_3)\otimes CF_{-\infty}(\Sigma_1,\Sigma_2)\to C^*(\La_1^+,\La_3^+)
  \end{alignat*}
defined on pairs of generators by:
  \begin{alignat*}{1}
     & \Fc^2(x_2,x_1)=\sum\limits_{\gamma_{3,1}^+,\bs{\delta}_i}\#\cM^0(\gamma_{3,1}^+;\bs{\delta}_1,x_1,\bs{\delta}_2,x_2,\bs{\delta}_3)\ep_1^-(\bs{\delta}_1)\ep_2^-(\bs{\delta}_2)\ep_3^-(\bs{\delta}_3)\cdot\gamma_{3,1}^+\\
     & \Fc^2(x_2,\gamma_1)=\sum\limits_{\gamma_{3,1}^+,\bs{\delta}_i}\#\cM^0(\gamma_{3,1}^+;\bs{\delta}_1,\gamma_1,\bs{\delta}_2,x_2,\bs{\delta}_3)\ep_1^-(\bs{\delta}_1)\ep_2^-(\bs{\delta}_2)\ep_3^-(\bs{\delta}_3)\cdot\gamma_{3,1}^+\\
     & \Fc^2(\gamma_2,x_1)=\sum\limits_{\gamma_{3,1}^+,\bs{\delta}_i}\#\cM^0(\gamma_{3,1}^+;\bs{\delta}_1,x_1,\bs{\delta}_2,\gamma_2,\bs{\delta}_3)\ep_1^-(\bs{\delta}_1)\ep_2^-(\bs{\delta}_2)\ep_3^-(\bs{\delta}_3)\cdot\gamma_{3,1}^+\\
     & \Fc^2(\gamma_2,\gamma_1)=\sum\limits_{\gamma_{3,1}^+,\bs{\delta}_i}\#\cM^0(\gamma_{3,1}^+;\bs{\delta}_1,\gamma_1,\bs{\delta}_2,\gamma_2,\bs{\delta}_3)\ep_1^-(\bs{\delta}_1)\ep_2^-(\bs{\delta}_2)\ep_3^-(\bs{\delta}_3)\cdot\gamma_{3,1}^+\\
  \end{alignat*} 
We have again to study breakings of index-$1$ pseudo-holomorphic curves in the moduli spaces \eqref{psixx}, \eqref{psixc}, \eqref{psicx}, and \eqref{psicc} in order to prove Theorem \ref{teofoncteur}.
\begin{figure}[ht]  
      \begin{center}\includegraphics[width=145mm]{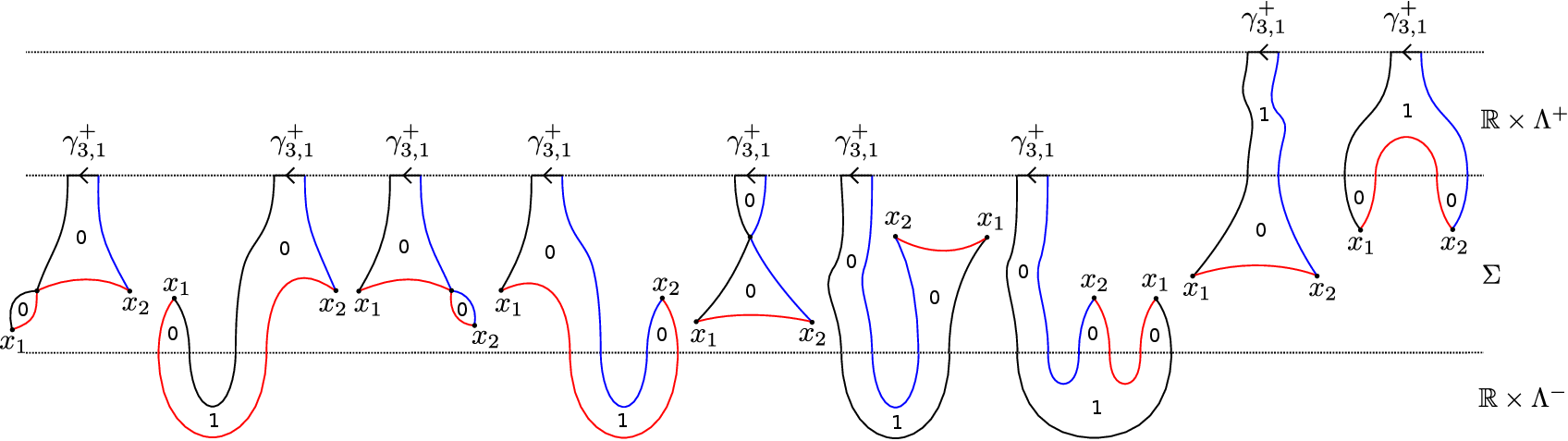}\end{center}
      \caption{Pseudo-holomorphic buildings in $\partial\overline{\cM^1}(\gamma_{3,1}^+;\bs{\delta}_1,x_1,\bs{\delta}_2,x_2,\bs{\delta}_3)$.}
      \label{degene psixx}
\end{figure}
For example, we describe below the boundary of the compactification of $\cM^1(\gamma_{3,1}^+;\bs{\delta}_1,x_1,\bs{\delta}_2,x_2,\bs{\delta}_3)$ (see Figure \ref{degene psixx}):
 \begin{alignat*}{1}
     & \partial\overline{\cM^1}(\gamma_{3,1}^+;\bs{\delta}_1,x_1,\bs{\delta}_2,x_2,\bs{\delta}_3)=\overline{\cM^1}^\partial(\gamma_{3,1}^+;\bs{\delta}_1,x_1,\bs{\delta}_2,x_2,\bs{\delta}_3)\\
     & \bigcup_{p\in\Sigma_1\cap\Sigma_2}\cM(\gamma_{3,1}^+;\bs{\delta}_1',p,\bs{\delta}_2'',x_2,\bs{\delta}_3)\times\cM(p;\bs{\delta}_1'',x_1,\bs{\delta}_2')\\
     & \bigcup_{q\in\Sigma_2\cap\Sigma_3}\cM(\gamma_{3,1}^+;\bs{\delta}_1,x_1,\bs{\delta}_2',q,\bs{\delta}_3'')\times\cM(q;\bs{\delta}_2'',x_2,\bs{\delta}_3')\\
     & \bigcup_{r\in\Sigma_1\cap\Sigma_3}\cM(\gamma_{3,1}^+;\bs{\delta}_1',r,\bs{\delta}_3'')\times\cM(r;\bs{\delta}_1'',x_1,\bs{\delta}_2,x_2,\bs{\delta}_3')\\
     & \bigcup_{\xi_{2,1}^-,\xi_{1,2}^-}\cM(\gamma_{3,1}^+;\bs{\delta}_1',\xi_{2,1}^-,\bs{\delta}_2''',x_2,\bs{\delta}_3)\times\widetilde{\cM}(\xi_{2,1}^-;\bs{\delta}_1'',\xi_{1,2}^-,\bs{\delta}_2'')\times\cM_{\Sigma_{12}}(\xi_{1,2}^-;\bs{\delta}_1''',x_1,\bs{\delta}_2')\\
     & \bigcup_{\xi_{3,2}^-,\xi_{2,3}^-}\cM(\gamma_{3,1}^+;\bs{\delta}_1,x_1,\bs{\delta}_2',\xi_{3,2}^-,\bs{\delta}_3''')\times\widetilde{\cM}(\xi_{3,2}^-;\bs{\delta}_2'',\xi_{2,3}^-,\bs{\delta}_3'')\times\cM_{\Sigma_{23}}(\xi_{2,3}^-;\bs{\delta}_2''',x_2,\bs{\delta}_3')\\
     & \bigcup_{\xi_{3,1}^-,\xi_{1,3}^-}\cM(\gamma_{3,1}^+;\bs{\delta}_1',\xi_{3,1}^-,\bs{\delta}_3''')\times\widetilde{\cM}(\xi_{3,1}^-;\bs{\delta}_1'',\xi_{1,3}^-,\bs{\delta}_3'')\times\cM_{\Sigma_{123}}(\xi_{1,3}^-;\bs{\delta}_1''',x_1,\bs{\delta}_2,x_2,\bs{\delta}_3')\\
     &\bigcup_{\xi_{3,1}^-,\xi_{1,2}^-,\xi_{2,3}^-}\cM(\gamma_{3,1}^+;\bs{\delta}_1',\xi_{3,1}^-,\bs{\delta}_3''')\times\widetilde{\cM}(\xi_{3,1}^-;\bs{\delta}_1'',\xi_{1,2}^-,\bs{\delta}_2'',\xi_{2,3}^-,\bs{\delta}_3'')\\
     &\hspace{65mm}\times\cM_{\Sigma_{12}}(\xi_{1,2}^-;\bs{\delta}_1''',x_1,\bs{\delta}_2')\times\cM_{\Sigma_{23}}(\xi_{2,3}^-;\bs{\delta}_2''',x_2,\bs{\delta}_3')\\
     &\\
     & \bigcup_{\xi_{3,1}^+}\widetilde{\cM^1}(\gamma_{3,1}^+;\beta_{1,1},\dots,\beta_{1,s},\xi_{3,1}^+,\beta_{3,1},\dots,\beta_{3,t})\times\cM^0(\xi_{3,1}^+;\bs{\delta}_1'',x_1,\bs{\delta}_2,x_2,\bs{\delta}_3')\\
     &\hspace{75mm}\times\prod_{i=1}^s\cM^0_{\Sigma_1}(\beta_{1,i};\bs{\delta}_{1,i})\times\prod_{j=1}^t\cM^0_{\Sigma_3}(\beta_{3,j};\bs{\delta}_{3,j})\\
     & \bigcup_{\xi_{2,1}^+,\xi_{3,2}^+}\widetilde{\cM^1}(\gamma_{3,1}^+;\beta_{1,1}\dots,\beta_{1,l},\xi_{2,1}^+,\beta_{2,1},\dots,\beta_{2,m},\xi_{3,2}^+,\beta_{3,1},\dots,\beta_{3,n})\times\cM^0(\xi_{2,1}^+;\bs{\delta}_1'',x_1,\bs{\delta}_2')\\
     &\hspace{1cm}\times\cM^0(\xi_{3,2}^+;\bs{\delta}_2''',x_2,\bs{\delta}_3')\times\prod_{i=1}^l\cM^0_{\Sigma_1}(\beta_{1,i};\bs{\delta}_{1,i})\times\prod_{i=1}^m\cM^0_{\Sigma_2}(\beta_{2,i};\bs{\delta}_{2,i})\times\prod_{i=1}^n\cM^0_{\Sigma_3}(\beta_{3,i};\bs{\delta}_{3,i})
  \end{alignat*}
where all unions except the last two are also for words of pure Reeb chords $\bs{\delta}_i',\bs{\delta}_i''$ and $\bs{\delta}_i'''$ of $\La_i^-$ such that $\bs{\delta}_i',\bs{\delta}_i''=\bs{\delta}_i$, or $\bs{\delta}_i'\bs{\delta}_i''\bs{\delta}_i'''=\bs{\delta}_i$, depending on cases. The second to last union is for:
\begin{itemize}
  \item $\beta_{1,i}\in\Rc(\La_1^+)$ for $1\leq i\leq s$,
  \item $\beta_{3,j}\in\Rc(\La_3^+)$ for $1\leq j\leq t$,
  \item $\bs{\delta}_1''$ and $\bs{\delta}_{1,i}$ for $1\leq i\leq s$ words of Reeb chords of $\La_1^-$,
  \item $\bs{\delta}_3'$ and $\bs{\delta}_{3,j}$ for $1\leq j\leq t$ words of Reeb chords of $\La_3^-$,
\end{itemize}
such that $\bs{\delta}_{1,1}\dots\bs{\delta}_{1,s}\bs{\delta}_1''=\bs{\delta}_1$ and $\bs{\delta}_3'\bs{\delta}_{3,1}\dots\bs{\delta}_{3,t}=\bs{\delta}_3$. The count of curves in the moduli space $\cM^0_{\Sigma_1}(\beta_{1,i};\bs{\delta}_{1,i})$ contributes to the coefficient $\langle\phi_{\Sigma_1}(\beta_{1,i}),\bs{\delta}_{1,i}\rangle$, with $\phi_{\Sigma_1}:\Ac(\La_1^+)\to\Ac(\La_1^-)$ is the chain map induced by $\Sigma_1$ (see Section \ref{morphisme_induit}). So the count of curves in
\begin{alignat*}{1}
 \cM^1(\gamma_{3,1}^+;\beta_{1,1},\dots,\beta_{1,s},\xi_{3,1}^+,\beta_{3,1},\dots,\beta_{3,t})\times\prod_{i=1}^s\cM^0_{\Sigma_1}(\beta_{1,i};\bs{\delta}_{1,i})\times\prod_{j=1}^t\cM^0_{\Sigma_3}(\beta_{3,j};\bs{\delta}_{3,j})
\end{alignat*}
contributes to $\langle\mu^1_{\ep_3^+,\ep_1^+}(\xi_{3,1}^+),\gamma^+_{3,1}\rangle$ because $\ep_i^+=\ep_i^-\circ\phi_{\Sigma_i}$. Finally, the last union is for
\begin{itemize}
  \item $\beta_{1,i}\in\Rc(\La_1^+)$ for $1\leq i\leq l$,
  \item $\beta_{2,i}\in\Rc(\La_2^+)$ for $1\leq i\leq m$,
  \item $\beta_{3,i}\in\Rc(\La_3^+)$ for $1\leq i\leq n$,
  \item $\bs{\delta}_1''$ and $\bs{\delta}_{1,i}$ for $1\leq i\leq l$ words of Reeb chords of $\La_1^-$,
  \item $\bs{\delta}_2'$, $\bs{\delta}_2'''$ and $\bs{\delta}_{2,i}$ for $1\leq i\leq m$ words of Reeb chords of $\La_2^-$,
  \item $\bs{\delta}_3'$ and $\bs{\delta}_{3,i}$ for $1\leq i\leq n$ words of Reeb chords of $\La_3^-$,
\end{itemize}
such that
\begin{itemize}
 \item $\bs{\delta}_{1,1}\dots\bs{\delta}_{1,l}\bs{\delta}_1''=\bs{\delta}_1$,
 \item $\bs{\delta}_2'\bs{\delta}_{2,1}\dots\bs{\delta}_{2,m}\bs{\delta}_2'''=\bs{\delta}_2$, and
 \item $\bs{\delta}_3'\bs{\delta}_{3,1}\dots\bs{\delta}_{3,n}=\bs{\delta}_3$.
\end{itemize}
 Again, the count of broken curves in 
\begin{alignat*}{1}
 \cM(\gamma_{3,1}^+;\beta_{1,1}\dots,\beta_{1,l},\xi_{2,1}^+,\beta_{2,1},\dots,\beta_{2,m},\xi_{3,2}^+,\beta_{3,1},\dots,\beta_{3,n})\times\prod_{i=1}^l\cM^0_{\Sigma_1}(\beta_{1,i};\bs{\delta}_{1,i})\\
 \times\prod_{i=1}^m\cM^0_{\Sigma_2}(\beta_{2,i};\bs{\delta}_{2,i})\times\prod_{i=1}^n\cM^0_{\Sigma_3}(\beta_{3,i};\bs{\delta}_{3,i})
\end{alignat*}
contributes to the coefficient $\langle\mu^2_{\ep_3^+,\ep_2^+,\ep_1^+}(\xi_{3,2}^+,\xi_{2,1}^+),\gamma_{3,1}^+\rangle$.

\begin{figure}[ht]  
      \begin{center}\includegraphics[width=135mm]{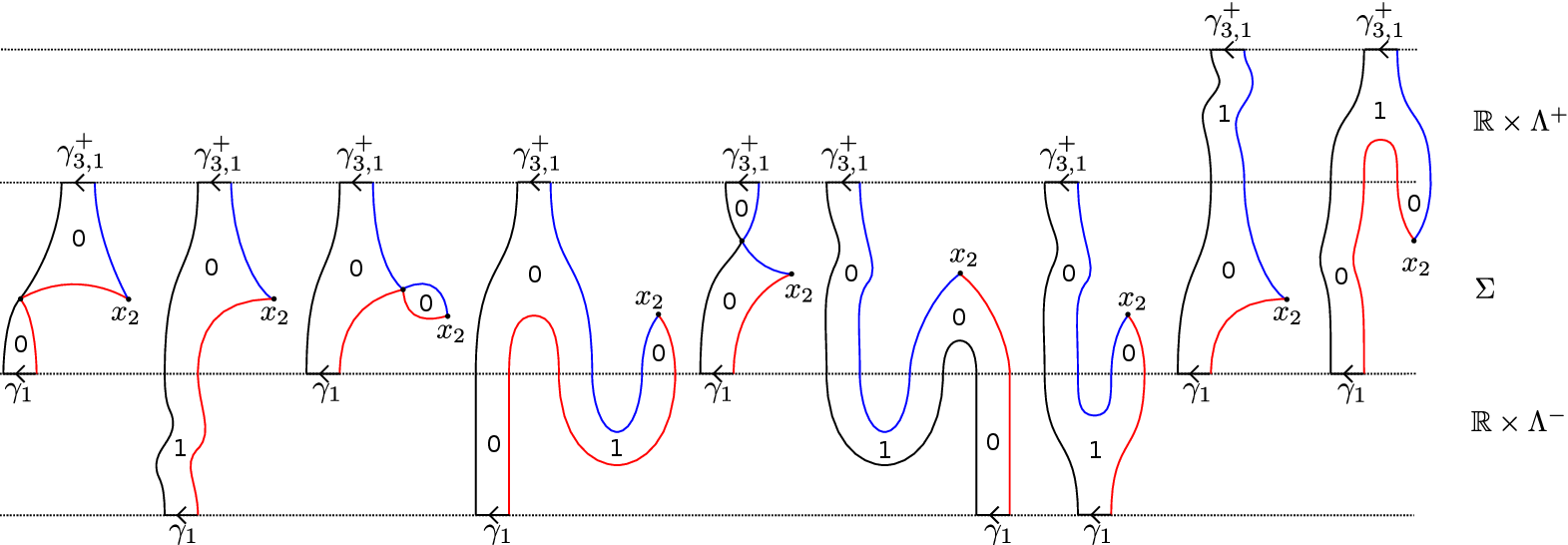}\end{center}
      \caption{Pseudo-holomorphic buildings in $\partial\overline{\cM^1}(\gamma_{3,1}^+;\bs{\delta}_1,\gamma_1,\bs{\delta}_2,x_2,\bs{\delta}_3)$.}
      \label{degene psixc}
\end{figure}

\begin{figure}[ht]  
      \begin{center}\includegraphics[width=12cm]{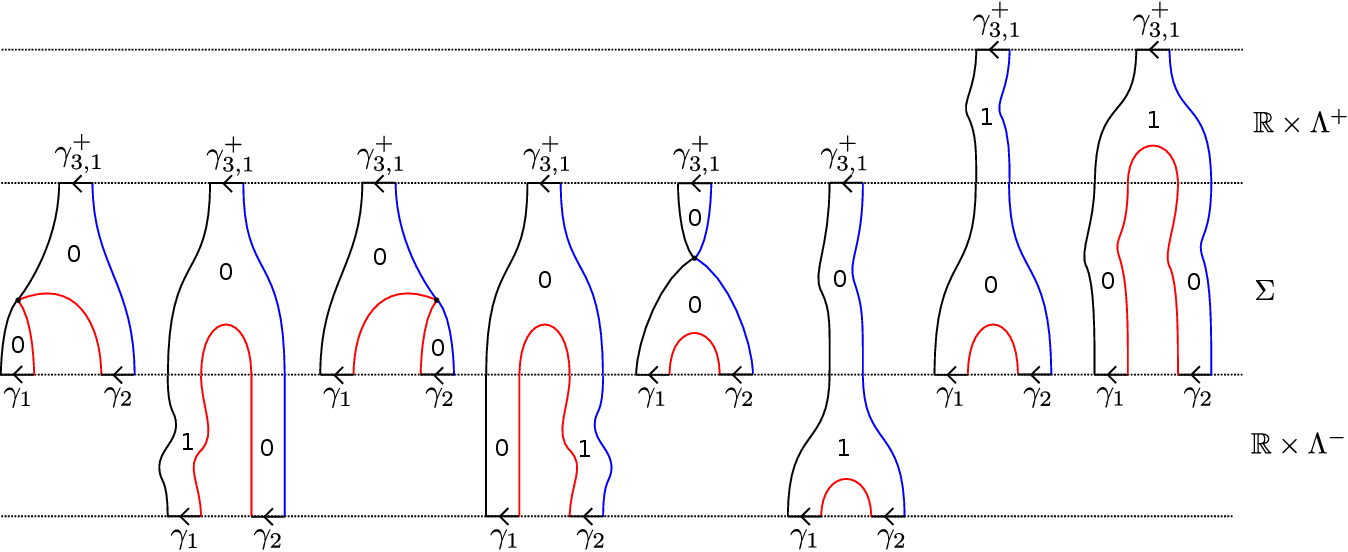}\end{center}
      \caption{Pseudo-holomorphic buildings in $\partial\overline{\cM^1}(\gamma_{3,1}^+;\bs{\delta}_1,\gamma_1,\bs{\delta}_2,\gamma_2,\bs{\delta}_3)$.}
      \label{degene psicc}
\end{figure}

By denoting $d_{++}$ for $\mu^1_{\ep_3^+,\ep_1^+}$, the study of breakings above implies that $\Fc^2$ satisfies the relation:
    \begin{alignat}{1}
      \Fc^2(x_2,d_{00}(x_1))+\Fc^2(x_2,d_{0-}(x_1))+\Fc^2(d_{00}(x_2),x_1)+\Fc^2(d_{-0}(x_2),x_1)+d_{+0}\circ m_{00}^0(x_2,x_1) \nonumber\\
      +d_{+-}\circ m_{00}^-(x_2,x_1)+d_{++}\circ\Fc^2(x_2,x_1)+\mu_{\ep_{3,2,1}^+}^2(d_{+0}(x_2),d_{+0}(x_1))=0 \label{degepsixx}
    \end{alignat}  
Analogously, the different types of buildings in $\partial\overline{\cM^1}_{\Sigma_{123}}(\gamma_{3,1}^+;\bs{\delta}_1,\gamma_1,\bs{\delta}_2,x_2,\bs{\delta}_3)$ are schematized on Figure \ref{degene psixc}, and this gives for $\Fc^2$ the relation:
    \begin{alignat}{1}
      \Fc^2(x_2,d_{0-}(\gamma_1))+\Fc^2(x_2,d_{--}(\gamma_1))+\Fc^2(d_{00}(x_2),\gamma_1)+\Fc^2(d_{-0}(x_2),\gamma_1)+d_{+0}\circ m_{0-}^0(x_2,\gamma_1)\nonumber \\
      +d_{+-}\circ m_{0-}^-(x_2,\gamma_1) +d_{++}\circ\Fc^2(x_2,\gamma_1)+\mu^2_{\ep_{3,2,1}^+}(d_{-0}(x_2),d_{+-}(\gamma_1))=0\label{degepsixc}
   \end{alignat}   
The symmetric relation for the pair $(\gamma_2,x_1)$ of asymptotics is of course also satisfied. Finally, pseudo-holomorphic buildings in $\partial\overline{\cM^1}_{\Sigma_{123}}(\gamma_{3,1}^+;\bs{\delta}_1,\gamma_1,\bs{\delta}_2,\gamma_2,\bs{\delta}_3)$ are schematized on Figure \ref{degene psicc} and thus we get:
\begin{alignat}{1}
      \Fc^2(\gamma_2,d_{0-}(\gamma_1))+\Fc^2(\gamma_2,d_{--}(\gamma_1))+\Fc^2(d_{0-}(\gamma_2),\gamma_1)+\Fc^2(d_{--}(\gamma_2),\gamma_1)+d_{+0}\circ m_{--}^0(\gamma_2,\gamma_1)\nonumber \\
      +d_{+-}\circ m_{--}^-(\gamma_2,\gamma_1)+d_{++}\circ\Fc^2(\gamma_2,\gamma_1)+\mu^2_{\ep_{3,2,1}^+}(d_{+-}(\gamma_2),d_{+-}(\gamma_1))=0 \label{degepsicc}
\end{alignat}
Combining Relations \eqref{degepsixx}, \eqref{degepsixc} and its symmetric one, and \eqref{degepsicc}, we deduce that $\Fc^2$ satisfies:
    \begin{alignat*}{1}
     \Fc^2(\cdot,\partial_{-\infty})+\Fc^2(\partial_{-\infty},\cdot)+\Fc_{31}^1\circ \mfm_2+d_{++}\circ\Fc^2+\mu^2_{\ep_{3,2,1}^+}(\Fc_{32}^1,\Fc_{21}^1)=0
    \end{alignat*}
The map induced by $\Fc^1$ in homology satisfies then $\Fc_{31}^1\circ \mfm_2+\mu^2_{\ep_{3,2,1}^+}(\Fc_{32}^1,\Fc_{21}^1)=0$, and so $\Fc^1$ preserves products in homology.
This concludes the proof of Theorem \ref{teofoncteur}.

Let us now prove Corollary \ref{isoES}. Given $\Sigma_1,\Sigma_2,\Sigma_3$ pairwise transverse exact Lagrangian cobordisms from $\La_i^-$ to $\La_i^+$ as before, Theorem \ref{teofoncteur} gives that the following diagram is commutative.	
	$$\begin{array}{ccc}
	HF_{-\infty}(\Sigma_2,\Sigma_3)\otimes HF_{-\infty}(\Sigma_1,\Sigma_2) & \xrightarrow[]{\Fc_{32}\otimes\Fc_{21}} & H\big(C(\La_2^+,\La_3^+)\big)\otimes H\big(C(\La_1^+,\La_2^+)\big)\\
	\quad&\quad&\quad\\
	\mfm_2\downarrow & & \downarrow\mu^2_{\ep^+_{3,2,1}}\\
	\quad&\quad&\quad\\
	HF_{-\infty}(\Sigma_1,\Sigma_3) & \xrightarrow[\,\;\;\Fc_{31}\,\;\;]{} & H\big(C(\La_1^+,\La_3^+))\\
	\end{array}$$
Now, assume that $\La_1^-=\emptyset$, then $\Sigma_1$ is an exact Lagrangian filling of $\La_1^+$. Then $\ep_1^+$ in this case is given by the DGA map $\phi_\Sigma:\Ac(\La_1^+)\to\Ac(\La_1^-)=\Z_2$. We denote $\Sigma:=\Sigma_1$, $\La:=\La_1^+$ and $\ep_\Sigma:=\ep_1^+$. Let us assume moreover that the cobordisms $\Sigma_2$, $\Sigma_3$ are appropriate Hamiltonian perturbations of $\Sigma$ as considered in Section \ref{perturbations} and that we recall now. First, remember that by the Weinstein Lagrangian neighborhood theorem, a neighborhood of $\Sigma\subset\R\times Y$ is identified with a neighborhood $U_0$ of the $0$-section in $T^*\Sigma$. For $\epsilon$ sufficiently small, then $\widetilde{\Sigma}_2:=\Phi_{H_D}^\epsilon(\Sigma)$ and $\widetilde{\Sigma}_3:=\Phi_{H_D}^{2\epsilon}(\Sigma)$ are identified with the graph of $d(\epsilon H_D)$ and $d(2\epsilon H_D)$ respectively in $U_0$, for the functions $\epsilon H_D$ and $2\epsilon H_D$ restricted to $\Sigma$. Define then $\Sigma_2\subset\R\times Y$ (resp. $\Sigma_3\subset\R\times Y$) to be the exact Lagrangian cobordism identified with the graph of $df_2$ (resp. $df_3$) in $U_0$, for $f_2,f_3:\Sigma\to\R$ Morse functions such that:
\begin{enumerate}
	\item $f_2$ is a small perturbation of $\epsilon H_D$ on $\Sigma$,
	\item $f_3$ is a small perturbation of $2\epsilon H_D$ on $\Sigma$,
	\item the critical points of $f_2$ and $f_3$ are all contained in $\Sigma\cap\big([-T,T]\times Y\big)$,
	\item the cylindrical positive ends $\R\times\La_2^+$ of $\Sigma_2$ and $\R\times\La_3^+$ of $\Sigma_3$ are such that the link $\La\cup\La_2^+\cup\La_3^+$ is a perturbed $3$-copy of $\La$ (see Section \ref{linearisation}).
\end{enumerate}

Then we have $HF_{-\infty}(\Sigma_1,\Sigma_2)=HF(\Sigma_1,\Sigma_2)\simeq H^*(\overline{\Sigma},\La)$, where the first equality comes from the fact that the negative ends are empty, and the second equality comes from Proposition \ref{prop_iso_can} (here we used that $\La=\partial_+\overline{\Sigma})$). We also have $HF_{-\infty}(\Sigma_2,\Sigma_3)\simeq HF_{-\infty}(\Sigma_1,\Sigma_3)\simeq H^*(\overline{\Sigma},\La)$. By Proposition \ref{prop_iso_can} we also have that
\begin{alignat*}{1}
	H\big(C(\La,\La_2^+)\big)\simeq H\big(C(\La_2^+,\La_3^+)\big)\simeq H\big(C(\La,\La_3^+)\big)\simeq LCH^*_{\ep_\Sigma}(\La)
\end{alignat*}
In this case, the vertical arrow on the left in the diagram above is given by a count of pseudo-holomorphic disks with $3$ punctures on the boundary asymptotic to intersection points (the negative ends are empty so there is no Reeb chord asymptotic). Moreover, these intersection points are all contained in a compact so these disks will also all be contained in a compact subset of $\R\times Y$. So by the discussion in \cite[Section 8.l]{S}, the map $\mfm_2$ computes the cup product on $H^*(\Sigma,\La)$. A proof of this fact is also contained in \cite{KMS}, the authors prove an isomorphism between Morse and Floer cohomology rings for the $0$-section of the cotangent bundle of a compact manifold. On the other side, we saw in Section \ref{linearisation} that the product $\mu^2_{\ep_\Sigma}$ on $LCH^*_{\ep_{\Sigma}}(\La)$ can be computed as the product on the $3$-copy of $\La$, which is the vertical map on the right in the diagram above.
We end by a remark on the degree. Recall that the Cthulhu complex is defined with a grading shift, i.e. we have
\begin{alignat*}{1}
	\Cth(\Sigma_1,\Sigma_2)=C(\La_1^+,\La_2^+)[2]\oplus CF(\Sigma_1,\Sigma_2)\oplus C(\La_1^-,\La_2^-)[1]
\end{alignat*}
Considering these shifts, the two horizontal maps in the diagram have degree $1$ and the two vertical maps have degree $0$.

\begin{rem}
	 We could use other types of Hamiltonian perturbations to compute the product. We give here some example with no details. If $\Sigma_1,\Sigma_2,\Sigma_3$ are pairwise transverse exact Lagrangian cobordisms such that $\Sigma_2$ (resp $\Sigma_3$) is a small Morse perturbation of $\Phi_{H_D}^{-\epsilon}(\Sigma_1)$ (resp $\Phi_{H_D}^{-2\epsilon}(\Sigma_1)$), then we have $HF_{-\infty}(\Sigma_i,\Sigma_{i+1})\simeq H^*(\overline{\Sigma}_1,\La_1^-)$. Moreover, by \cite[Proposition 7.5]{CDGG2}, $H^*C(\La_i^\pm,\La_{i+1}^\pm)\simeq LCH_{n-1-*}^{\ep_i^\pm,\ep_{i+1}^\pm}(\La_1^\pm)$. In the same setting but if the cobordisms $\Sigma_i$ are Lagrangian fillings, then Theorem \ref{teofoncteur} illustrates the fact that the Legendrian contact homology is endowed with a unital ring structure, corresponding to the cohomology ring of $\overline{\Sigma}_1$, as observed in \cite[Remark 5.9]{NRSSZ}.
\end{rem}
\begin{rem}\label{remEL}
	Actually, (a unital version of) Corollary \ref{isoES} appears in a paper of Ekholm and Lekili \cite{EL}. Considering an exact Lagrangian filling $\Sigma$ in a Liouville domain $X$ of a Legendrian $\La$ in $\partial X$, the \cite[Theorem 53]{EL} states that there is an $A_\infty$-quasi-isomorphism between the Floer complex of $L$ after adding a unit and the \textit{Legendrian $A_\infty$-algebra} $LA^*_{\ep_\Sigma}(\La)$. This Legendrian $A_\infty$-algebra can be computed using perturbed $k$-copies of $\La$. The case corresponding to Corollary $1$ is when this $k$-copy $\La=\La_1,\La_2,\dots,\La_k$ is such that $\La_i$ is a perturbation of $\La+i\epsilon R$ where $R$ the Reeb vector field, then $LA^*_{\ep_\Sigma}(\La)$ is the $A_\infty$-algebra $\Z_2\oplus LCH_{\ep_\Sigma}^*(\La)$, where the extra $\Z_2$ is used to make the algebra unital. The $A_\infty$-structure on $CF(\Sigma)$ is also computed using a system of parallel copies $\Sigma=\Sigma_1,\Sigma_2,\dots,\Sigma_{d+1}$ of the filling. The maps giving the $A_\infty$-structure are then the same as the maps $\mfm_k$ (defined in Section \ref{infini} below) for the case of fillings.
	Moreover, the family of maps we construct in Theorem \ref{ffon} in this paper, for the case of fillings, recover the maps defining the functor in the proof of \cite[Theorem 53]{EL} computed with the parallel copies such that the positive end consists in the $k$-copy $\La=\La_1,\La_2,\dots,\La_k$ as above.
\end{rem}

\section{Example}\label{example}
\begin{figure}[ht]  
	\begin{center}\includegraphics[width=9cm]{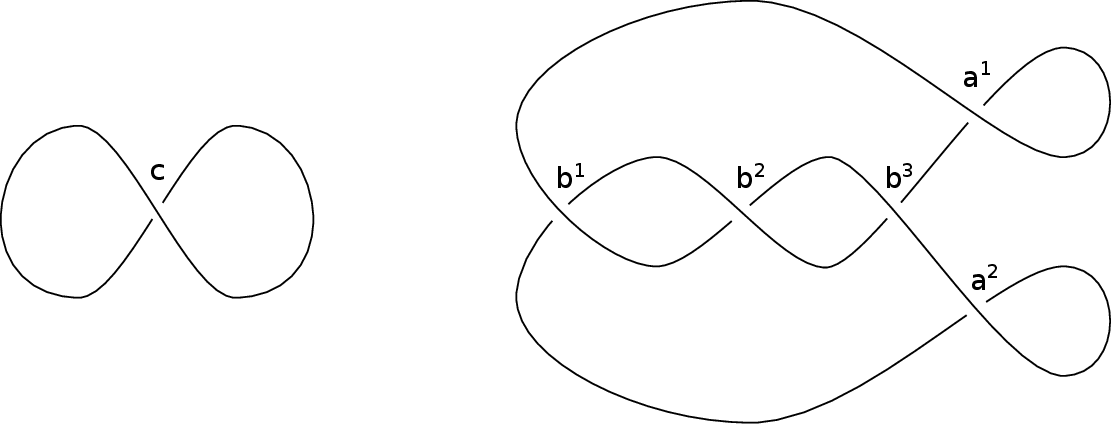}\end{center}
	\caption{Left: Lagrangian projection of $\La_0^-$, right: Lagrangian projection of $\La_0^+$.}
	\label{untr}
\end{figure}

We give in this section a very simple example of computation of the product using Theorem \ref{teofoncteur}. Consider a cobordism $\La_1^-\prec_{\Sigma_1}\La_1^+$ where $\La_1^-$ is the Legendrian unknot and $\La_1^+$ the right-handed trefoil in $\R^3$ (see Figure \ref{untr}).
The unknot has one Reeb chord $c$ in degree $1$, and the Chekanov-Eliashberg DGA $\Ac(\La_1^-)$ admits a unique augmentation $\ep^-$ which is trivial ($\ep^-(c)=0$). The trefoil has $5$ Reeb chords with $|b^i|=0$ and $|a^i|=1$. The differential is given by 
\begin{alignat*}{1}
	&\partial a^1=1+b^1+b^3+b^1b^2b^3\\
	&\partial a^2=1+b^1+b^3+b^3b^2b^1\\
	&\partial b^i=0,\; i=1,2,3
\end{alignat*}
The DGA $\Ac(\La_1^+)$ admits five augmentations, which are all geometric (see \cite{EHK}). Let us assume that the cobordism $\Sigma_1$ is such that $\ep^-\circ\phi_{\Sigma_1}=\ep^+$ where $\ep^+$ is the augmentation of $\Ac(\La_1^+)$ defined by $\ep^+(b^1)=1$ and $\ep^+(b^2)=\ep^+(b^3)=0$.

Consider now a small Morse perturbation $\Sigma_2$ of $\Phi_{H_D}^\epsilon(\Sigma_1)$ by a Morse function $f_2$ as described in Section \ref{perturbations} and in the proof of Corollary \ref{isoES}. We have $\La_2^-\prec_{\Sigma_2}\La_2^+$, where $\La_2^\pm$ is a perturbation of $\La_1^\pm+\epsilon\frac{\partial}{\partial z}$. Moreover, the critical points of $f_2$ are contained in $\Sigma_1\cap([-T,T]\times \R^3)$ and we assume that there is no minimum (it is possible since $\Sigma_1$ is a punctured torus). For a small enough perturbation, the DGAs $\Ac(\La_1^\pm)$ and $\Ac(\La_2^\pm)$ are the same (by canonical identification of the Reeb chords). Furthermore, the Reeb chords from $\La_2^\pm$ to $\La_1^\pm$ are in bijection with the Reeb chords of $\La_1^\pm$ (Section \ref{linearisation}), we denote by $\gamma_{21}$ the chord from $\La_2$ to $\La_1$ corresponding to the chord $\gamma$ of $\La_1$.

We have $C(\La_1^+,\La_2^+)=\Z_2\langle a_{21}^1,a_{21}^2,b^1_{21},b^2_{21},b^3_{21}\rangle$ and using the augmentation $\ep^+$ we can compute 
\begin{alignat*}{1}
	&d_{++}(b^1_{21})=d_{++}(b^3_{21})=a_{21}^1+a_{21}^2\\
	&d_{++}(b^2_{21})=0\\
	&d_{++}(a_{21}^i)=0,\; i=1,2
\end{alignat*}
and we get $H\big((C(\La_1^+,\La_2^+),d_{++})\big)\simeq\Z_2\langle[a_{21}^1]\rangle\oplus\Z_2\langle [b^1_{21}+b^3_{21}],[b^2_{21}]\rangle\simeq LCH^*_{\ep^+}(\La_1^+)$.

On the other side we have $C(\La_1^-,\La_2^-)=\Z_2\langle c_{21}\rangle$, and $d_{--}=0$, thus $H\big((C(\La_1^-,\La_2^-),d_{--})\big)\simeq\Z_2\langle[c_{21}]\rangle\simeq LCH^*_{\ep^-}(\La_1^-)$.

Using a third copy $\Sigma_3$ being a perturbation of $\Phi^{2\epsilon}_{H_D}(\Sigma_1)$ by a Morse function $f_3$, such that $f_3-f_2$ is also Morse and these two Morse functions have no minima, we compute that the non trivial components of
\begin{alignat*}{1}
	\mu_{\ep^+}^2:C(\La_2^+,\La_3^+)\otimes C(\La_1^+,\La_2^+)\to C(\La_1^+,\La_3^+)
\end{alignat*}
in homology are $\mu_{\ep^+}^2([b_{32}^1+b_{32}^3],[b_{21}^2])=\mu_{\ep^+}^2([b_{32}^2],[b_{21}^1+b_{21}^3])=[a_{31}^1]$ (which under canonical identification of the generators is the product on $LCH_{\ep^+}^*(\La_1^+)$).

We now want to compute $HF_{-\infty}(\Sigma_1,\Sigma_2)$ and the product structure on it. Recall that $\Fc^1:CF_{-\infty}(\Sigma_1,\Sigma_2)\to C(\La_1^+,\La_2^+)$ is a quasi-isomorphism, so $HF_{-\infty}(\Sigma_1,\Sigma_2)$ is of rank $3$. Also, there exists $[A]=[A_0+A_-]\in HF_{-\infty}(\Sigma_1,\Sigma_2)$, with $A_0\in CF(\Sigma_1,\Sigma_2)$ and $A_-\in C(\La_1^-,\La_2^-)$, such that $\Fc^1([A])=[a_{21}^1]$ in homology. This gives
\begin{alignat*}{1}
	[d_{+0}(A_0)+d_{+-}(A_-)]=[a^1_{21}]
\end{alignat*}
which implies $|A_0|=|a_{21}^1|+1=2$ and $|A_-|=|a_{21}^1|=1$ as $d_{+0}$ is a degree $-1$ map and $d_{+-}$ a degree $0$ map. By Proposition \ref{prop_iso_can}, $A_0$ corresponds to a linear combination of critical points of $f_2$ of Morse index $0$, but we have assumed that $f_2$ has no minimum so there is no such $A_0$. Then the only possibility is $A_-=c_{21}$ and so $c_{21}$ is a non-trivial cycle in $HF_{-\infty}(\Sigma_1,\Sigma_2)$, with $\Fc([c_{21}])=[a_{21}^1]$. Then, there must also exist $[x_{21}],[y_{21}]\in HF_{-\infty}(\Sigma_1,\Sigma_2)$, with $x_{21},y_{21}\in CF(\Sigma_1,\Sigma_2)$ of degree $1$ such that $\Fc([x_{21}])=[b_{21}^1+b_{21}^3]$ and $\Fc([y_{21}])=[b_{21}^2]$. Therefore, by Theorem \ref{teofoncteur}, the product in homology
\begin{alignat*}{1}
	\mfm_2: HF_{-\infty}(\Sigma_2,\Sigma_3)\otimes HF_{-\infty}(\Sigma_1,\Sigma_2)\to HF_{-\infty}(\Sigma_1,\Sigma_3)
\end{alignat*}
is given by $\mfm_2([x_{32}],[y_{21}])=\mfm_2([y_{32}],[x_{21}])=[c_{31}]$.

Of course, if we fill the unknot by a disk then $\Sigma_i$ are Lagrangian fillings of the trefoil. In this case, the homology $HF_{-\infty}(\Sigma_1,\Sigma_2)$ is generated by $x_{21}$, $y_{21}$ and a minimum $m_{21}\in CF(\Sigma_1,\Sigma_2)$. By Corollary \ref{isoES}, the non trivial components of the product on $HF_{-\infty}$ are the non trivial components of the cup product on the punctured torus
\begin{alignat*}{1}
\cup:H^1(\Sigma_1,\La_1^+)\otimes H^1(\Sigma_1,\La_1^+)\to H^2(\Sigma_1,\La_1^+)
\end{alignat*}

\section{An $A_\infty$-structure}\label{infini}

The goal of this section is to show that the product structure can be expanded to an $A_\infty$-structure induced by a fixed sequence of pairwise transverse cobordisms.
More precisely, let us consider a fixed $(d+1)$-tuple of transverse cobordisms $\Sigma_1,\dots,\Sigma_{d+1}$. For every $1\leq k\leq d$ and every $(k+1)$-tuple $i_1,\dots,i_{k+1}$ of distinct indices in $\{1,\dots,d+1\}$, we will construct a map $\mfm_k$: 
\begin{alignat*}{1}
\mfm_k\colon CF_{-\infty}(\Sigma_{i_k},\Sigma_{i_{k+1}})\otimes\dots\otimes CF_{-\infty}(\Sigma_{i_1},\Sigma_{i_2})\to CF_{-\infty}(\Sigma_{i_1},\Sigma_{i_{k+1}})
\end{alignat*}
such that the family of maps $\{\mfm_k\}_{1\leq k\leq d}$ satisfies for every $1\leq k\leq d$ and every $(k+1)$-tuple of distinct indices $i_1,\dots,i_{k+1}$:
\begin{alignat*}{1}
\sum\limits_{\substack{1\leq j\leq k\\ 0\leq n\leq k-j}}\mfm_{k-j+1}(\id^{\otimes k-j-n}\otimes\mfm_j\otimes\id^{\otimes n})=0
\end{alignat*}
where in the sum above, for $1\leq j\leq k$ and $0\leq n\leq k-j$, we have 
\begin{alignat*}{1}
	\mfm_j: CF_{-\infty}(\Sigma_{i_{n+j}},\Sigma_{i_{n+j+1}})\otimes\dots\otimes CF_{-\infty}(\Sigma_{i_{n+1}},\Sigma_{i_{n+2}})\to CF_{-\infty}(\Sigma_{i_{n+1}},\Sigma_{i_{n+j+1}})
\end{alignat*}
and $\mfm_{k-j+1}$ has domain
\begin{alignat*}{1}
CF_{-\infty}(\Sigma_{i_{k}},\Sigma_{i_{k+1}})\otimes\dots\otimes CF_{-\infty}(\Sigma_{i_{n+1}},\Sigma_{i_{n+j+1}})\otimes\dots\otimes CF_{-\infty}(\Sigma_{i_{1}},\Sigma_{i_{2}})
\end{alignat*}
and codomain $ CF_{-\infty}(\Sigma_{i_{1}},\Sigma_{i_{k+1}})$.

In order to simplify notations when defining these maps in the following section, we will assume without loss of generality that the $(k+1)$-tuple of distinct indices $i_1,\dots,i_{k+1}$ is equal to $1,\dots,k+1$. Then, for each $1\leq k\leq d$, we have $\mfm_k=m^0_k+m^-_k$ where $m^0_k$ takes values in $CF(\Sigma_1,\Sigma_{k+1})$ and $m^-_k$ takes values in $C^*(\La_1^-,\La_{k+1}^-)$. We will define those two components separately.

\subsection{Definition of the operations}\label{definf}

Let $\Sigma_1,\dots,\Sigma_{d+1}$, for $d\geq2$, be transverse Lagrangian cobordisms from $\La_i^-$ to $\La_i^+$ for $i=1,\dots, d+1$ such that the algebras $\Ac(\La_i^-)$ admit augmentations $\ep_i^-$. For $1\leq k\leq d$, we define $\mfm_k$ as follows. First, for $k=1$, the map $\mfm_1:CF_{-\infty}(\Sigma_1,\Sigma_2)\to CF_{-\infty}(\Sigma_1,\Sigma_2)$ is the differential $\partial_{-\infty}$ on the Floer complex. For $k=2$, it is the product on Floer complexes as defined in Section \ref{defprod}. Then, for $3\leq k\leq d$, the component $m_k^0$ is naturally the generalization of $m_2^0$ and is thus defined by a count of rigid pseudo-holomorphic disks with boundary on non cylindrical parts of the cobordisms, and with $k+1$ mixed asymptotics. Indeed, we define:
\begin{alignat*}{1}
m_k^0\colon CF^*_{-\infty}(\Sigma_k,\Sigma_{k+1})\otimes\dots\otimes CF^*_{-\infty}(\Sigma_1,\Sigma_2)\to CF^*(\Sigma_1,\Sigma_{k+1})
\end{alignat*}
by
\begin{alignat}{1}
m_k^0(a_k,\dots,a_1)&=\sum\limits_{\substack{x^+\in\Sigma_1\cap\Sigma_{k+1}\\ \bs{\delta}_1,\dots,\bs{\delta}_{k+1}}}\#\cM^0(x^+;\bs{\delta}_1,a_1,\bs{\delta}_2,\dots,a_k,\bs{\delta}_{k+1})\cdot\ep^-\cdot x^+ \label{prod0}
\end{alignat}
where $\bs{\delta}_i$ are words of Reeb chords of $\La_i^-$ for $1\leq i\leq k+1$, and the term ``$\ep^-$'' means that we augment all the pure Reeb chords with the corresponding augmentations, i.e. $\ep^-$ should be replaced by $\ep_1^-(\bs{\delta}_1)\ep_2^-(\bs{\delta}_2)\dots\ep_{k+1}^-(\bs{\delta}_{k+1})$ in the formula. Also, the choice of Lagrangian label for the moduli spaces involved in the definition of $m_k^0$ is $(\Sigma_1,\dots,\Sigma_{k+1})$. Now let us define $m_k^-$. As in the case $k=2$, this map is defined by a count of unfinished pseudo-holomorphic buildings, except when all the asymptotics are Reeb chords, and so we define it as a composition of maps.
First, consider the map 
\begin{alignat*}{1}
f^{(k)}\colon CF_{-\infty}(\Sigma_k,\Sigma_{k+1})\otimes\dots\otimes CF_{-\infty}(\Sigma_1,\Sigma_2)\to C_{n-1-*}(\La_{k+1}^-,\La_1^-)
\end{alignat*}
defined for a $k$-tuple of asymptotics $(a_k,\dots,a_1)$ with $a_i\in CF_{-\infty}(\Sigma_i,\Sigma_{i+1})$ by:
\begin{alignat*}{1}
f^{(k)}(a_k,\dots,a_1)=\sum\limits_{\substack{\gamma_{1,k+1}\\ \bs{\delta}_1,\dots,\bs{\delta}_{k+1}}}\#\cM^0_{\Sigma_{1,\dots,k+1}}(\gamma_{1,k+1};\bs{\delta}_1,a_1,\bs{\delta}_2,\dots,a_k,\bs{\delta}_{k+1})\cdot\ep^-\cdot\gamma_{1,k+1}
\end{alignat*}
These maps $f^{(k)}$ are generalizations of the maps $f^{(1)}$ and $f^{(2)}$ defined in Section \ref{defprod}. For $k\geq2$, in the case where all the mixed asymptotics are Reeb chords $(\gamma_k,\dots,\gamma_1)$, with $\gamma_i\in C^*(\La_i^-,\La_{i+1}^-)$, we have $f^{(k)}(\gamma_k,\dots,\gamma_1)=0$ for energy reasons. However, recall that for $k=1$ and $\gamma\in C^*(\La_1^-,\La_{2}^-)$ a Reeb chord, we set by convention $f^{(1)}(\gamma)=\gamma$ (and not $f^{(1)}(\gamma)=0$).

\begin{figure}[ht]  
	\begin{center}\includegraphics[width=8cm]{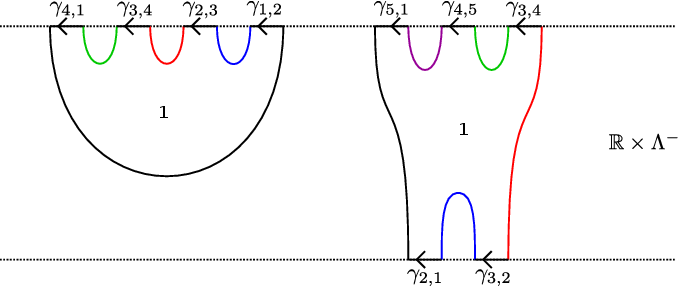}\end{center}
	\caption{Left: a curve contributing to $b^{(3)}(\gamma_{3,4},\gamma_{2,3},\gamma_{1,2})$; right: a curve contributing to $b^{(4)}(\gamma_{4,5},\gamma_{3,4},\gamma_{3,2},\gamma_{2,1})$.}
	\label{banane_general}
\end{figure}

\noindent Now we generalize the bananas $b$ (Section \ref{cthulhu}) and $b^{(2)}$ (Section \ref{defprod}) in a family of maps $b^{(k)}$, $1\leq k\leq d$. For $j>i$, recall that we denote $\mathfrak{C}^*(\La_i^-,\La_j^-)=C_{n-1-*}(\La_j^-,\La_i^-)\oplus C^*(\La_i^-,\La_j^-)$. We define for all $1\leq k\leq d$:
\begin{alignat*}{1}
b^{(k)}\colon\mathfrak{C}^*(\La_k^-,\La_{k+1}^-)\otimes\dots\otimes\mathfrak{C}^*(\La_1^-,\La_2^-)\to C^*(\La_1^-,\La_{k+1}^-)
\end{alignat*}
by 
\begin{alignat*}{1}
b^{(k)}(\gamma_k,\dots,\gamma_1)=\sum\limits_{\substack{\gamma_{k+1,1}\\ \bs{\delta}_1,\dots,\bs{\delta}_{k+1}}}\#\widetilde{\cM^1}(\gamma_{k+1,1};\bs{\delta}_1,\gamma_1,\bs{\delta}_2,\dots,\gamma_k,\bs{\delta}_{k+1})\cdot\ep^-\cdot\gamma_{k+1,1}
\end{alignat*}
where the choice of Lagrangian label is $(\R\times\La_1^-,\dots,\R\times\La_{k+1}^-)$, and the $\bs{\delta}_i$ are still words of Reeb chords of $\La_i^-$ and are negative asymptotics. This formula for $k=1$ gives $b^{(1)}=b+d_{--}$ and for $k=2$ it is the same formula as in Section \ref{defprod} to define $b^{(2)}$.
Remark also that as for $k=2$, for any $k$-tuple $(\gamma_k,\dots,\gamma_1)$ of chords $\gamma_i\in C^*(\La_i^-,\La_{i+1}^-)$ with $k\geq2$, the map $b^{(k)}$ is equal to the map $\mu_{\ep_{k+1}^-,\dots,\ep_1^-}^k$ restricted to $C(\La_k^-,\La_{k+1}^-)\otimes\dots\otimes C(\La_1^-,\La_2^-)$ in the augmentation category $\Aug_-(\La_1^-\cup\dots\cup\La_{k+1}^-)$. On Figure \ref{banane_general} are schematized examples of curves involved in the definition of the banana maps.

\begin{figure}[ht]  
	\begin{center}\includegraphics[width=8cm]{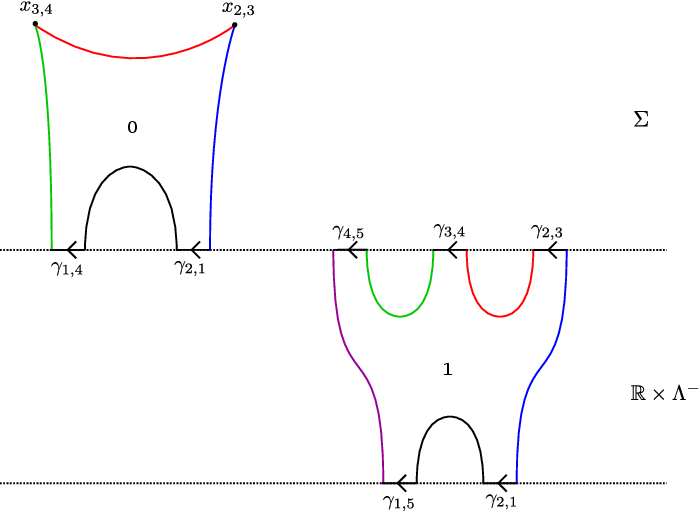}\end{center}
	\caption{Left: a curve contributing to $f^{(3)}(x_{3,4},x_{2,3},\gamma_{2,1})$; right: a curve contributing to $\Delta^{(4)}(\gamma_{4,5},\gamma_{3,4},\gamma_{2,3},\gamma_{2,1})$.}
	\label{exemple_f_3_delta_4}
\end{figure}

\noindent Finally, for $3\leq k\leq d$, we generalize the maps $\Delta^{(1)}:=\delta_{--}$ and $\Delta^{(2)}$ (defined in Section \ref{Rel-xx}) by:
\begin{alignat*}{1}
\Delta^{(k)}\colon \mathfrak{C}^*(\La_k^-,\La_{k+1}^-)\otimes\dots\otimes \mathfrak{C}^*(\La_1^-,\La_2^-)\to C_{n-1-*}(\La_{k+1}^-,\La_1^-)
\end{alignat*}
defined for $(\gamma_k,\dots,\gamma_1)$ a $k$-tuple of Reeb chords, with $\gamma_i\in\mathfrak{C}^*(\La_i^-,\La_{i+1}^-)$, by:
\begin{alignat*}{1}
\Delta^{(k)}(\gamma_k,\dots,\gamma_1)=\sum\limits_{\substack{\gamma_{1,k+1}\\ \bs{\delta}_i}}\#\widetilde{\cM^1}(\gamma_{1,k+1};\bs{\delta}_1,\gamma_1,\bs{\delta}_2,\dots,\gamma_k,\bs{\delta}_{k+1})\cdot\ep^-\cdot\gamma_{1,k+1}
\end{alignat*}
where the Lagrangian label is $(\R\times\La_1^-,\dots,\R\times\La_{k+1}^-)$. If $(\gamma_k,\dots,\gamma_1)$ is a $k$-tuple of Reeb chords $\gamma_i\in C^*(\La_i^-,\La_{i+1}^-)$, we have $\Delta^{(k)}(\gamma_k,\dots,\gamma_1)=0$ for energy reasons. Remark that chords from $\La_i^-$ to $\La_{i+1}^-$ that are asymptotics of curves in moduli spaces involved in the definition of $\Delta^{(k)}$ are positive asymptotics, while chords from $\La_{i+1}^-$ to $\La_i^-$ are negative asymptotics. These maps $\Delta^{(k)}$ are not directly involved in the definition of $m_k^-$ but they will be useful in order to express algebraically the unfinished pseudo-holomorphic buildings appearing in the study of breakings of pseudo-holomorphic curves. We schematized on Figure \ref{exemple_f_3_delta_4} examples of curves contributing to $f^{(k)}$ and $\Delta^{(k)}$.
Now we can finally define the map:
\begin{alignat*}{1}
m_k^-&\colon CF^*_{-\infty}(\Sigma_k,\Sigma_{k+1})\otimes\dots\otimes CF^*_{-\infty}(\Sigma_1,\Sigma_2)\to C^*(\La_1^-,\La_{k+1}^-)
\end{alignat*}
by setting
\begin{alignat}{1}
m_k^-(a_k,\dots,a_1)=\sum\limits_{\substack{1\leq j\leq k\\ i_1+\dots+i_j=k}}b^{(j)}\big(f^{(i_j)}(a_k,...,a_{k-i_j+1}),\dots,f^{(i_1)}(a_{i_1},...,a_1)\big)\label{prod-}
\end{alignat}
for a $k$-tuple of generators $(a_k,\dots,a_1)$, and recall the following conventions on the $f^{(i)}$'s in the formula:
\begin{alignat}{1}
&f^{(1)}(a_i)=a_i \,\,\mbox{ if }\,\, a_i=\gamma_{i+1,i}\label{conv1}\\
&f^{(s)}(\gamma_{i+1,i},\gamma_{i,i-1},\dots,\gamma_{i-s+2,i-s+1})=0\,\, \mbox{ for }\,\, 1<s\leq i\leq d\label{conv2}
\end{alignat}
Remark that the formulas \eqref{prod0} and \eqref{prod-} in the case $k=1$ give $\mfm_1=\partial_{-\infty}$ and in the case $k=2$ recover the product $\mfm_2$ as defined in Section \ref{defprod}. 

\subsection{Proof of Theorem \ref{teo3}}\label{proofteo3}
In order to show the $A_\infty$-relations, again we study breakings of pseudo-holomorphic curves. The $A_\infty$-relations for the maps $\{\mfm_k\}_{1\leq k\leq d}$ can be rewritten as follows. For all $1\leq k\leq d$:
\begin{alignat*}{1}
\sum\limits_{\substack{1\leq j\leq k\\ 0\leq n\leq k-j}}m^0_{k-j+1}(\id^{\otimes k-j-n}\otimes \mfm_j&\otimes\id^{\otimes n})\\
&+\sum\limits_{\substack{1\leq j\leq k\\ 0\leq n\leq k-j}}m^-_{k-j+1}(\id^{\otimes k-j-n}\otimes \mfm_j\otimes\id^{\otimes n})=0
\end{alignat*}
First we start by showing that
\begin{alignat}{1}
\sum\limits_{\substack{1\leq j\leq k\\ 0\leq n\leq k-j}}m^0_{k-j+1}(\id^{\otimes k-j-n}\otimes \mfm_j&\otimes\id^{\otimes n})=0\label{premier terme}
\end{alignat}
and then we will prove that
\begin{alignat}{1}
\sum\limits_{\substack{1\leq j\leq k\\ 0\leq n\leq k-j}}m^-_{k-j+1}(\id^{\otimes k-j-n}\otimes \mfm_j\otimes\id^{\otimes n})=0\label{second terme}
\end{alignat}

\subsubsection{Proof of Relation \eqref{premier terme}.}
To show this relation we need to understand the different types of pseudo-holomorphic buildings contributing to the maps in the sum. For a $k$-tuple $(a_k,\dots,a_1)$ of asymptotics, each term of \eqref{premier terme} is either of the form
\begin{alignat}{1}
m^0_{k-j+1}\big(a_k,\dots,m_j^0(a_{n+j},\dots,a_{n+1}),a_n,\dots,a_1\big)\label{premier terme1}
\end{alignat}
or
\begin{alignat}{1}
m^0_{k-j+1}\big(a_k,\dots,m_j^-(a_{n+j},\dots,a_{n+1}),a_n,\dots,a_1\big)\label{premier terme2}
\end{alignat}
The pseudo-holomorphic buildings contributing to \eqref{premier terme1} are of height $0|1|0$. The middle level of each building contains two curves which have a common asymptotic on an intersection point, and can be glued on a pseudo-holomorphic disk in the moduli space
\begin{alignat}{1}
\cM^1(x^+;\bs{\delta}_1,a_1,\bs{\delta}_2,a_2,\dots,\bs{\delta}_{k},a_k,\bs{\delta}_{k+1})\label{espace1}
\end{alignat}
The same happens for pseudo-holomorphic buildings contributing to the terms in \eqref{premier terme2}. Such a building is of height $1|1|0$ and has components that can be glued on an index-1 curve in the moduli space \eqref{espace1}. Indeed, $m^-$ counts unfinished buildings (as soon as one asymptotic at least is an intersection point, otherwise it counts just one banana) of height $1|1|0$ such that the curve with boundary on the negative cylindrical ends is a banana which has for output a positive chord $\gamma_{n+j+1,n+1}\in\Rc(\La^-_{n+j+1},\La^-_{n+1})$. The map $m^0$ applied to the remaining asymptotics and $\gamma_{n+j+1,n+1}$ is then given by the count of index-0 pseudo-holomorphic curves in the middle level. The unfinished buildings contributing to $m^-$ and the curves contributing to $m^0$ together give a pseudo-holomorphic building, and the corresponding glued curve is an index-1 pseudo-holomorphic curve in \eqref{espace1}.

Now, in order to establish Relation \eqref{premier terme}, we must study the boundary of the compactification of this moduli space.
As for the case $d=2$, the pseudo-holomorphic buildings arising as limit of a one parameter family of disks in \eqref{espace1} must satisfy some conditions that we recall here:
\begin{enumeratec}
	\item each curve in the building must have positive energy,
	\item each curve in the building has a non negative Fredholm index,
	\item the building is asymptotic (the asymptotics that are not nodes) to $x^+,a_1,\dots,a_k$, so in particular contains at least one non trivial curve with boundary on the compact parts of the cobordisms, having the intersection point $x^+$ as asymptotic.
	\item if the building consists of the pseudo-holomorphic disks $\{u_i\}$, the relation $\sum \ind(u_i)+\nu=1$ must be satisfied (where $\nu$ is the number of pair of nodes asymptotic to intersection points, see Section \ref{comp}), which implies that there are basically two types of buildings to consider:
	\begin{enumeratec}
		\item buildings with two rigid components and boundary on the compact parts having a common asymptotic at an intersection point (pair of nodes),
		\item buildings with several disjoint rigid components with boundary on the compact parts, each having a node asymptotic to a negative Reeb chord, and one index-1 component with boundary on the negative ends, having these Reeb chords as positive asymptotics, and possibly some negative Reeb chord asymptotics among $a_1,\dots,a_k$.
	\end{enumeratec}	
\end{enumeratec}

In the case (4)(a), the limit building is of height $0|1|0$ and the middle level contains two index-$0$ disks having a node asymptotic to an intersection point $q\in CF(\Sigma_{n+1},\Sigma_{n+j+1})$ for some $1\leq j\leq k$ and $0\leq n\leq k-j$. One disk is asymptotic to $(x^+,a_1,\dots,a_{n},q,a_{n+j+1},\dots,a_k)$, the other one is asymptotic to $(q,a_{n+1},\dots,a_{n+j})$ in this cyclic order when following the boundary of the curves counter-clockwise. Such a pseudo-holomorphic building contributes then to:
\begin{alignat*}{1}
m^0_{k-j+1}\big(a_k,\dots,m^0_j(a_{n+j},\dots,a_{n+1}),a_n,\dots,a_1\big)
\end{alignat*}

\begin{figure}[ht]  
	\begin{center}\includegraphics[width=9cm]{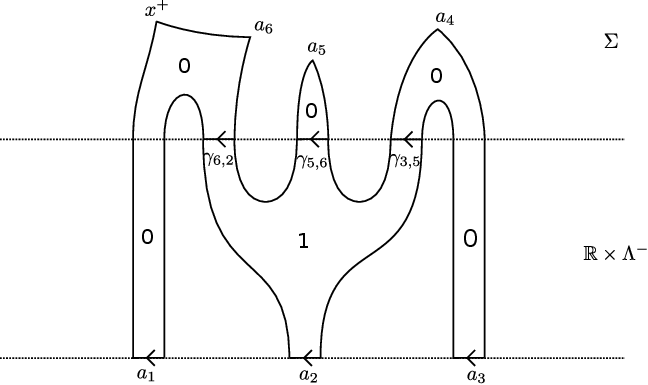}\end{center}
	\caption{Example of pseudo-holomorphic building arising as limit of a family of disks in $\cM^1(x^+;a_1,a_2,a_3,a_4,a_5,a_6)$.}
	\label{degef}
\end{figure}

In the case (4)(b), the limit building is of height $1|1|0$. The middle level contains some rigid curves and the bottom level contains one index-1 banana. For some $1\leq j\leq k$ and $0\leq n\leq k-j$, there is a disk in the middle level asymptotic to $x^+,a_1,\dots,a_{n},\gamma_{n+j+1,n+1}$, $a_{n+j+1},\dots,a_k$ in this cyclic order and which contributes to the map $m^0$. Then, there are some Reeb chords $\gamma_{\alpha_s,\alpha_{s+1}}$ for $1\leq s\leq r$, such that $(\alpha_s)_s$ is a strictly increasing finite sequence of length $r\leq j$, with $n+1\leq\alpha_s\leq n+j+1$ such that each of the other disks in the middle level has one negative puncture asymptotic to a chord $\gamma_{\alpha_s,\alpha_{s+1}}$ and $\alpha_{s+1}-\alpha_s$ other asymptotics. Such a curve contributes to the map $f^{(\alpha_{s+1}-\alpha_s)}$. The banana in the bottom level has positive Reeb chords asymptotics to $\gamma_{n+j+1,n+1}$ (the output), and $\gamma_{\alpha_s,\alpha_{s+1}}$ (which are inputs), and possibly negative Reeb chords among the asymptotics $(a_{n+j},\dots,a_{n+1})$ which are not asymptotics of curves in the middle level (see Figure \ref{degef} for an example of such kind of breaking). So finally such a pseudo-holomorphic building contributes to:
\begin{alignat*}{1}
m^0_{k-j+1}\Big(a_k,\dots,\sum_{\substack{1\leq s\leq j\\ i_1+\dots+i_s=j}}b^{(s)}\big(f^{(i_s)}(a_{n+j},\dots),\dots,f^{(i_1)}(\dots)\big),a_n,\dots,a_1\Big)
\end{alignat*}
which is by definition equal to
\begin{alignat*}{1}
m^0_{k-j+1}\big(a_k,\dots,m^-_j(a_{n+j},\dots,a_{n+1}),a_n,\dots,a_1\big)
\end{alignat*}

We have thus described every types of pseudo-holomorphic buildings arising as limits of curves in the moduli space $\cM^1(x^+;\bs{\delta}_1,a_1,\bs{\delta}_2,a_2,\dots,\bs{\delta}_{k},a_k,\bs{\delta}_{k+1})$. These buildings are in bijection with the elements in the boundary of the compactification of the moduli space. This compactification being a $1$-dimensional manifold with boundary, its boundary components arise in pair, which gives 0 modulo $2$. This implies Relation \eqref{premier terme}.
\begin{rem}
	As in Section \ref{PROD}, $\partial$-breaking can occur when a family of curves in the moduli space $\cM^1(x^+;\bs{\delta}_1,a_1,\bs{\delta}_2,a_2,\dots,\bs{\delta}_{k},a_k,\bs{\delta}_{k+1})$ breaks on a pure Reeb chord $\gamma\in\Rc(\La_i^-)$. In such a type of breaking, we get a pseudo-holomorphic building of height $1|1|0$ with one rigid component in the middle level and a disk contributing to $\partial^i(\gamma)$ in the bottom level, where $\partial^i$ is the differential on the Chekanov-Eliashberg algebra asssociated to $\La_i^-$ (Section \ref{LCHdef}). Then as we apply the augmentations $\ep_j^-$ to all pure negative Reeb chords, and as $\ep_i^-\circ\partial^i=0$, the contribution of such a pseudo-holomorphic building vanishes. In the following section, $\partial$-breaking also occurs, but its contribution vanishes for the same reason so we don't mention it. 
\end{rem}

\subsubsection{Proof of Relation \eqref{second terme}.}
\begin{figure}[ht]  
	\begin{center}\includegraphics[width=13cm]{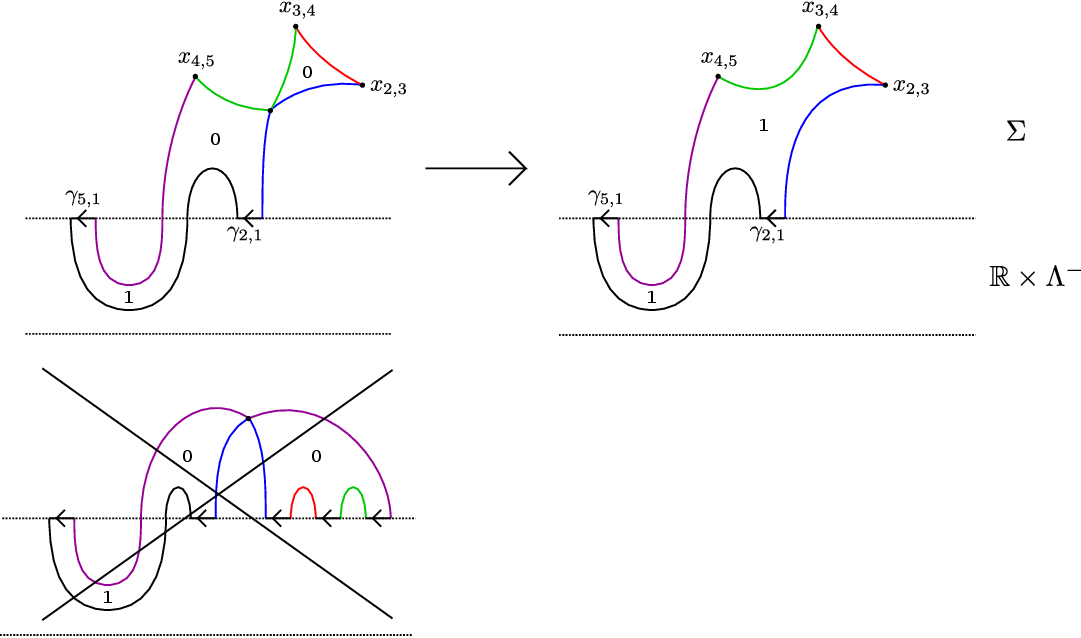}\end{center}
	\caption{On the top: example of unfinished building contributing to $b^{(1)}\circ f^{(4)}(x_{4,5},x_{3,4},x_{2,3},\gamma_{2,1})$ and the corresponding glued curve. On the bottom: impossible breaking for energy reasons.}
	\label{brisure_impossible}
\end{figure}

In this section, in order to cause less confusion, we do not write the pure chords asymptotics in the moduli spaces anymore. As before, the left-hand side of Relation \eqref{second terme}, with inputs a $k$-tuple of asymptotics $(a_k,\dots,a_1)$, splits into two sums:
\begin{alignat}{1}
\sum\limits_{\substack{1\leq j\leq k\\ 0\leq n\leq k-j}}m^-_{k-j+1}\big(a_k,\dots,a_{n+j+1},m^0_j(a_{n+j},\dots,a_{n+1}),a_n,\dots,a_1\big)\label{second terme1}
\end{alignat}
and
\begin{alignat}{1}
\sum\limits_{\substack{1\leq j\leq k\\ 0\leq n\leq k-j}}m^-_{k-j+1}\big(a_k,\dots,a_{n+j+1},m^-_j(a_{n+j},\dots,a_{n+1}),a_n,\dots,a_1\big)\label{second terme2}
\end{alignat}
First, let us look for example the term $$b^{(1)}\circ f^{(k-j+1)}\big(a_k,\dots,m^0_j(a_{n+j},\dots,a_{n+1}),a_n,\dots a_1\big)$$ appearing in \eqref{second terme1}. If $(a_k,\dots,a_1)$ is a $k$-tuple of Reeb chords, then as we already saw this term vanishes for energy reasons. So let us assume that at least one $a_i$ is an intersection point. Then, the term we consider is given by a count of unfinished buildings of height $1|1|0$ with two components in the middle level having a common asymptotic at an intersection point, and a banana in the bottom level. Gluing the middle level components gives an unfinished building in (see Figure \ref{brisure_impossible}):
\begin{alignat*}{1}
\widetilde{\cM^1}(\gamma_{k+1,1};\gamma_{1,k+1})\times\cM^1(\gamma_{1,k+1};a_1,a_2,\dots,a_k)
\end{alignat*}

\begin{figure}[ht]  
	\begin{center}\includegraphics[width=13cm]{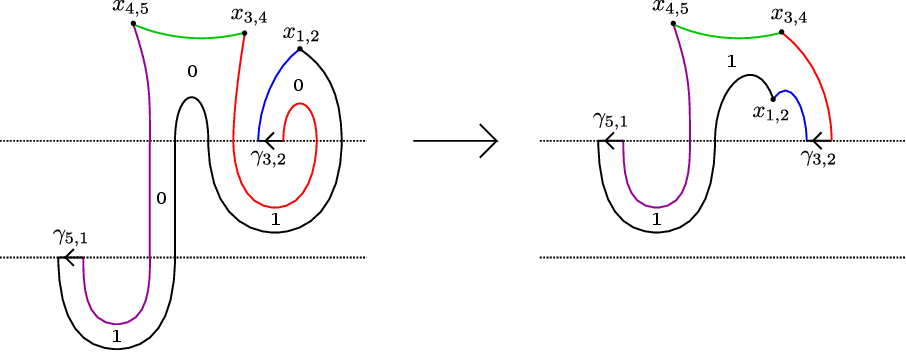}\end{center}
	\caption{Unfinished pseudo-holomorphic building contributing to $b^{(1)}\circ f^{(3)}(x_{4,5},x_{3,4},b^{(1)}\circ f^{(2)}(\gamma_{3,2},x_{1,2}))$ and the corresponding glued curve.}
	\label{recoll1}
\end{figure}

\noindent Let us take another term of \eqref{second terme1}, for example:
\begin{alignat*}{1}
b^{(2)}\Big(f^{(k-n-j)}\big(a_k,\dots,a_{n+j+1}\big),f^{(n+1)}\big(m^0_j(a_{n+j},\dots,a_{n+1}),a_n,\dots a_1\big)\Big)
\end{alignat*}
This one counts again unfinished buildings of height $1|1|0$ with the following conditions:
\begin{enumerate}
	\item assume $a_k,\dots,a_{n+j+1}$ are Reeb chords: if $n+j+1<k$, then $f^{(k-n-j)}(a_k,\dots,a_{n+j+1})=0$, and if $n+j+1=k$, then $f^{(1)}(a_k)=a_k$. In this latter case, the term above counts unfinished buildings with two components in the middle level having a common asymptotic to an intersection point, and one banana in the bottom level having 3 mixed Reeb chords asymptotics: one output positive Reeb chord asymptotic at a chord $\gamma_{k+1,1}$, and two inputs which are a positive Reeb chord $\gamma_{1,k}$ (output of $f^{(n+1)}$), and the chord $a_k$ as a negative asymptotic.
	\item if $a_k,\dots,a_{n+j+1}$ contains at least one intersection point, we get then an unfinished building with three components in the middle level: two of them have a common asymptotic to an intersection point and the other is disjoint from them and has asymptotics to $a_k,\dots,a_{n+j+1}$ and a Reeb chord $\gamma_{n+j+1,k+1}$ (output of $f^{(k-n-j)}$). The bottom level contains a banana with 3 mixed Reeb chord asymptotics which are all positive: a Reeb chord $\gamma_{1,n+j+1}$ (the output of $f^{(n+1)}$) and the chord $\gamma_{n+j+1,k+1}$ as inputs, and a chord $\gamma_{k+1,1}$ as output.
\end{enumerate}
In the two cases above, we assume that at least one asymptotic among $a_{n+j},\dots,a_{n+1},a_n,\dots a_1$ is an intersection point, otherwise the term vanishes for energy reasons. Then, in each case the two components in the middle level having a common asymptotic can be glued and thus after gluing we get an unfinished building in:
\begin{alignat}{1}
\widetilde{\cM^1}(\gamma_{k+1,1};\gamma_{1,k},a_k)\times\cM^1(\gamma_{1,k};a_1,a_2,\dots,a_{k-1})\label{cas1}
\end{alignat}
for the case (1) above, and in
\begin{alignat}{1}
\widetilde{\cM^1}(\gamma_{k+1,1};\gamma_{1,n+j+1},\gamma_{n+j+1,k+1})&\times\cM^1(\gamma_{1,n+j+1};a_1,a_2,\dots,a_{n+j})\label{cas2}\\
&\times\cM^0(\gamma_{n+j+1,k+1};a_{n+j+1},a_{n+j+2},\dots,a_k)\nonumber
\end{alignat}
for the case (2) above. More generally, each term of the sum \eqref{second terme1} takes the form
$$b^{(j)}\big(f^{(i_j)}\otimes\dots\otimes f^{(i_s)}(\id^{\otimes p}\otimes m_q^0\otimes\id^{\otimes r})\otimes\dots\otimes f^{(i_1)}\big)$$
with $p+q+r=i_s$. Hence, analogously to the two special terms described above, the unfinished buildings contributing to \eqref{second terme1}, are composed by several rigid curves in the middle level so that two of them have a common asymptotic to an intersection point, and one index-$1$ banana in the bottom level, having for positive input asymptotics the output chords of the maps $f^{(i_\alpha)}$, and potentially some negative chords among $a_1,\dots,a_k$ (which happens when $f^{(i_\alpha)}=f^{(1)}$ and the corresponding input is a Reeb chord), and finally a positive Reeb chord asymptotic $\gamma_{k+1,1}$ as output. These buildings are in the boundary of the compactification of products of moduli spaces of type \eqref{cas2}, with possibly more or no (as for \eqref{cas1}) rigid components in the middle level. In such products, there is only one moduli space of non-rigid pseudo-holomorphic disks. These disks have:
\begin{enumerate}
	\item boundary on the compact parts of the Lagrangian cobordisms that without loss of generality we label $\Sigma_1,\dots,\Sigma_{k+1}$ (in order to simplify the indices notation in the description of the boundary below),
	\item punctures asymptotic to intersection points and chords in the complexes $CF_{-\infty}(\Sigma_i,\Sigma_{i+1})$,
	\item one negative puncture asymptotic to a Reeb chord $\gamma_{1,k+1}$ (output of the map $f^{(k)}$).
\end{enumerate}
We denote by
\begin{alignat*}{1}
\cM^1_{\Sigma_{1,\dots,k+1}}(\gamma_{1,k+1};a_1,a_2,\dots,a_k)
\end{alignat*}
such non compact moduli spaces, with $a_i\in CF_{-\infty}(\Sigma_i,\Sigma_{i+1})$, and say that these moduli spaces are of \textit{type A}. The discussion above implies that in order to deduce Relation \eqref{second terme}, we will have to study the boundary of the compactification of moduli spaces of type $A$. Before that, let us describe the kinds of buildings that contribute to the terms of the sum \eqref{second terme2}.
One of the terms in \eqref{second terme2} is for example:
\begin{alignat*}{1}
b^{(1)}\circ f^{(k-1)}\big(a_k,\dots,a_3,b^{(1)}\circ f^{(2)}(a_2,a_1)\big)
\end{alignat*}
which vanishes if $a_2,a_1$ are both Reeb chords or if $a_k,\dots,a_3$ are all Reeb chords. If it does not vanish, such a composition of maps is given by a count of unfinished buildings of height $2|1|0$ such that the components in the middle level and in the bottom level on floor $-1$, corresponding to curves contributing to $f^{(k-1)}\big(a_k,\dots,a_3,b^{(1)}\circ f^{(2)}(a_2,a_1)\big)$, form together a pseudo-holomorphic building and therefore can be glued. After gluing, we get an unfinished building in
\begin{alignat}{1}
\widetilde{\cM^1}(\gamma_{k+1,1};\gamma_{1,k+1})\times\cM^1(\gamma_{1,k+1};a_1,a_2,\dots,a_k)\label{cas3}
\end{alignat}
(see Figure \ref{recoll1})
where the first moduli space is a rigid banana and the second is of type $A$. Another term appearing in \eqref{second terme2} is for example:
\begin{alignat*}{1}
b^{(2)}\big(f^{(k-2)}(a_k,\dots,a_3),f^{(1)}\circ b^{(1)}\circ f^{(2)}(a_2,a_1)\big)=b^{(2)}\big(f^{(k-2)}(a_k,\dots,a_3),b^{(1)}\circ f^{(2)}(a_2,a_1)\big)
\end{alignat*}
where the equality comes from the convention \eqref{conv1}. The unfinished pseudo-holomorphic buildings contributing to this term are of height $2|1|0$ again, but this time the components in the bottom level (on floor $-1$ and floor $-2$) form a building. After gluing the components of this building we get an unfinished building in:
\begin{alignat}{1}
\widetilde{\cM^2}(\gamma_{k+1,1};\gamma_{1,3},\gamma_{3,k+1})\times\cM^0(\gamma_{1,3};a_1,a_2)\times\cM^0(\gamma_{3,k+1};a_3,a_4,\dots,a_k)\label{cas4}
\end{alignat}
The first moduli space is a moduli space of non rigid bananas, whereas the two last moduli spaces are rigid and contribute respectively to the maps $f^{(2)}$ and $f^{(k-2)}$.

More generally, each term of \eqref{second terme2} is given by a count of unfinished buildings of height $2|1|0$ such that either
\begin{enumerate}
	\item[(A)] the components in the middle level and in the bottom level floor $-1$ form a building, or
	\item[(B)] the components in the bottom levels form a building.
\end{enumerate}
In the first case, after gluing we get an unfinished building in a product of moduli spaces of type \eqref{cas3} (with possibly several rigid curves with boundary on the non-cylindrical parts). In such a product, the only non compact moduli space is of type $A$. In the second case (which appears only when the interior $m^-$ is composed with $f^{(1)}$, because $f^{(1)}\circ m^-=m^-$ by convention) we get an unfinished building in a product of type \eqref{cas4}. In this case, the non-rigid disks are index-2 bananas: disks with boundary on the negative cylindrical ends of the cobordisms, with input punctures asymptotics to positive Reeb chords (which are outputs of maps $f^{(i)}$, $i\geq2$)) and negative Reeb chords (in the case of $f^{(1)}$), and an output puncture asymptotic to a positive Reeb chord. We will say that such moduli spaces of index-2 bananas are of \textit{type B}.

To sum up, we have seen that each term on the left-hand side of Relation \eqref{second terme} is defined by a count of unfinished pseudo-holomorphic buildings arising at the boundary of the compactification of some product of moduli spaces (as \eqref{cas1}, \eqref{cas2}, \eqref{cas3}, \eqref{cas4} and with possibly more rigid curves with boundary on the compact parts). In any case, the non compact components in such products are moduli spaces of type $A$ or of type $B$. So let us now describe the boundary of the compactification of such moduli spaces.
\vspace{2mm}

\noindent\textit{Type A:} We describe now $\partial\overline{\cM^1}_{\Sigma_{1,\dots,k+1}}(\gamma_{1,k+1};a_1,a_2,\dots,a_k)$. Remark first that this moduli space is empty if the asymptotics $a_i$ are all Reeb chords in $CF_{-\infty}(\Sigma_i,\Sigma_{i+1})$.
The pseudo-holomorphic buildings arising as limits of a one parameter family of disks of type $A$ must satisfy:
\begin{enumeratec}
	\item each curve in the building has positive energy,
	\item each curve in the building has a non negative Fredholm index,
	\item the building is asymptotic (the asymptotics that are not nodes) to $\gamma_{1,k+1},a_1,\dots,a_k$,
	\item if the building consists of the pseudo-holomorphic disks $\{u_i\}$, the relation $\sum \ind(u_i)+\nu=1$ must be satisfied. This implies that the following types of buildings can appear:
	\begin{enumeratec}
		\item buildings of height $0|1|0$ with two rigid components and boundary on the compact parts having a common asymptotic at an intersection point (pair of nodes),
		\item buildings of height $1|1|0$ with several disjoint rigid components with boundary on the compact parts, and one index-1 component with boundary on the negative ends. We can divide this case into two subcases:
		\begin{enumeratec}
			\item the output chord $\gamma_{1,k+1}$ is an asymptotic of a disk with boundary on the compact parts in the building,
			\item the chord $\gamma_{1,k+1}$ is an asymptotic of the non trivial disk with boundary on the negative ends of the cobordisms.
		\end{enumeratec}
	\end{enumeratec}
\end{enumeratec}

In the case (4)(a), the component containing the output $\gamma_{1,k+1}$ has a node asymptotic to an intersection point $q\in CF(\Sigma_{n+1},\Sigma_{n+j+1})$ for some $1\leq j\leq k$ and $0\leq n\leq k-j$. This disk is asymptotic to $\gamma_{1,k+1},a_1,\dots,a_n,q,a_{n+j+1},\dots,a_k$ in this cyclic order, and thus contributes to the map $f^{(k-j+1)}$ (with output $\gamma_{1,k+1}$). The other component is asymptotic to $q,a_{n+1},\dots,a_{n+j}$ and contributes to $\langle m^0(a_{n+j},\dots,a_{n+1}),q\rangle$ (coefficient of $q$ in $m^0(a_{n+j},\dots,a_{n+1})$). Thus in this case the building contributes to
\begin{alignat*}{1}
f^{(k-j+1)}\big(a_k,\dots,m^0_j(a_{n+j},\dots,a_{n+1}),a_n,\dots,a_1\big).
\end{alignat*}

In the case (4)(b)(i), the disk containing the output $\gamma_{1,k+1}$ is asymptotic to $$\gamma_{1,k+1},a_1,\dots,a_n,\gamma_{n+j+1,n+1},a_{n+j+1},\dots,a_k$$ with $\gamma_{n+j+1,n+1}\in C(\La_{n+1}^-,\La_{n+j+1}^-)$ negative Reeb chord which is a node corresponding thus to the positive output Reeb chord asymptotic of the disk in level $-1$. This disk in the middle level contributes to $\langle f^{(k-j+1)}(a_k,\dots,a_{n+j+1},\gamma_{n+j+1,n+1},a_n,\dots,a_1),\gamma_{1,k+1}\rangle$. Then, the disk in level $-1$ is a banana with output $\gamma_{n+j+1,n+1}$ and may have inputs at
\begin{itemize}
	\item negative Reeb chord asymptotics among $a_{n+1},\dots,a_{n+j}$,
	\item positive Reeb chord asymptotics which are nodes and are therefore negative asymptotics for (rigid) disks in the middle level having asymptotics among $a_{n+1},\dots,a_{n+j}$.
\end{itemize}
The banana together with the rigid components in the middle level not containing $\gamma_{1,k+1}$ contribute to $\langle m_j^-(a_{n+j},\dots,a_{n+1}),\gamma_{n+j+1,n+1}\rangle$. Putting these together, the buildings of case (4)(a)(i) contribute to:
\begin{alignat*}{1}
f^{(k-j+1)}\big(a_k,\dots,m_j^-(a_{n+j},\dots,a_{n+1}),a_n,\dots,a_1\big)
\end{alignat*}

In the case (4)(b)(ii), the disk containing the output $\gamma_{1,k+1}$ is an index-$1$ curve contributing to a map $\Delta$. Then, every curve in the middle level (such a curve exists because as observed before we can assume that at least one asymptotic is an intersection point otherwise the moduli spaces of type $A$ are empty) has in particular a negative Reeb chord asymptotic not in the Floer complexes (chords $\gamma_{i,i'}$ with $i<i'$) which corresponds to the output of a map $f$. These chords are nodes being also positive Reeb chords inputs for the map $\Delta$. Of course, the disk in level $-1$ can also have negative Reeb chords asymptotics among $a_1,\dots,a_k$. Using the conventions \eqref{conv1} and \eqref{conv2}, we have that the pseudo-holomorphic buildings of case (4)(b)(ii) contribute to:
\begin{alignat*}{1}
\Delta^{(s)}\big(f^{(i_s)}(a_k,...,a_{k-i_{s}+1}),\dots,f^{(i_1)}(a_{i_1},...,a_1)\big)\,\,\mbox{ with }\,\,i_s+\dots+i_1=k
\end{alignat*}

Finally, all these possibilities of breakings give the relation:
\begin{alignat}{1}
\sum\limits_{\substack{1\leq j\leq k\\ 0\leq n\leq k-j}}f^{(k-j+1)}\big(\id^{\otimes k-j-n}\otimes m^0_j\otimes\id^{\otimes n}\big)&+\sum\limits_{\substack{1\leq j\leq k-1\\ 0\leq n\leq k-j}}f^{(k-j+1)}\big(\id^{\otimes k-j-n}\otimes m^-_j\otimes\id^{\otimes n}\big)\nonumber\\
&+\sum\limits_{\substack{1\leq s\leq k\\ i_1+\dots+i_s=k}}\Delta^{(s)}\big(f^{(i_s)}\otimes\dots\otimes f^{(i_1)}\big)=0\label{brisuref}
\end{alignat}
with conventions \eqref{conv1} and \eqref{conv2}.
\vspace{2mm}

\noindent\textit{Type B: Index-$2$ bananas.} Let us assume without loss of generality that such an index-$2$ banana has boundary on the Lagrangian label $\R\times\La^-_{1,\dots,k+1}$, with output a positive Reeb chord asymptotic at $\gamma_{k+1,1}\in \Rc(\La_{k+1}^-,\La_{1}^-)$ and other Reeb chord asymptotics at Reeb chords $\gamma_{i}\in\mathfrak{C}^*(\La_i^-,\La^-_{i+1})$, for $1\leq i\leq k$. We want to describe $\partial\overline{\cM^2}(\gamma_{k+1,1};\gamma_1,\dots,\gamma_{k})$. An index-$2$ banana in $\cM^2(\gamma_{k+1,1};\gamma_1,\dots,\gamma_{k})$ is a pseudo-holomorphic disk with boundary on the negative cylindrical ends of the cobordisms, so it can break on a pseudo-holomorphic building with boundary on the negative cylindrical ends too, in particular, each (non trivial) component of the building has index at least $1$. So, an index-$2$ banana can only break into a building of height $2|0|0$ such that each floor in the bottom level contains one non trivial component, and both have a common asymptotic to a Reeb chord (node which is a positive asymptotic for the component on floor $-2$ and negative asymptotic for the component on floor $-1$). We can distinguish two cases:
\begin{enumerate}
	\item the output $\gamma_{k+1,1}$ is an asymptotic of the component on floor $-1$,
	\item the output $\gamma_{k+1,1}$ is an asymptotic of the component on floor $-2$.
\end{enumerate}

In case (1), the component on floor $-1$ is asymptotic to $$\gamma_{k+1,1},\gamma_1,\dots,\gamma_{n},\gamma_{n+j+1,n+1},\gamma_{n+j+1},\dots,\gamma_k$$ with $\gamma_{n+j+1,n+1}\in \Rc(\La^-_{n+j+1},\La_{n+1}^-)$ negative asymptotic which is a node, for some $1\leq j\leq k$ and $0\leq n\leq k-j$. Then, the disk on floor $-2$ has asymptotics to $\gamma_{n+j+1,n+1},\gamma_{n+1},\dots,\gamma_{n+j}$ with $\gamma_{n+j+1,n+1}$ positive output. So the two components of the buildings are bananas and together contribute to
\begin{alignat*}{1}
b^{(k-j+1)}\big(\gamma_k,\dots,b^{(j)}(\gamma_{n+j},\dots,\gamma_{n+1}),\gamma_n,\dots,\gamma_1\big)
\end{alignat*}

In case (2), the component on floor $-2$ is asymptotic to $$\gamma_{k+1,1},\gamma_1,\dots,\gamma_{n},\gamma_{n+1,n+j+1},\gamma_{n+j+1},\dots,\gamma_k$$ with $\gamma_{n+1,n+j+1}\in\Rc(\La^-_{n+1},\La_{n+1+j}^-)$ positive asymptotic which is a node, for some $1\leq j\leq k$ and $0\leq n\leq k-j$. Then, the disk on floor $-1$ has asymptotics to $\gamma_{n+1,n+j+1},\gamma_{n+1},\dots,\gamma_{n+j}$ with $\gamma_{n+1,n+j+1}$ negative output. Thus, this disk on floor $-1$ contributes to $$\langle \Delta(\gamma_{n+j},\dots,\gamma_{n+1}),\gamma_{n+1,n+j+1}\rangle$$ and the disk on floor $-2$ is a banana with output $\gamma_{k+1,1}$. Hence, the building in case (2) contributes to
\begin{alignat*}{1}
b^{(k-j+1)}\big(\gamma_k,\dots,\Delta^{(j)}(\gamma_{n+j},\dots,\gamma_{n+1}),\gamma_n,\dots,\gamma_1\big)
\end{alignat*}

We described above all the types of pseudo-holomorphic buildings in the boundary of the compactification of $\widetilde{\cM^2}(\gamma_{k+1,1};\gamma_1,\dots,\gamma_{k})$ and we deduce from this the relation:
\begin{alignat}{1}
&\sum\limits_{\substack{1\leq j\leq k\\ 0\leq n\leq k-j}}b^{(k-j+1)}\big(\id^{\otimes k-j-n}\otimes(b^{(j)}+\Delta^{(j)})\otimes\id^{\otimes n}\big)=0\label{brisureb}
\end{alignat}
By combining Relations \eqref{brisuref} and \eqref{brisureb}, we get the following:
\begin{alignat*}{1}
&\sum\limits_{\substack{1\leq s\leq k\\ i_1+\dots+i_s=k\\1\leq\alpha\leq s}}b^{(s)}\Big(f^{(i_s)}\otimes...\otimes f^{(i_{\alpha+1})}\otimes\sum\limits_{\substack{1\leq j\leq i_\alpha \\0\leq n\leq i_\alpha-j}}f^{(i_\alpha-j+1)}\big(\id^{\otimes i_\alpha-n-j}\otimes m^0_j\otimes\id^{\otimes n}\big)\otimes...\otimes f^{(i_1)}\Big)\\
&+\sum\limits_{\substack{1\leq s\leq k\\ i_1+\dots+i_s=k\\1\leq\alpha\leq s}}b^{(s)}\Big(f^{(i_s)}\otimes...\otimes\sum\limits_{\substack{1\leq j\leq i_\alpha-1\\ 0\leq n\leq i_\alpha-j}}f^{(i_\alpha-j+1)}\big(\id^{\otimes i_\alpha-n-j}\otimes m^-_j\otimes\id^{\otimes n}\big)\otimes...\otimes f^{(i_1)}\Big)\\
&+\sum\limits_{\substack{1\leq s\leq k\\ i_1+\dots+i_s=k\\1\leq\alpha\leq s}}b^{(s)}\Big(f^{(i_s)}\otimes...\otimes f^{(i_{\alpha+1})}\otimes\sum\limits_{\substack{1\leq j\leq i_\alpha\\ n_1+\dots+n_j=i_\alpha}}\Delta^{(j)}\big(f^{(n_j)}\otimes...\otimes f^{(n_1)}\big)\otimes...\otimes f^{(i_1)}\Big)\\
&+\sum\limits_{\substack{1\leq s\leq k\\ i_1+\dots+i_s=k}}\sum\limits_{\substack{1\leq j\leq s\\ 0\leq n\leq s-j}}b^{(s-j+1)}\Big(f^{(i_s)}\otimes...\otimes b^{(j)}\big(f^{(i_{n+j})}\otimes...\otimes f^{(i_{n+1})}\big)\otimes...\otimes f^{(i_1)}\Big)\\
&+\sum\limits_{\substack{1\leq s\leq k\\ i_1+\dots+i_s=k}}\sum\limits_{\substack{1\leq j\leq s\\ 0\leq n\leq s-j}}b^{(s-j+1)}\Big(f^{(i_s)}\otimes...\otimes \Delta^{(j)}\big(f^{(i_{n+j})}\otimes...\otimes f^{(i_{n+1})}\big)\otimes...\otimes f^{(i_1)}\Big)=0
\end{alignat*}
where the sum of the first three lines equals zero because of \eqref{brisuref}, as well as the sum of the two last lines because of \eqref{brisureb}. The sums in the third and fifth lines are equal so cancel each other (on $\Z_2$) and we get then:
\begin{alignat*}{1}
&\sum\limits_{\substack{1\leq s\leq k\\ i_1+\dots+i_s=k\\1\leq\alpha\leq s}}b^{(s)}\Big(f^{(i_s)}\otimes...\otimes\sum\limits_{\substack{1\leq j\leq i_\alpha \\0\leq n\leq i_\alpha-j}}f^{(i_\alpha-j+1)}\big(\id^{\otimes i_\alpha-n-j}\otimes m^0_j\otimes\id^{\otimes n}\big)\otimes\cdots\otimes f^{(i_1)}\Big)\\
&+\sum\limits_{\substack{1\leq s\leq k\\ i_1+\dots+i_s=k\\1\leq\alpha\leq s}}b^{(s)}\Big(f^{(i_s)}\otimes\cdots\otimes\sum\limits_{\substack{1\leq j\leq i_\alpha-1\\ 0\leq n\leq i_\alpha-j}}f^{(i_\alpha-j+1)}\big(\id^{\otimes i_\alpha-n-j}\otimes m^-_j\otimes\id^{\otimes n}\big)\otimes\cdots\otimes f^{(i_1)}\Big)\\
&+\sum\limits_{\substack{1\leq s\leq k\\ i_1+\dots+i_s=k}}\sum\limits_{\substack{1\leq j\leq s\\ 0\leq n\leq s-j}}b^{(s-j+1)}\Big(f^{(i_s)}\otimes\cdots\otimes b^{(j)}\big(f^{(i_{n+j})}\otimes\cdots\otimes f^{(i_{n+1})}\big)\otimes\cdots\otimes f^{(i_1)}\Big)=0\\
\end{alignat*}
The first sum corresponds to 
\begin{alignat*}{1}
\sum\limits_{\substack{1\leq j\leq k\\ 0\leq n\leq k-j}}m_{k-j+1}^-\big(\id^{\otimes k-j-n}\otimes m^0_j\otimes\id^{\otimes n}\big)
\end{alignat*}
and the two last sums to
\begin{alignat*}{1}
\sum\limits_{\substack{1\leq j\leq k\\ 0\leq n\leq k-j}}m_{k-j+1}^-\big(\id^{\otimes k-j-n}\otimes m^-_j\otimes\id^{\otimes n}\big)
\end{alignat*}
because the sum in the third line is equal to the missing terms 
\begin{alignat*}{1}
b^{(s-j+1)}\Big(f^{(i_s)}\otimes\cdots\otimes f^{(i_{n+j+1})}\otimes m^-_j\otimes\cdots\otimes f^{(i_1)}\Big)
\end{alignat*}
in the second line, which gives the relation. This ends the proof of Theorem \ref{teo3}.

\subsection{$A_\infty$-functor}\label{secfonc}
In this subsection, we naturally generalize the maps $\Fc^1$ and $\Fc^2$ with higher order maps. Again, let us consider a fixed $(d+1)$-tuple of pairwise transverse exact Lagrangian cobordisms $\Sigma_1,\dots,\Sigma_{d+1}$. For every $1\leq k\leq d$ and every $(k+1)$-tuple $i_1,\dots,i_{k+1}$ of distinct indices in $\{1,\dots,d+1\}$, we will construct a map 
\begin{alignat*}{1}
	\Fc^k: CF_{-\infty}(\Sigma_{i_k},\Sigma_{i_{k+1}})\otimes\dots\otimes CF_{-\infty}(\Sigma_{i_1},\Sigma_{i_2})\to C^*(\La_{i_1}^+,\La_{i_{k+1}}^+)
\end{alignat*}
such that the family of maps $\{\Fc^k\}_{1\leq k\leq d}$ satisfies the $A_\infty$-functor relation, i.e. for every $1\leq k\leq d$ and every $(k+1)$-tuple $i_1,\dots,i_{k+1}$ of distinct indices, we have:
\begin{alignat}{1}
\sum_{\substack{1\leq j\leq k\\ 0\leq n\leq k-j}}\Fc^{k-j+1}(\id^{\otimes k-j-n}\otimes\mfm_j\otimes\id^{\otimes n})+\sum_{\substack{1\leq s\leq k\\j_1+\dots+j_s=k}}\mu^s(\Fc^{j_s}\otimes\dots\otimes\Fc^{j_1})=0\label{foncteur}
\end{alignat}
where
\begin{itemize}
	\item for $1\leq j\leq k$ and $0\leq n\leq k-j$ we have
	\begin{alignat*}{1}
	\mfm_j:CF_{-\infty}(\Sigma_{i_{n+j}},\Sigma_{i_{n+j+1}})\otimes\dots\otimes CF_{-\infty}(\Sigma_{i_{n+1}},\Sigma_{i_{n+2}})\to CF_{-\infty}(\Sigma_{i_{n+1}},\Sigma_{i_{n+j+1}})
	\end{alignat*}
	and $\Fc^{k-j+1}$ has domain
	\begin{alignat*}{1}
	CF_{-\infty}(\Sigma_{i_{k}},\Sigma_{i_{k+1}})\otimes\dots\otimes CF_{-\infty}(\Sigma_{i_{n+1}},\Sigma_{i_{n+j+1}})\otimes\dots\otimes CF_{-\infty}(\Sigma_{i_{1}},\Sigma_{i_{2}})
	\end{alignat*}
	and codomain $CF_{-\infty}(\Sigma_{i_{1}},\Sigma_{i_{k+1}})$,
	\item for all $1\leq s\leq k$, all indices $j_1,\dots,j_s$ such that $j_1+\dots+j_s=k$ and all $1\leq\alpha\leq s$, the domain of $\Fc^{j_\alpha}$ is 
	\begin{alignat*}{1}
	CF_{-\infty}(\Sigma_{i_{j_1+...+j_{\alpha}}},\Sigma_{i_{j_1+...+j_{\alpha}+1}})\otimes\dots\otimes CF_{-\infty}(\Sigma_{i_{j_1+...+j_{\alpha-1}+1}},\Sigma_{i_{{j_1+...+j_{\alpha-1}+2}}})
	\end{alignat*}
	and the codomain is $CF_{-\infty}(\Sigma_{i_{j_1+...+j_{\alpha-1}+1}},\Sigma_{i_{{j_1+...+j_{\alpha}+1}}})$,
	\item $\mu^s$ are the $A_\infty$-maps of the augmentation category $\Aug_{-}(\La_1^+\cup\dots\cup\La_{k+1}^+)$. In the formula above, we have
	\begin{alignat*}{1}
		\mu^s:CF_{-\infty}(\Sigma_{i_{j_1+...+j_{s-1}+1}},\Sigma_{i_{{k+1}}})\otimes\dots\otimes CF_{-\infty}(\Sigma_{i_{1}},\Sigma_{i_{j_1+1}})\to CF_{-\infty}(\Sigma_{i_{1}},\Sigma_{i_{k+1}})
	\end{alignat*}
	and we did not write the augmentations in the index in order for the formula to stay readable.
\end{itemize} 

\begin{figure}[ht]  
	\begin{center}\includegraphics[width=5cm]{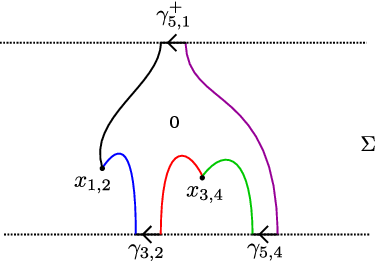}\end{center}
	\caption{Example of curve contributing to $\Fc^4(\gamma_{5,4},x_{3,4},\gamma_{3,2},x_{1,2})$.}
	\label{foncteur4}
\end{figure}

Now, to simplify notations once again, let us assume that the $(k+1)$-tuple $i_1,\dots,i_{k+1}$ of distinct indices is equal to $1,\dots,k+1$.
The maps $\Fc^k$, for $k\geq3$, are defined analogously to the maps $\Fc^1$ and $\Fc^2$, by a count of rigid pseudo-holomorphic disks, but with more mixed asymptotics (see for example Figure \ref{foncteur4}). So we define
\begin{alignat*}{1}
\Fc^k\colon CF_{-\infty}(\Sigma_k,\Sigma_{k+1})\otimes\dots\otimes CF_{-\infty}(\Sigma_1,\Sigma_2)\to C^*(\La_1^+,\La_{k+1}^+)
\end{alignat*}
by
\begin{alignat}{1}
\Fc^k(a_k,\dots,a_1)=\sum_{\substack{\gamma^+\in\Rc(\La_{k+1}^+,\La_1^+)\\ \bs{\delta}_1,\dots,\bs{\delta}_{k+1}}}\#\cM^0_{\Sigma_{1,2,\dots,k+1}}(\gamma^+;\bs{\delta}_1,a_1,\dots,\bs{\delta}_k,a_k,\bs{\delta}_{k+1})\cdot\ep^-\cdot\gamma^+\label{def_foncteur}
\end{alignat}
where as always the $\bs{\delta}_i$'s are words of Reeb chords of $\La_i^-$, and again the term $\ep^-$ should be replaced by $\prod \ep_i^-(\bs{\delta}_i)$.
In order to show the $A_\infty$-functor relation, we study degeneration of curves in the moduli space $\cM^1(\gamma^+;\bs{\delta}_1,a_1,\dots,\bs{\delta}_k,a_k,\bs{\delta}_{k+1})$. A family of curves in such a moduli space can break into:
\begin{enumerate}
	\item a pseudo-holomorphic building of height $0|1|0$ such that the middle level contains two index-$0$ curves which have a common asymptotic at an intersection point. These buildings contribute thus to $$\langle \Fc^{k-j+1}(\id^{\otimes k-j-n}\otimes m^0_j\otimes\id^{\otimes n}),\gamma^+\rangle$$
	\item a pseudo-holomorphic building of height $1|1|0$ with possibly several index-$0$ curves in the middle level and one index-$1$ banana in the bottom level (for index reasons, this level can not contain any other non trivial curve). These buildings contribute to $$\langle\Fc^{k-j+1}(\id^{\otimes k-j-n}\otimes m^-_j\otimes\id^{\otimes n}),\gamma^+\rangle$$
	\item a pseudo-holomorphic building of height $0|1|1$ with again possibly several index-$0$ disks in the middle level and one index-$1$ curve in the top level. These buildings contribute to $$\langle\mu^s(\Fc^{j_s}\otimes\cdots\otimes\Fc^{j_1}),\gamma^+\rangle$$
\end{enumerate}
As boundary of a 1-dimensional manifold, the sum of all these contributions gives $0$ modulo 2, and this implies the relation \eqref{foncteur}.

\bibliographystyle{alpha}
\bibliography{references.bib}
\end{document}